\documentclass{amsart}
\usepackage{amssymb}
\usepackage{amsmath}
\usepackage{mathdots}
\usepackage[dvipdfm,colorlinks,backref,pagebackref,hypertexnames=false,linkcolor=blue,urlcolor=green,citecolor=red]{hyperref}
\newcommand{\FF}{{\mathcal{F}}}

\newcommand{\N}{{\mathbb{N}}}
\newcommand{\Z}{{\mathbb{Z}}}
\newcommand{\F}{{\mathbb{F}}}
\newcommand{\C}{{\mathbb{C}}}
\newcommand{\T}{{\mathbb{T}}}

\newcommand{\la}{\langle}
\newcommand{\ra}{\rangle}

\newcommand{\al}{\alpha}
\newcommand{\lam}{\lambda}
\newcommand{\Sig}{\Sigma}

\newcommand{\dfour}[4]{\begin{footnotesize}\ensuremath{\hspace{-0.1cm}\begin{array}{ccc}&#2&\\#1\hspace{-0.2cm}&#3&\hspace{-0.2cm}#4\end{array}\hspace{-0.1cm}}\end{footnotesize}}
\newcommand{\Esix}[6]{\begin{footnotesize}\ensuremath{\hspace{-0.1cm}\begin{array}{ccccc}&&#2&&\\#1\hspace{-0.2cm}&#3\hspace{-0.2cm}&#4&\hspace{-0.2cm}#5&\hspace{-0.2cm}#6\end{array}\hspace{-0.1cm}}\end{footnotesize}}
\newcommand{\Eeight}[8]{\begin{footnotesize} \ensuremath{\hspace{-0.1cm} \begin{array}{ccccccc}&&#2&&&&\\#1\hspace{-0.2cm}&#3\hspace{-0.2cm}&#4&\hspace{-0.2cm}#5&\hspace{-0.2cm}#6&\hspace{-0.2cm}#7&\hspace{-0.2cm}#8\end{array}\hspace{-0.1cm}}\end{footnotesize}}

\renewcommand{\epsilon}{\varepsilon}

\newtheorem{theorem}{Theorem}[section]
\newtheorem{lemma}[theorem]{Lemma}
\newtheorem{proposition}[theorem]{Proposition}

\newtheorem{definition}[theorem]{Definition}

\begin{document}

\title[On the Character Degrees of Sylow $p$-subgroups of $E(p^f)$]{On the Character Degrees of Sylow $p$-subgroups of Chevalley Group of Type $E(p^f)$}

\author{Tung Le and Kay Magaard}

\address{T.L.: Institute of Mathematics, University of Aberdeen, Aberdeen AB24 3UE, U.K.}

\address{K.M.: School of Mathematics, University of Birmingham, Edgbaston, Birmingham B15 2TT, U.K.}

\email{T.L: t.le@abdn.ac.uk}

\email{K.M.: k.magaard@bham.ac.uk}

\date{August 12, 2010}



%%%%%%%%%%%%%%%%%%%%%%%%%%%%%%%%%%%%%%%%%%%%%%%%%%%%%%%%%%%%%%%%%%%%%%%%%%%%%%%%%%%%%%%%%%

\begin{abstract}

Let $\F_q$ be a field of characteristic $p$ with $q$ elements. It is known that the degrees of the irreducible characters of the Sylow $p$-subgroup of $GL_n(\F_q)$ are powers of $q,$ see Isaacs \cite{Is2}. On the other hand Sangroniz \cite{San} showed that this is true for a Sylow
$p$-subgroup of a classical group defined over $\F_q$ if and only if $p$ is odd. For the classical groups of Lie type
$B$, $C$ and $D$ the only bad prime is $2$. For the exceptional groups there are others.
In this paper we construct irreducible characters for the Sylow $p$-subgroups of the Chevalley groups $D_4(q)$ with $q=2^f$
of degree $q^3/2$. Then we use an analogous construction for $E_6(q)$ with $q=3^f$ to obtain characters of degree $q^7/3$, and for $E_8(q)$ with $q=5^f$ to obtain characters of degree $q^{16}/5.$ This helps
to explain why the primes $2,$ $3$ and $5$ are bad for the Chevalley groups of type $E$ in terms of the representation theory of
the Sylow $p$-subgroup.
\end{abstract}

\maketitle

\section{Introduction}

Let $G$ be a Chevalley group defined over a field $\F_q$ of order $q$ and characteristic $p > 0$.
By $\alpha_0$ we denote the highest root of the root system $\Sigma$ of $G$. It is well known that $\alpha_0$
is a positive integral linear combination of the fundamental roots of $\Sigma$. So without loss $\alpha_0 = \sum_{i = 1}^r a_i\alpha_i$
where the $\alpha_i$ are fundamental roots of $\Sigma$. Recall that $p$ is a bad prime for $G$ if $p$ is a divisor of some $a_i$.

It is well known that if $G$ classical then the only possible bad prime for $G$ is $2$. On the other hand if $G$ is exceptional of type $E$, then the prime $3$ is also bad. The ``badness'' of the prime evidences itself in the classification of the unipotent conjugacy classes of $G$. Here we aim to explain why the primes $3$ and $5$ are bad for groups of type $E$ in terms to the representation theory of the Sylow $p$-subgroup of $G = E_6(q)$ with prime 3 and $G=E_8(q)$ with prime 5. Let $UE_k(q)$ denote the unipotent radical of the standard Borel subgroup of $E_k(q)$ for $k=6$ and $8$;
i.e. the subgroup generated by all the positive root groups of $G$.
By $U_k$ we denote the quotient $UE_k(q)/K_{k-1}$, where $K_{k-1}$ is the normal subgroup of $UE_k(q)$ generated by all root groups  $X_\alpha$
such that $\alpha$ has height $k-1$ or more. Clearly any character of $U_k$ inflates to a character of $UE_k(q)$.
Abusing terminology slightly we call the image under the natural projection of a root group of $UE_k(q)$, a root
group of $U_k$. We observe that $Z(U_k)$ is generated by the root groups of height $k-2$ and
hence $|Z(U_k)| = q^{k-1}$. We define the family

\begin{center} ${\mathcal F_k} := \{ \chi \in {\rm Irr}(U_k) \ : \ X_\alpha \not \subset {\rm Ker}(\chi) \mbox{ for all } X_\alpha \subset Z(U_k) \}.$ \end{center}

%%%%%%%%%%%%%%%%%%%%%%%%%%%%%%%%%%% MAIN INTRO THEOREM I %%%%%%%%%%%%%%%%%%%%%%%%%%%%%%%

\begin{theorem} \label{MainIntro}
The following are true.
\begin{itemize}

\item[(a)] If $q = 3^f$, then for all $ \chi \in {\mathcal F_6}$ we have $ \chi(1) \in \{ q^7, q^7/3 \}$. Moreover
${\mathcal F_6}$ contains exactly $(q-1)^5(q^2 - (q-1)/2)$ characters of degree $q^7$ and exactly $3^2(q-1)^6/2$ characters of degree $q^7/3$.

\item[(b)] If $q = 5^f$, then for all $ \chi \in {\mathcal F_8}$ we have $ \chi(1) \in \{ q^{16}, q^{16}/5 \}$. Moreover
${\mathcal F_8}$ contains exactly $(q-1)^8(q^3+q^2+q+3/4)$ characters of degree $q^{16}$ and exactly $5^2(q-1)^8/4$ characters of degree $q^{16}/5$.

\end{itemize}
\end{theorem}

%%%%%%%%%%%%%%%%%%%%%%%%%%%%%%%%%%% END OF MAIN INTRO THEOREM I %%%%%%%%%%%%%%%%%%%%%%%%%%%%%%%

We remark that $9(q-1)^6/2,$ $(q-1)^5(q^2 - (q-1)/2),$ $(q-1)^8(q^3+q^2+q+3/4)$ and $25(q-1)^8/4$  are not in $\Z[q].$ On the other hand we remark also that $|{\mathcal F_6}| = (q-1)^5q^2 \in \Z[q]$ and every character in ${\mathcal F}_6$ has degree $q^7$ whenever $p \neq 3,$ and that $|\FF_8|=(q-1)^7q^{4}\in \Z[q]$ and every character in $\FF_8$ has degree $q^{16}$ whenever $p\neq 5.$  Taken together these remarks
provide evidence for a generalization of Higman's conjecture for groups of type $UE_i(q)$, $i = 6,7,8$, see for example \cite{Hig},
namely that $|{\rm Irr}(UE_i(q))| \not \in \Z[q]$ if and only if $p$ is a bad prime for $E_i(q)$ .

To prove our main theorem we begin by analyzing our construction of the irreducible characters of the Sylow $2$-subgroup
of $D_4(2^f)$ from \cite{FTK}. Our starting point is the quotient of $UD_4(q)/K_4$ where $UD_4(q)$ is the unipotent
radical of the standard Borel subgroup of the universal Chevalley group $D_4(q)$ and $K_4$ is the normal
subgroup of $UD_4(q)$ generated by the root groups of roots of height $4$ and $5$. We showed that when $p=2$ we have
a family of characters of degree $q^3/2$ of size $4(q-1)^4$. As $UD_4(q)$ is a quotient of $UE_i(q)$ for $i =6,7,8$ we
also have families of irreducible characters of degree $q^3/2$ for groups of type $UE_i(q)$, $i =6,7,8$ and $q$ is even.
%We remark here that for $p=5$  we expect to show that $UE_8(q)/K_6$, where $K_6$ is the normal subgroup of $UE_8(q)$ generated by root subgroups of roots of height $6$ or more, has a family of characters of degree $q^{16}/5$ in a subsequent paper.

Our construction is fairly elementary. Starting with large elementary abelian normal subgroups we
construct our characters via induction, using Clifford theory. To compute the necessary stabilizers we critically
use Proposition \ref{prop:mainFp} and Lemma \ref{lem:Reduction}. Throughout this paper we fix a nontrivial homomorphism $\phi:(\F_q,+)\longrightarrow \C^\times.$ For each $a\in\F_q,$ we define $\phi_a(x):=\phi(ax)$ for all $x\in\F_q,$ and denote $\F_q^\times:=\F_q-\{1\}.$ Hence, $\{\phi_a:a\in \F_q^\times\}$ are all non-principal irreducible characters of $\F_q.$

\begin{definition}
For $a\in \F_q,$ we define $\T_a :=\{t^p-a^{p-1}t: t\in\F_q\}.$
\end{definition}
\noindent We note that $\T_0=\F_q.$

%%%%%%%%%%%%%%%%%%%%%%%%%%%%%%%%% FINITE FIELDS' PROPERTY %%%%%%%%%%%%%%%%%%%%%%%%%%%

\begin{proposition} \label{prop:mainFp}
The following are true.
\begin{itemize}

\item[(a)] $t^p-a^{p-1}t=\prod_{c\in \F_p}(t-ca).$

\item[(b)] If $a\in \F_q^\times$, then $\T_a$ is an additive subgroup of $\F_q$ of index $p$.

\item[(c)] For each $a\in\F_q^\times,$ there exists $b\in\F_q^\times$ such that $b\T_a=ker\,\phi.$ Furthermore, $cb\T_a=ker(\phi)$ iff $c\in\F_p^\times.$

\item[(d)] $\{ \T_a:a\in\F_q^\times\}=\{ker\, \phi_a:a\in \F_q^\times\}$ are all subgroups  of index $p$ in $\F_q.$
\end{itemize}

\end{proposition}
%
%
%%%%%%%%%%%%%%%%%%%%%%%%%%%%%%%%%%%%%%%% END OF FINITE FIELDS' PROPERTY %%%%%%%%%%%%%%%%%%%%%%%%%
%
%
%
\noindent{\it Proof.} Part (a) is clear since the degree of the polynomial $t^p -a^{p-1}t$ is $p$ and
the $\F_p$-multiples of $a$ are clearly zeros. As $\F_q$ is of characteristic $p$, the map $\psi_a: \F_q \longrightarrow \F_p$ defined by $\psi_a(t) = t^p - a^{p-1}t$ is $\F_p$-linear. By Part (a) the kernel of the map is $1$-dimensional and
thus (b) follows.  Evidently (d) follows from (c). We defer the proof (c) to Subsection \ref{proof-prop:mainFp}.~$\Box$

%\smallskip
%%It is noted that for some $r\in \F_q,$ if the equation $t^p-rt=0$ has a nontrivial solution $t_0$ then $ct_0$ is another solution for all $c\in \F_p.$ Therefore, by Proposition \ref{prop:mainFp} (a), $t^p-rt=0$ has $p$ distinct solutions iff $r\in \{a^{p-1}:a\in \F_q^\times\}.$ Since $\F_q^\times$ is cylic and $\gcd(p-1,q-1)=p-1,$ we have the number ${|\{a^{p-1}:a\in \F_q^\times\}|=\frac{q-1}{p-1}}.$

%Since the element $b$ in Proposition \ref{prop:mainFp} (c) is uniquely determined by $a$ and $\phi$ up to a scalar of $\F_p^\times,$ we define.

% The element $b$ in Proposition \ref{prop:mainFp} (c) is uniquely determined by $a$ and $\phi$ up to a scalar of %$\F_p^\times,$ and thus we define.

\begin{definition}
For each $a\in\F_q^\times,$ we pick $a_{\phi}$ such that  $a_{\phi}\T_a=ker\,\phi.$
\end{definition}

By Proposition \ref{prop:mainFp} (c), $a_{\phi}$ exists and but is only determined up to a scalar in the prime field.
In the definition above we make some arbitrary choice which will not change throughout the paper.

% We use Proposition \ref{prop:mainFp} to construct characters of $U$ of degree not a power of $q.$

Throughout we fix notation as follows. Let $G$ be a group.  $G^\times:=G-\{1\},$ $Irr(G)$ the set of all complex irreducible characters of $G,$ and $Irr(G)^\times:=Irr(G)-\{1_G\}.$
For $H, K\leq G,$ and $\xi\in Irr(H),$ define $Irr(G/K):=\{\chi \in Irr(G):K\subset ker(\chi)\},$ %denote the set of all irreducible characters of $G$ with $K$ in the kernel.
 $Irr(G,\xi):=\{\chi\in Irr(G):(\chi,\xi^G)\neq 0\},$ and
$Irr(G/K,\xi):=Irr(G/K)\cap Irr(G,\xi).$
%$G=H\rtimes K,$ then for each character $\xi$ of $K,$ we denote the inflation of $\xi$ to $G$ by $\xi_G.$
Furthermore, for a character $\chi$ of $G,$ we denote its restriction to $H$ by $\chi|_H.$

%%%%%%%%%%%%%%%%%%%%%%%%%%%%%%%% REDUCTION LEMMA %%%%%%%%%%%%%%%%%%%%

\begin{lemma}
\label{lem:Reduction}
Let $N \unlhd G$ and $1\in X$ be a transversal of $N$ in $G$. Suppose $N=ZYM$ where $Y\unlhd N,$
%%is abelian
$Z\subset Z(N)$, $M\leq N$ and $X\subset N_G(ZY).$ If there is $\lambda\in Irr(ZY)$ such that
$Y\subset ker(\lambda),$                                    %% together with $Y$ normal in $H$ giving the weldefined for the quotient $H/Y$
%% $\{{}^x\lambda|_Y: x\in X\}=Irr(Y),$                        %% injective for induction map and that's why need $Y$ abelian
and ${}^u\lambda\neq {}^v\lambda$ for all $u\neq v\in X,$   %% irreducible induction and injective
then the following are true.
\begin{itemize}
\item[(a)] For all $\chi\in Irr(N/Y,\lambda),$ $\chi^G\in Irr(G).$ Moreover, if $\chi_1\neq \chi_2\in Irr(N/Y,\lambda),$ then ${\chi_1}^G\neq {\chi_2}^G.$

\item[(b)] The induction map from $Irr(N/Y,\lambda)$ to $Irr(G,\lambda)$ is bijective.

\end{itemize}
\end{lemma}
\noindent{\it Proof.} See Subsection \ref{proof-lem:Reduction}.~$\Box$

%%%%%%%%%%%%%%%%%%%%%%%%%%%%%% END OF REDUCTION LEMMA %%%%%%%%%%%%%%%%%%%%%%%%%%

\medskip

We recall that a $p$-group $P$ is monomial, i.e. for each $\chi\in Irr(P),$ there exist a subgroup $H$ of $P$ and a linear character $\lambda$ of $H$ such that $\chi=\lam^P.$ To construct irreducible characters whose degrees are not powers of $q =p^f, f>1$ we construct subgroups $H\unlhd P$ and $T \leq P$ such that $T$ is a transversal of $H$. Then we find a linear character $\lambda$ of $H$ such that the order of the stabilizer $Stab_T(\lambda)$ of $T$ is not a power of $q$. Moreover we insure that $\lambda$ is extendable to the inertial group $I_P(\lambda)=HStab_T(\lambda).$
Let $\lambda_I$ denote some extension of $\lambda$ to $I_P(\lambda)$. By Clifford theory the induction of $\lambda_I$ to $P$ is irreducible, of degree not a power of $q$. The existence of a suitable pair $(H,\lambda)$ is based on Proposition \ref{prop:mainFp} because a polynomial of the form $x^p+a^{p-1}x,$ $a\neq 0,$ appears in the formulae of the action of elements of $T$ on the characters of $H$.

We will now highlight the main steps of the constructions of our characters. We have deferred all of our proofs to Section \ref{Section-All-proofs}.

%***************************************************************************************
%***************************************************************************************

\section{Sylow $2$-subgroups of the Chevalley groups $D_4(2^f)$} \label{Section-D4}

Let $\F_q$ be a field of order $q$ and characteristic 2. Let $\Sig:=\la \alpha_1,\alpha_2,\alpha_3,\alpha_4\ra$ be the root system
of type $D_4,$ see Carter \cite{cart2}, Chapter 3. The Dynkin diagram of $\Sig$ is

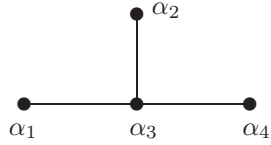
\begin{center}
\setlength{\unitlength}{1cm}
\begin{picture}(4,2)
\thinlines
\put(0.5,0.5){\circle*{0.17}}
\put(2,0.5){\circle*{0.17}}
\put(3.5,0.5){\circle*{0.17}}
\put(2,1.7){\circle*{0.17}}
\put(0.5,0.5){\line( 1, 0){1.5}}
\put(2,0.5){\line( 1, 0){1.5}}
\put(2,1.7){\line( 0, -1){1.2}}
\put(0.3,0.1){$\alpha_1$}
\put(1.9,0.1){$\alpha_3$}
\put(3.4,0.1){$\alpha_4$}
\put(2.2,1.7){$\alpha_2$}
\end{picture}
\end{center}

The positive roots are those roots which can be written as linear
combinations of the simple roots $\alpha_1, \alpha_2, \alpha_3,
\alpha_4$ with nonnegative coefficients and we write $\Sig^+$ for the
set of positive roots. We use the notation \dfour1121 for the
root $\alpha_1 + \alpha_2 + 2 \alpha_3 + \alpha_4$ and we use a
similar notation for the remaining positive roots. The $12$ positive
roots of $\Sigma$ are given in Table~\ref{tab:posrootsD4}.

\begin{table}[!ht]
\caption{Positive roots of the root system $\Sigma$ of type $D_4.$}
\label{tab:posrootsD4}

\begin{center}
\begin{tabular}{c|llll}
\hline
\rule{0cm}{0.4cm}
Height & Roots &&&
\rule[-0.1cm]{0cm}{0.4cm}\\
\hline% \cline{1-3} \hline
\rule{0cm}{0.5cm}
5 & $\alpha_{12} := $ \dfour1121 &&&
\rule[-0.2cm]{0cm}{0.4cm}\\
\hline
\rule{0cm}{0.4cm}
4 & $\alpha_{11} := $ \dfour1111 &&&
\rule[-0.2cm]{0cm}{0.4cm}\\
\hline
\rule{0cm}{0.4cm}
3 & $\alpha_8 := $ \dfour1110 & $\alpha_9 := $ \dfour1011 & $\alpha_{10} := $ \dfour0111 &
\rule[-0.2cm]{0cm}{0.4cm}\\
\hline
\rule{0cm}{0.4cm}
2 & $\alpha_5 := $ \dfour1010 & $\alpha_6 := $ \dfour0110 & $\alpha_7 := $ \dfour0011 &
\rule[-0.2cm]{0cm}{0.4cm}\\
\hline
\rule{0cm}{0.4cm}
1 & $\alpha_1$  & $\alpha_2$  & $\alpha_3$   & $\alpha_4$
\rule[-0.2cm]{0cm}{0.4cm}\\
\hline
\end{tabular}
\end{center}
\end{table}

For $\alpha \in \Sigma$ we denote the corresponding root subgroup of the Chevalley group $G$ by $X_\alpha$ whose
elements we label by $x_\alpha(t)$ where $t \in \F_q$. We note that $X_\alpha \cong  (\F_q, +)$.

We recall that the commutator formula $[x_{\alpha}(r), x_{\beta}(s)] = x_{\alpha+\beta}(-C_{\alpha,\beta}rs)$ if $\alpha +\beta \in\Sig,$ and $= 1$ otherwise, see Carter \cite{cart2}, Theorem~5.2.2. Since $p=2$ we have $-1=1$ in $\F_q,$ so all
non-zero coefficients $C_{\alpha,\beta}$ are equal to $1$. For positive roots, we use the abbreviation
$x_i(t) := x_{\alpha_i}(t),$ $i=1,2,\dots,12.$
All nontrivial commutators are given in Table~\ref{tab:commrelD4}.

\begin{table}[!ht]
\caption{Commutator relations for type $D_4.$}
\label{tab:commrelD4}

%\begin{center}
\begin{tabular}{ll}
$\left[x_1(t), x_3(u)\right]=x_5(tu),$ &
$\left[x_1(t), x_6(u)\right]=x_8(tu),$ \\
$\left[x_1(t), x_7(u)\right]=x_9(tu),$ &
$\left[x_1(t), x_{10}(u)\right]=x_{11}(tu),$ \\
$\left[x_2(t), x_3(u)\right]=x_6(tu),$  &
$\left[x_2(t), x_5(u)\right]=x_8(tu),$\\
$\left[x_2(t), x_7(u)\right]=x_{10}(tu),$ &
$\left[x_2(t), x_9(u)\right]=x_{11}(tu),$ \\
$\left[x_3(t), x_4(u)\right]=x_7(tu),$ &
$\left[x_3(t), x_{11}(u)\right]=x_{12}(tu),$ \\
$\left[x_4(t), x_5(u)\right]=x_9(tu),$  &
$\left[x_4(t), x_6(u)\right]=x_{10}(tu),$\\
$\left[x_4(t), x_8(u)\right]=x_{11}(tu),$  &
$\left[x_5(t), x_{10}(u)\right]=x_{12}(tu),$\\
$\left[x_6(t), x_9(u)\right]=x_{12}(tu),$ &
$\left[x_7(t), x_8(u)\right]=x_{12}(tu)$
\end{tabular}
%\end{center}
\end{table}

The group $UD_4$ generated by all $X_\alpha$ for $\al\in\Sig^+$ is a Sylow $2$-subgroup of the Chevalley group $D_4(q).$ Each element
$u \in UD_4$ can be written uniquely as

\begin{center}
$u = x_1(t_1) x_2(t_2) x_4(t_4) x_3(t_3) x_5(t_5) \cdots x_{12}(t_{12})$
where $x_i(t_i)\in X_i.$
\end{center}
So we write $\prod_{i=1}^{12}x_i(t_i)$ as this order. We note that our ordering of the roots is slightly non-standard as the positions
of $x_3$ and $x_4$ are reversed.

%%say some about triality if we need
We define $\FF_4 := \{ \chi \in Irr(UD_4(q)) \ : \  \chi|_{X_i}=\chi(1)\phi_{a_i} \mbox{ for each } a_8,a_9,a_{10}\in\F_q^\times \}.$
If $\Psi$ is a representation affording $\chi \in \FF_4$, then we have $\Psi([x_8(t_8),x_4(t_4)])=[\Psi(x_8(t_8)),\Psi(x_4(t_4))]=[\phi_{a_8}(t_8)\Psi(1),\Psi(x_4(t_4))]=\Psi(1)$ for all $t_4,t_8\in\F_q.$ Therefore, $X_{11}=[X_8,X_4]\subset ker(\chi).$ Use the same argument for $X_{12}=[X_8,X_7]\subset ker(\chi).$ Thus only the factor group $U = UD_4/X_{12}X_{11}$ acts on a module affording $\chi$.  Therefore, we may work with $U$ which has has order $q^{10},$ and $Z(U)=X_8X_9X_{10}.$

\begin{center}
\begin{picture}(125,105)
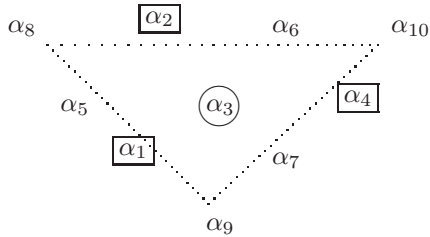
\put(0,0){Figure $UD_4(q)$: Relations of Roots}
\put(0,90){$\alpha_8$}
\put(50,90){\framebox(15,10){$\alpha_2$}}
\put(100,90){$\alpha_6$}
\put(145,90){$\alpha_{10}$}
\put(20,60){$\alpha_5$}
\put(40,40){\framebox(15,10){$\alpha_1$}}
\put(75,15){$\alpha_9$}
\put(100,40){$\alpha_7$}
\put(125,60){\framebox(15,10){$\alpha_4$}}
\put(80,62){\circle{15}}
\put(75,60){$\alpha_3$}
\qbezier[30](15,85)(95,85)(140,85)
\qbezier[30](140,85)(108,55)(76,25)
\qbezier[30](76,25)(46,55)(15,85)
\end{picture}
\end{center}

\smallskip
Let $H:=[U,U]=X_5X_6X_7X_8X_9X_{10},$ and $T=X_1X_2X_4.$ It is clear that $H,$ $HX_3$ and $T$ are elementary abelian. The group $U$ can be visualized in the above figure. The roots in boxes are in $T,$ the others outside are in $H,$ and $\alpha_3,$ which is neither in $T$ nor in $H,$ is in a circle. The disconnected lines demonstrate the relations of sum of roots equal to roots in center, e.g. $\alpha_2+\alpha_5=\alpha_6+\alpha_1=\alpha_8,...$
Each edge of the triangle contains four roots, two on the outside and two
on the inside. Each of the three vertices of the triangle together
with the four adjacent inside roots forms a hook of length 2, see \cite{FTK}.
The corresponding hook group is special of order $q^{1+4}$. The group generated by each hook group $q^{1+4}$ and $X_3$ is isomorphic to a Sylow $2$-subgroup of the general linear group $GL_4(\F_q).$

%It also says about the hook's structure which is defined in \cite{FTK}. Therefore, the figure gives a triangle.

To study those characters $\chi,$ we start with a linear character $\lambda$ of $H$ such that $\lambda|_{X_i}\neq 1_{X_i}$ for $i=8,9,10.$

\begin{definition}
For  $a_8,a_9,a_{10}\in\F_q^\times$ and $b_5,b_6,b_7\in \F_q,$  we define
\begin{itemize}

\item[(a)] $\lam^{a_8,a_9,a_{10}}_{b_5,b_6,b_7}(\prod_{i=5}^{10}x_i(t_i)):=\phi(\sum_{i=5}^7b_it_i+\sum_{j=8}^{10}a_jt_j).$

\item[(b)] $S_{567}:=\{x_{567}(t):=x_5(a_{10}t)x_6(a_9t)x_7(a_8t):t\in \F_q\}.$

\item[(c)] $S_{124}:=\{x_{124}(t):=x_1(a_{10}t)x_2(a_9t)x_4(a_8t):t\in\F_q\}.$

\item[(d)] $A:=a_8a_9a_{10}$ and $t_0:=\frac{1}{A}(b_5a_{10}+b_6a_9+b_7a_8).$

\item[(e)] $F_{124}:=\{1,x_{124}(t_0)\}.$

\item[(f)] $F_3:=\{1\}$ if $t_0=0,$  and $F_3:=\{1,x_3(\frac{(t_0)_\phi}{A})\}$ otherwise.
\end{itemize}

\end{definition}

It is easy to check that $S_{567},$ $S_{124},$ $F_{124},$ $F_3$ are subgroups of $U.$ If $t_0=0,$ then $F_{124}=F_3=\{1\},$ otherwise $F_{124}\cong F_3\cong (\F_2,+).$ Since $S_{124},S_{567}\cong (\F_q,+),$ their linear characters are in the form $\phi_{b_i}(x_i(t))=\phi(b_it)$ where $i\in\{124,567\}$ for all $b_i,t\in \F_q.$ For each $\xi\in Irr(F_{124}),$ $\xi=\phi_{b_{124}}|_{F_{124}}$ for some $\phi_{b_{124}}\in Irr(S_{124}),$ $b_{124}\in\F_q.$  If $F_{124}$ is nontrivial, we choose $b_{124}\in \{0,a_{124}\}\cong (\F_2,+)$ where $\phi(a_{124}t_0)=-1.$ The same for $F_3\leq X_3,$ for each $\xi\in Irr(F_3),$ $\xi=\phi_{b_3}|_{F_3}$ for some $\phi_{b_3}\in Irr(X_3)$ and $b_3\in \{0,a_3\}\cong (\F_2,+)$ such that $\phi(a_3\frac{(t_0)_\phi}{A})=-1$ if $(t_0)_\phi$ exists.

%It is clear that $\lam^{a_8,a_9,a_{10}}_{b_5,b_6,b_7}\in Irr(H).$
For each $a_8,a_9,a_{10}\in\F_q^\times,$ there are $q^3$ linear characters $\lam^{a_8,a_9,a_{10}}_{-,-,-}$ of $H.$ By the definition of $t_0,$ there are $q^2$ of them such that $t_0=0$ and the others $q^2(q-1)$ linears such that $t_0\neq 0.$ Therefore, there are $q^2$ cases where $F_{124},F_3$ are trivial and $q^2(q-1)$ cases where $F_{124},F_3$ are of order $2.$

For all $x_1(t_1)x_2(t_2)x_4(t_4)\in T,$ we have
\begin{center}
${}^{x_1(t_1)x_2(t_2)x_4(t_4)}(\lam^{a_8,a_9,a_{10}}_{b_5,b_6,b_7})=\lam^{a_8,a_9,a_{10}}_{b_5+a_8t_2+a_9t_4,b_6+a_8t_1+a_{10}t_4,b_7+a_9t_1+a_{10}t_2}.$
\end{center}
Hence, $T$ acts on the set of linears $\{\lambda^{a_8,a_9,a_{10}}_{-,-,-}\}.$ It is easy to check that $t_0$ is invariant under this action.
All properties of $\lam^{a_8,a_9,a_{10}}_{b_5,b_6,b_7}$ are known as follows.

\begin{lemma} \label{lem:2extensions}
Set $\lambda:=\lambda_{b_5,b_6,b_7}^{a_8,a_9,a_{10}}.$ The following are true.
\begin{itemize}
\item[(a)] $S_{124}=Stab_{T}(\lambda)$ and  $S_{567}=\{x\in X_5X_6X_7: |\lambda^U(x)|=\lambda^U(1) \}.$ Moreover, $\lambda^U|_{S_{567}}=\lambda^U(1)\phi_{At_0}.$

\item[(b)] $\lambda$ extends to $HX_3F_{124}$ and $HF_3S_{124}.$ Let $\lambda_1,\lambda_2$  be extensions of $\lambda$ to $HX_3F_{124}.$ The inertia groups $I_U(\lambda_1)=HX_3F_{124}.$

\item[(c)] ${\lambda_1}^U={\lambda_2}^U\in Irr(U)$ iff $\lambda_1|_{F_3}=\lambda_2|_{F_3}$ and  $\lambda_1|_{F_{124}}=\lambda_2|_{F_{124}}.$

\end{itemize}

\end{lemma}
\noindent{\it Proof.} See Subsection \ref{proof-lem:2extensions}.

\medskip
{\bf Remark} When $q$ is odd, both $\{x\in X_5X_6X_7: |\lambda^U(x)|=\lambda^U(1)=q^4\}$ and $Stab_{T}(\lambda)$ are trivial. Thus, $\lambda$ extends to $HX_3$ and induces irreducibly to $U$ of degree $q^3.$

When $t_0\neq 0,$ the statement in Lemma \ref{lem:2extensions} (c) makes sense since the dihedral subgroup $\la F_{124},F_3\ra\subset I_U(\lambda_1).$ By Lemma \ref{lem:2extensions} (b), $X_3,S_{124}\subset I_U(\lambda)$ but $\lambda$ is not able to extend to $HX_3S_{124}$ since $[X_3,S_{124}]\nsubseteq ker(\lambda).$ %It is also remarked that if $q$ is odd, $Stab_{X_1X_2X_4}(\lambda)$ and $\{x\in X_5X_6X_7:|\lambda^U(x)|=\lambda^U(1)\}$ are trivial.

\medskip
By Lemma \ref{lem:2extensions} (a), $T$ acts on the set of $q^3$ linears $\lambda^{a_8,a_9,a_{10}}_{-,-,-}$ into $q$ orbits, each has size $q^2.$ By Lemma \ref{lem:2extensions} (b), all $q^3$ linears $\lambda^{a_8,a_9,a_{10}}_{-,-,-}$ extend to $HX_3$ and  we obtain $q^4$ linear extensions, in which there are $q^3$ linears with $t_0=0$ and  $q^3(q-1)$ linears with $t_0\neq 0.$

% tell more about irreducible constituents of $\lambda^U$ here, then restate theorem in a better form of using $\xi_{124}$ and $\xi_3$ as restriction of $\eta_i$ to subgroups $F_{124}$ and $F_3.$

If $t_0=0,$ $F_{124}$ is trivial. By Lemma \ref{lem:2extensions} (b), $\lambda$ extends to $I_U(\lambda)=HX_3\unlhd U,$ as $\eta.$ $T$ is a transversal of $HX_3$ in $U$ and acts regularly on these $q^3$ linears $\eta$ with $t_0=0.$ Therefore, $\eta^U\in Irr(U)$ of degree $q^3$ only depends on $a_8,a_9,a_{10},$ so we denote it by $\chi^{a_8,a_9,a_{10}}_{8,9,10,q^3}\in Irr(U).$ %It is also noted that
This character is the unique $\chi\in \FF_4$ of degree $q^3$  such that $\chi|_{X_i}=\chi(1)\phi_{a_i}$ where $i=8,9,10.$ Furthermore, by Lemma \ref{lem:2extensions} (a), this is the unique constituent $\chi$ of $(\lambda|_{X_8X_9X_{10}})^U$ such that $S_{567}\subset ker(\chi).$

If $t_0\neq 0,$ $F_{124},F_3$ are isomorphic to $\F_2.$ By Lemma \ref{lem:2extensions} (b), $\lambda$ extends to $HX_3F_{124}$ as $\lambda_1,$ and ${\lambda_1}^U\in Irr(U)$ of degree $\frac{q^3}{2}.$ For each $t_0\neq 0,$ by Lemma \ref{lem:2extensions} (c), all constituents ${\lambda_1}^U$ of $\lambda^U$ only depend on the restrictions of $\lambda_1$ to $F_{124}$ and $F_3.$
Therefore, we denote these constituents of $\lambda^U$ by $\chi_{8,9,10,\frac{q^3}{2}}^{b_{124},b_3,t_0,a_8,a_9,a_{10}}$ where $b_{124},b_3\in \F_2,$ $t_0,a_8,a_9,a_{10}\in\F_q^\times.$
For each $a_8,a_9,a_{10}\in \F_q^\times,$ there are $4(q-1)$ characters $\chi\in \FF_4$ of degree $\frac{q^3}{2}$ such that $\chi|_{X_i}=\chi(1)\phi_{a_i}$ where $i=8,9,10.$

The next theorem lists generic character values of all $\chi\in Irr(U)$ such that $\chi|_{X_i}=\chi(1)\phi_{a_i}$ where $i=8,9,10.$

\begin{theorem} \label{thm:Main2Chars}

For $a_8,a_9,a_{10}\in\F_q^\times,$ suppose $\chi\in Irr(U)$ such that $\chi|_{X_i}=\chi(1)\phi_{a_i}$ where $i=8,9,10.$
Set $Z=F_{124}S_{567}X_8X_9X_{10}$ and the Kronecker $\delta_{i,j}=\left\{ \begin{array}{l}1 \mbox{ if } i=j, \\0 \mbox{ otherwise}\end{array}\right..$ The following are true.
\begin{itemize}

\item[(a)] If $\chi(1)=q^3,$ then $\chi=\chi_{8,9,10,q^3}^{a_8,a_9,a_{10}}$ and
\begin{center}
$\chi(\prod_{i=1}^{10}x_i(t_i))=\delta_{0,t_1}\delta_{0,t_2}\delta_{0,t_4}\delta_{0,t_3}\delta_{a_8t_5,a_{10}t_7}\delta_{a_8t_6,a_{9}t_7}q^3\phi(\sum_{i=8}^{10}a_it_i).$
\end{center}

\item[(b)] If $\chi(1)=\frac{q^3}{2},$ then $\chi=\chi_{8,9,10,\frac{q^3}{2}}^{b_{124},b_3,t_0,a_8,a_9,a_{10}}$ for some $b_{124},b_3\in \F_2,$ $t_0\in\F_q^\times$  and
$\chi(\prod_{i=1}^{10}x_i(t_i))=$
\[\left\{
\begin{array}{l}
\frac{q^3}{2}\phi(b_{124}\frac{t_1}{a_{10}}+At_0\frac{t_7}{a_8}+\sum_{i=8}^{10}a_it_i) \mbox{ if } \prod_{i=1}^{10}x_i(t_i)\in Z, \mbox{ and}
\\
\delta_{a_8t_1,a_{10}t_4}\delta_{a_8t_2,a_9t_4}\delta_{t_3,t_0^\phi}\frac{q^2}{2}\phi(b_{124}\frac{t_1}{a_{10}}+b_3t_3+At_0\frac{t_7}{a_8}+(*)+\sum_{i=8}^{10}a_it_i)\mbox{ otherwise,}
\end{array}\right.\]
where $t_0^\phi=\frac{(t_0)_\phi}{A}$ and $(*)=\frac{A^2}{(t_0)_\phi}(\frac{t_5}{a_{10}}+\frac{t_7}{a_8})(\frac{t_6}{a_9}+\frac{t_7}{a_8}).$
\end{itemize}
\end{theorem}
\noindent{\it Proof.} See Subsection \ref{proof-thm:Main2Chars}.

%************************************************************************************
%************************************************************************************

\section{Sylow $3$-subgroups of the Chevalley groups $E_6(3^f)$} \label{Section-E6}

Let $\F_q$ be a field of order $q$ and characteristic 3. We study $E_6(q)$ by its Lie root system. Let $\Sig:=\la \alpha_1,\alpha_2,\alpha_3,\alpha_4, \alpha_5, \alpha_6\ra$ be the root system
of $E_6,$ see Carter \cite{cart2}, Chapter 3. The Dynkin diagram of $\Sig$ is

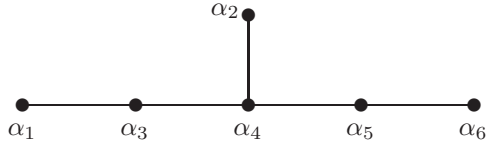
\begin{center}
\setlength{\unitlength}{1cm}
\begin{picture}(8,2)
\thinlines
\put(0.5,0.5){\circle*{0.17}}
\put(2,0.5){\circle*{0.17}}
\put(3.5,0.5){\circle*{0.17}}
\put(5,0.5){\circle*{0.17}}
\put(6.5,0.5){\circle*{0.17}}
\put(3.5,1.7){\circle*{0.17}}

\put(0.5,0.5){\line( 1, 0){1.5}}
\put(2,0.5){\line( 1, 0){1.5}}
\put(3.5,0.5){\line( 1, 0){1.5}}
\put(5,0.5){\line( 1, 0){1.5}}
\put(3.5,1.7){\line( 0, -1){1.2}}

\put(0.3,0.1){$\alpha_1$}
\put(1.8,0.1){$\alpha_3$}
\put(3.3,0.1){$\alpha_4$}
\put(4.8,0.1){$\alpha_5$}
\put(6.3,0.1){$\alpha_6$}
\put(3,1.7){$\alpha_2$}
\end{picture}
\end{center}

The positive roots are those roots which can be written as integral linear
combinations of the simple roots $\alpha_1, \alpha_2, ..., \alpha_6$ with nonnegative coefficients. We write $\Sig^+$ for the
set of positive roots. Here, $|\Sig^+|=36.$ We use the notation \Esix122321 for the
root $\alpha_1 + 2\alpha_2 + 2 \alpha_3 + 3\alpha_4+2\alpha_5+\alpha_6$ and we use a
similar notation for the remaining positive roots. Let $X_\alpha := \langle x_\alpha(t) \, | \, t \in \F_q \rangle$ be the
root subgroup corresponding to $\alpha \in \Sig.$ The group generated by all $X_\alpha$ for $\al\in\Sig^+$ is a
Sylow $3$-subgroup of the Chevally group $E_6(q),$ which we call $UE_6.$

In this section, we are going to construct irreducible characters $\chi$ of degree $\frac{q^7}{3}$ by considering the following special family of irreducible characters of $UE_6.$
 $$\FF_6:=\{\chi\in Irr(UE_6):\chi|_{X_\alpha}=\chi(1)\phi_a,\ ht(\alpha)=4,\ a\in \F_q^\times\}.$$
Let $\psi$ be an affording representation of some $\chi\in \FF_6.$ Using the same argument as in Section \ref{Section-D4} for all positive roots $\alpha$ with height greater than 4 to obtain $X_\alpha\subset ker(\chi).$ %Therefore,
Let $K_5$ be the normal subgroup of $UE_6$ generated by all root subgroups of height greater than $4$. Thus only the factor group $U := UE_6/K_5$ acts on a module affording $\chi.$ By the nature of the canonical map from $UE_6$ to $U,$ we can identify all root groups of root heights less than or equal 4 to their image groups. There are 21 roots $\alpha\in\Sig^+$ with $ht(\alpha)\leq 4.$ These $21$ positive
roots are given in Table~\ref{tab:posrootsE6}. Therefore, the group $U$ has order $q^{21}$ and $Z(U)=X_{17}X_{18}X_{19}X_{20}X_{21}=\la X_\beta:ht(\beta)=4\ra.$

For positive roots, we use the abbreviation
$x_i(t) = x_{\alpha_i}(t),$ $i=1,2,\dots,21.$ Each element $u \in U$ can be written uniquely as
\begin{center}
$u = x_2(t_2) x_1(t_1) x_3(t_3) x_4(t_4) x_5(t_5)  \cdots x_{21}(t_{21})$
where $x_i(t_i)\in X_i.$
\end{center}
So we write $\prod_{i=1}^{21}x_i(t_i)$ as this order. It is noted that there is a permutation of $x_2.$

\begin{table}[!ht]
\caption{Positive roots of the root system $\Sigma$ of type $E_6.$}
\label{tab:posrootsE6}

\begin{center}
\begin{tabular}{c|lll}
\hline
\rule{0cm}{0.4cm}
Height & Roots &&
\rule[-0.1cm]{0cm}{0.4cm}\\
\hline% \cline{1-3} \hline
\rule{0cm}{0.5cm}
4 &$\alpha_{20} := $ \Esix010111 &$\alpha_{21} := $ \Esix001111& \\
&&&\\
&$\alpha_{17} := $ \Esix111100 & $\alpha_{18} := $ \Esix101110&$\alpha_{19} := $ \Esix011110
\rule[-0.2cm]{0cm}{0.4cm}\\
\hline
\rule{0cm}{0.4cm}
3 & $\alpha_{15} := $ \Esix001110&$\alpha_{16} := $ \Esix000111&\\
&&&\\
&$\alpha_{12} := $ \Esix101100 &$\alpha_{13} := $ \Esix011100&$\alpha_{14} := $ \Esix010110
\rule[-0.2cm]{0cm}{0.4cm}\\
\hline
\rule{0cm}{0.4cm}
2 & $\alpha_{10} := $ \Esix000110&$\alpha_{11} := $ \Esix000011&\\
&&&\\
&$\alpha_{7} := $ \Esix101000&$\alpha_{8} := $ \Esix010100&$\alpha_{9} := $ \Esix001100
\rule[-0.2cm]{0cm}{0.4cm}\\
\hline
\rule{0cm}{0.4cm}
1 & $\alpha_2$\hspace{1cm}  $\alpha_1$   & $\alpha_3$ \hspace{1cm} $\alpha_4$ &$\alpha_5$ \hspace{1cm} $\alpha_6$
\rule[-0.2cm]{0cm}{0.4cm}\\
\hline
\end{tabular}
\end{center}
\end{table}

 For each $\alpha\in \Sig,$ since the lengths of $\alpha$-chains of roots through a root are at most $1,$ the commutator formula $[x_{\alpha}(r), x_{\beta}(s)] = x_{\alpha + \beta}(-C_{\alpha,\beta}rs)$ if $\alpha + \beta \in\Sig,$ and $= 1$ otherwise, see Cater \cite{cart2}, Theorem~5.2.2. For each extraspecial pair $(\alpha,\beta),$ we choose the coefficient $C_{\alpha,\beta}:=-1.$
By computing directly or using MAGMA \cite{MAG} with the following codes, all nontrivial commutators are given in Table~\ref{tab:commrelE6}.
%\vspace{0.1in}
%\\

W$:=$RootDatum("E6");
%\\

R$:=$PositiveRoots(W); A$:=$R[1..21];
%\\

for i in [7..21] do
%\\

${}^{}$\hspace{0.1in} for j in [1..(i-1)] do
%\\

${}^{}$\hspace{0.2in} if (R[i]-R[j]) in A then
%\\

${}^{}$\hspace{0.3in} k$:=$RootPosition(W,R[i]-R[j]);
%\\

${}^{}$\hspace{0.3in} if k le j then print k,"+",j,"=",i,"(",LieConstant$\_$C(W,1,1,k,j),")"; end if;
%\\

${}^{}$\hspace{0.2in} end if;
%\\

${}^{}$\hspace{0.1in} end for;
%\\

end for;

\begin{table}[!ht]
\caption{Commutator relations for type $E_6.$}
\label{tab:commrelE6}
%\begin{center}

\begin{tabular}{lll}
$\left[x_1(t), x_3(u)\right]=x_7(tu),$ &
$\left[x_2(t), x_4(u)\right]=x_8(tu),$ &
$\left[x_3(t), x_4(u)\right]=x_9(tu),$ \\
$\left[x_4(t), x_5(u)\right]=x_{10}(tu),$ &
$\left[x_5(t), x_6(u)\right]=x_{11}(tu),$ &
$\left[x_1(t), x_9(u)\right]=x_{12}(tu),$\\
$\left[x_4(t), x_7(u)\right]=x_{12}(-tu),$ &
$\left[x_2(t), x_9(u)\right]=x_{13}(tu),$ &
$\left[x_3(t), x_8(u)\right]=x_{13}(tu),$ \\
$\left[x_2(t), x_{10}(u)\right]=x_{14}(tu),$ &
$\left[x_5(t), x_8(u)\right]=x_{14}(-tu),$ &
$\left[x_3(t), x_{10}(u)\right]=x_{15}(tu),$ \\
$\left[x_5(t), x_9(u)\right]=x_{15}(-tu),$ &
$\left[x_4(t), x_{11}(u)\right]=x_{16}(tu),$ &
$\left[x_6(t), x_{10}(u)\right]=x_{16}(-tu),$ \\
$\left[x_1(t), x_{13}(u)\right]=x_{17}(tu),$ &
$\left[x_7(t), x_8(u)\right]=x_{17}(tu),$ &
$\left[x_2(t), x_{12}(u)\right]=x_{17}(tu),$\\
$\left[x_1(t), x_{15}(u)\right]=x_{18}(tu),$ &
$\left[x_7(t), x_{10}(u)\right]=x_{18}(tu),$ &
$\left[x_5(t), x_{12}(u)\right]=x_{18}(-tu),$ \\
$\left[x_2(t), x_{15}(u)\right]=x_{19}(tu),$ &
$\left[x_3(t), x_{14}(u)\right]=x_{19}(tu),$ &
$\left[x_5(t), x_{13}(u)\right]=x_{19}(-tu),$\\
$\left[x_2(t), x_{16}(u)\right]=x_{20}(tu),$ &
$\left[x_8(t), x_{11}(u)\right]=x_{20}(tu),$ &
$\left[x_6(t), x_{14}(u)\right]=x_{20}(-tu),$ \\
$\left[x_3(t), x_{16}(u)\right]=x_{21}(tu),$ &
$\left[x_9(t), x_{11}(u)\right]=x_{21}(tu),$ &
$\left[x_6(t), x_{15}(u)\right]=x_{21}(-tu).$
\end{tabular}
%\end{center}
\end{table}

\begin{center}
\begin{picture}(190,180)
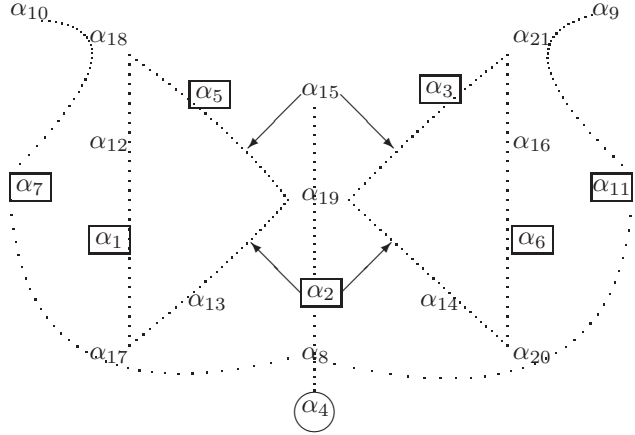
\put(40,0){Figure $UE_6(q):$ Relations of Roots.}
\put(-10,170){$\alpha_{10}$}
\put(20,160){$\alpha_{18}$}
\put(100,140){$\alpha_{15}$}
\put(180,160){$\alpha_{21}$}
\put(210,170){$\alpha_9$}
\put(20,120){$\alpha_{12}$}
\put(57,60){$\alpha_{13}$}
\put(100,100){$\alpha_{19}$}
\put(100,40){$\alpha_8$}
\put(180,40){$\alpha_{20}$}
\put(180,120){$\alpha_{16}$}
\put(145,60){$\alpha_{14}$}
\put(20,40){$\alpha_{17}$}

\put(-10,100){\framebox(15,10){$\alpha_7$}}
\put(20,80){\framebox(15,10){$\alpha_1$}}
\put(58,135){\framebox(15,10){$\alpha_5$}}
\put(100,60){\framebox(15,10){$\alpha_2$}}
\put(145,137){\framebox(15,10){$\alpha_3$}}
\put(180,80){\framebox(15,10){$\alpha_6$}}
\put(210,100){\framebox(15,10){$\alpha_{11}$}}

\put(105,20){\circle{15}}
\put(100,20){$\alpha_4$}

\qbezier[30](35,155)(35,100)(35,45)
\qbezier[30](35,155)(75,127)(95,100)
\qbezier[30](35,45)(75,77)(95,100)

\qbezier[30](178,155)(178,100)(178,45)
\qbezier[30](178,155)(148,132)(118,100)
\qbezier[30](118,100)(148,72)(178,45)

%\qbezier[30](35,135)(100,115)(178,135)
\put(100,140){\vector(-1,-1){20}}
\put(115,140){\vector(1,-1){20}}

%\qbezier[30](35,25)(100,58)(178,25)
\put(100,65){\vector(-1,1){19}}
\put(115,65){\vector(1,1){19}}

\qbezier[30](-10,100)(-5,15)(95,40)
\qbezier[30](225,100)(220,10)(105,40)

\qbezier[30](2,168)(45,165)(-10,110)
\qbezier[30](210,170)(170,160)(225,110)

\qbezier[20](105,135)(105,100)(105,72)
\qbezier[10](105,60)(105,40)(105,28)

\end{picture}
\end{center}

\smallskip
Let $H:=\la X_\alpha:\alpha_4\neq\alpha\in \Sig^+,(\alpha,\alpha_4)>0\ra=H_4H_3H_2$ where $H_4:=Z(U),$ $H_3:=\prod_{i=12}^{16}X_i,$ $H_2:=\prod_{i=8}^{10}X_i,$ and $T:=\la X_2,X_1,X_3,X_5,X_6\ra=X_2X_1X_3X_7X_5X_6X_{11}.$ It is clear that $|H|=q^{13},|T|=q^7,$ $H_k$ is generated by all root groups of root height $k$ in $H,$ and $T$ is a transversal of $HX_4$ in $U$.
Both $H$ and $HX_4$ are elementary abelian and normal in $U,$ and $T$ is isomorphic to $UA_2(q)\times UA_2(q)\times UA_1(q),$ where $UA_k(q)$ is the unipotent subgroup of the standard Borel subgroup of the general linear group $GL_{k+1}(q).$ We can visualize the group $U$ in the above figure. The roots in boxes are in $T,$ the others outside are in $H,$ and $\alpha_4,$ which is neither in $H$ nor in $T,$ is in a circle.  The disconnected lines demonstrate the relations between roots to give a sum root in center, e.g. $\alpha_7+\alpha_{10}=\alpha_{18},$ $\alpha_7+\alpha_8=\alpha_{17}...$ In addition, we have two triangles, as same as in Section \ref{Section-D4} of $UD_4(q)$, namely $(\alpha_{17},\alpha_{18},\alpha_{19})$ and $(\alpha_{19},\alpha_{20},\alpha_{21}).$ These two triangles share a common pair of roots $(\alpha_2,\alpha_{15})$ where $\alpha_2+\alpha_{15}=\alpha_{19}.$

We consider $\lambda\in Irr(H)$ such that $\lambda|_{X_i}=\phi_{a_i}\neq 1_{X_i}$ for $17\leq i\leq 21.$ Since the maximal split torus of $E_6(q)$ acts transitively on $\oplus_{i=17}^{21} Irr(X_i)^\times,$ it is allowed to assume that $\lambda|_{X_i}=\phi$ for $17\leq i\leq21.$ So we set $\lambda=\lambda_{b_8,b_9,b_{10}}^{b_{12},b_{13},b_{14},b_{15},b_{16}}\in Irr(H)$ such that $\lambda|_{X_i}=\phi_{b_i}$ where  $b_i\in\F_q$
 for all $8\leq i\leq 16,i\neq 11.$

\begin{definition}
For $b_8,b_9,b_{10},b_{12},b_{13},b_{14},b_{15},b_{16}\in\F_q,$ we define
\begin{itemize}

\item[(a)] $S_1:=\{s_1(t,r,s):=x_2(t)x_1(t)x_3(-t)x_5(t)x_6(-t)x_7(r)x_{11}(s):t,r,s\in\F_q\}.$

\item[(b)] $S_2:=\{s_2(t):=s_1(t,2t^2,2t^2)
% x_2(t)x_1(t)x_3(-t)x_5(t)x_6(-t)x_7(2t^2)x_{11}(2t^2)
:t\in\F_q\}.$
%where $r_7=2t^2-b_{14}t-b_{16}t+b_{15}t$ and $r_{11}=2t^2+b_{12}t+b_{13}t-b_{15}t$

\item[(c)] $R_3:=\{r_3(t):=x_{12}(t)x_{13}(-t)x_{14}(-t)x_{15}(t)x_{16}(t):t\in \F_q\}.$

\item[(d)] $R_2:=\{r_2(t):=x_8(-t)x_9(t)x_{10}(t):t\in\F_q\}.$

\item[(e)] $B_3:=b_{12}-b_{13}-b_{14}+b_{15}+b_{16}.$

\item[(f)] $B_2:=b_{10}+b_9-b_8.$%+b_{16}(b_{13}-b_{12}-b_{15})+b_{12}(b_{14}-b_{15})$

\item[(g)] If $B_2=c^2\in \F_q^\times,$ $F_2:=\{1,s_2(\pm c)\}$ and $F_4:=\{1,x_4(\pm c_\phi)\}.$% otherwise $F_2:=F_4:=\{1\}.$

\end{itemize}
\end{definition}

We note that $R_k\leq H_k$ for $k=2,3,$ $F_2\leq S_2\leq S_1\leq T,$ and $F_4\leq X_4.$  Since $R_k\cong \F_q,$ for each $a\in \F_q$ we define $\phi_a(r_k(t))=\phi_a(t)$ for all $r_k(t)\in R_k.$ Hence, $Irr(R_k)=\{\phi_a:a\in\F_q\}.$ Since $S_2 \cong \F_q,$ we can define $\phi_a(s_2(t)) = \phi_a(t)$ for all $s_2(t) \in S_2$.
When $B_2=c^2\in\F_q^\times,$ for each linear $\xi\in Irr(F_2)$ there is $b_2\in \{0,\pm a_2\}\cong (\F_3,+)$ such that $\xi=\phi_{b_2}|_{F_2}$ where $\phi_{b_2}\in Irr(S_2)$ and $\phi(a_2c)\neq1.$ Use the same argument for $F_4,$ for each $\xi\in Irr(F_4)$ there is $b_4\in\{0,\pm a_4\}\cong (\F_3,+)$ such that $\xi=\phi_{b_4}|_{F_4},$ where $\phi_{b_4}\in Irr(X_4)$ and $\phi(a_4c_\phi)\neq 1.$

Let $\overline{H_3}$ be the normal closure of $H_3$ in $HX_4S_1.$ Since $HX_4$ is abelian, $X_4\subset Stab_{U}(\lambda).$ All properties of $\lambda=\lambda_{b_8,b_9,b_{10}}^{b_{12},b_{13},b_{14},b_{15},b_{16}}$ are known as follows.

\begin{lemma} \label{lem:3extensions}
The following are true
\begin{itemize}
\item[(a)] $R_3=\{x\in H_3: |\lambda^U(x)|=\lambda^U(1)\}$ and $S_1 =Stab_{T}(\lambda|_{H_4H_3}).$ Moreover, $\lambda^U|_{R_3}=\lambda^U(1)\phi_{B_3}.$

\item[(b)] If $B_3\neq 0,$ then $Stab_T(\lambda)=\{1\}.$ Hence, if $\eta$ is an extension of $\lambda$ to $HX_4,$ then $I_U(\eta)=HX_4.$

%If $B_3=0,$ then otherwise $Stab_T(\lambda)=\{1\}.$

\item[(c)]  If $B_3=0,$ then there exists $x\in T$ such that ${}^{x}\lambda=\lambda^{0,0,0,0,0}_{b_8',b_9',b_{10}'}$ for some $b_8',b_9',b_{10}'\in \F_q.$ %and $Stab_{HT}({}^x\lambda)=HS_2.$
     Furthermore, $\overline{H_3}\subset ker({}^x\lambda)^{HX_4S_1}$ and the induction map from $Irr(HX_4S_1,{}^x\lambda)$ to $Irr(U,\lambda)$ is bijective.
\end{itemize}
\end{lemma}
\noindent{\it Proof.} See Subsection \ref{proof-lem:3extensions}$.~\Box$

\medskip
{\bf Remark}
If $\gcd(q,3)=1,$ then $\{x\in H_3:|\lambda^U(x)|=\lambda^U(1)\}$ and $Stab_T(\lambda)$ are trivial. Thus $\lambda$
extends to $HX_4$ and hence induces up to $U$ irreducibly.

\medskip
By Lemma \ref{lem:3extensions} (a), it is easy to see that $T$ acts invariant on $B_3=B_3(\lambda),$ i.e. $B_3(\lambda)=B_3({}^x\lambda)$ for all $x\in T.$ As above we fix the actions of $\lambda|_{X_i}=\phi,17\leq i\leq 21,$ $H$ has $q^8$ linears, in which there are $q^7$ linears with $B_3=0$ and $q^7(q-1)$ linears with $B_3\neq 0.$

By Lemma \ref{lem:3extensions} (b), these $q^7(q-1)$ linears of $H$ with $B_3\neq 0$ extend to $HX_4$ to be $q^8(q-1)$ linears and induce irreducibly to $U$ of degree $[U:HX_4]=q^7.$ Therefore, there are $\frac{q^8(q-1)}{q^7}=q(q-1)$ irreducibles in this case and they are parametrized by $(b_4,B_3)$. %where $B_3\in\F_q^\times$ and $b_4\in \F_q.$
So we denote them by $\chi_{q^7}^{b_4,B_3}$ where $b_4\in \F_q$ and $B_3\in\F_q^\times.$

Since $H\unlhd U,$ we have $\lambda,{}^x\lambda\in Irr(H)$ and $Irr(U,\lambda)=Irr(U,{}^x\lambda)$ for all $x\in T,$ hence, by Lemma \ref{lem:3extensions} (c), we suppose that $\lambda=\lambda^{0,0,0,0,0}_{b_8,b_9,b_{10}}.$ Since $[U:HX_4S_1]=q^4$ and the induction map from $HX_4S_1$ to $U$ is irreducible from all constituents of $Irr(HX_4S_1,\lambda),$ the above $q^7$ linears of $H$ with $B_3=0$ are corresponding with these $q^3$ linears $\lambda^{0,0,0,0,0}_{b_8,b_9,b_{10}}$ when we observe them at level of $HX_4S_1.$

\begin{lemma}
\label{lem:3ext_B1=0}
The following are true.
\begin{itemize}

\item[(a)] $R_2=\{x\in H_2: |\lambda^{HX_4S_1}(x)|=\lambda^{HX_4S_1}(1)\}$ and $S_2 =Stab_{S_1}(\lambda).$ Moreover, $\lambda^{HX_4S_1}|_{R_2}=\lambda^{HX_4S_1}(1)\phi_{B_2}.$

\item[(b)] If $B_2\notin \{c^2:c\in\F_q^\times\}$ and let $\eta$ be an extension of $\lambda$ to $HX_4,$ then $I_{HX_4S_1}(\eta)=HX_4.$ Therefore, $S_2$ acts transitively and faithfully on all extensions of $\lambda$ to $HX_4.$

\item[(c)] If $B_2=c^2\in\F_q^\times,$ then $\lambda$ extends to $HX_4F_2$ and $HF_4S_2.$ Let $\lambda_1,\lambda_2$ be extensions of $\lambda$ to $HX_4F_2.$ Then $I_{HX_4S_1}(\lambda_1)=HX_4F_2.$ Moreover, ${\lambda_1}^{HX_4S_1}={\lambda_2}^{HX_4S_1}$ iff $\lambda_1|_{F_2}=\lambda_2|_{F_2}$ and $\lambda_1|_{F_4}=\lambda_2|_{F_4}.$

\end{itemize}

\end{lemma}
\noindent{\it Proof.} See Subsection \ref{proof-lem:3ext_B1=0}$.~\Box$

\medskip
{\bf Remark}
When $B_2 = c^2 \neq 0$, we see that $HX_4F_3\unlhd U$ and $HF_4S_2\ntrianglelefteq U,$ and both have index $\frac{q^7}{3}$ in $U$. By Lemma \ref{lem:3extensions}(c) and Lemma \ref{lem:3ext_B1=0} (c) all constituents of $\lambda^U$ have degree $\frac{q^7}{3}.$ Hence, if $\eta$ is an extension of $\lambda$ to $HF_4S_2,$ then $\eta^U\in Irr(U,\lambda).$  We have $X_4,S_2\subset I_U(\lambda)$ and $\lambda$ extends to $HX_4F_3$ and $HF_4S_2,$ but $\lambda$ is not able to extend to $HX_4S_2.$

\medskip
The group $HX_4$ has $q^4$ linear characters $\lambda$ such that $\lambda|_H=\lambda_{b_8,b_9,b_{10}}^{0,0,0,0,0}.$ Since $\F_q^\times$ is even and cyclic, there are $\frac{q^3(q+1)}{2}$ linears with $B_2\notin \{c^2:c\in\F_q^\times\},$ and $\frac{q^3(q-1)}{2}$ linears with $B_2\in  \{c^2:c\in\F_q^\times\}.$ Hence, by Lemma \ref{lem:3ext_B1=0} (b), there are $\frac{q^3(q+1)}{2|S_1|}=\frac{q+1}{2}$ irreducibles of degree $|S_1|=q^3$ which are parametrized by $B_2\notin \{c^2:c\in\F_q^\times\}.$ By Lemma \ref{lem:3extensions} (c), we obtain $\frac{q+1}{2}$ irreducibles of degree $q^3[U:HX_4S_1]=q^7$ which are denoted by $\chi_{q^7}^{B_2}$ where  $B_2\in \F_q-\{c^2:c\in\F_q^\times\}.$ Therefore, together with characters $\chi_{q^7}^{b_4,B_3}$ as computed above, $\FF_6$ has exactly $(q-1)q+\frac{q+1}{2}$ irreducible characters $\chi$ of degree $q^7$ such that $\chi|_{X_i}=\chi(1)\phi$ for all $X_i\subset Z(U).$

By Lemma \ref{lem:3ext_B1=0} (c),  let $\lambda_1$ be an extension of $\lambda$ to $HX_4F_2,$ then ${\lambda_1}^{HX_4S_1}$ is irreducible of degree $[HX_4S_1:HX_4F_2]=\frac{q^3}{3}.$ These ${\lambda_1}^{HX_4S_1}$ only depend on $B_2$ and their restrictions to $F_2,F_4.$ Hence, by Lemma \ref{lem:3extensions} (c), ${\lambda_1}^U\in Irr(U)$ of degree $\frac{q^7}{3}$ is denoted by $\chi_{\frac{q^7}{3}}^{b_2,b_4,B_2}$ where $b_2,b_4\in \F_3$ and $B_2\in\{c^2:c\in \F_q^\times\}.$ Therefore, $\FF_6$ has exactly $\frac{9(q-1)}{2}$ irreducibles of degree $\frac{q^7}{3}$ such that $\chi|_{X_i}=\chi(1)\phi$ for all $X_i\subset Z(U).$

By the transitivity of the conjugate action of the maximal split torus $T_0$ of the Chevalley group $E_6(q)$ on $\oplus_{i=17}^{21}Irr(X_i)^\times,$ there are $(q-1)^5(q^2-q+\frac{q+1}{2})$ characters $\chi\in \FF_6$ of degree $q^7,$ and $\frac{9(q-1)^6}{2}$ characters $\chi\in \FF_6$ of degree $\frac{q^7}{3}$ such that $\chi|_{X_i}=\chi(1)\phi_{a_i},$ where $a_i\in\F_q^\times,17\leq i\leq 21.$ This gives the proof for the next theorem.

\begin{theorem} \label{thm:Main3Chars}

Let $\chi\in \FF_6.$ The following are true.
%For each $a_i\in \F_q^\times, 17\leq i\leq 21,$ suppose there is $\chi\in Irr(U)$ such that $\chi|_{X_i}=\chi(1)\phi_{a_i},17\leq i\leq 21.$ The following are true.
\begin{itemize}
\item[(a)] If $\chi(1)=q^7,$ then there exists $t\in T_0$ such that ${}^t\chi$ is either $\chi_{q^7}^{b_4,B_3}$ or $\chi_{q^7}^{B_2},$ for some $b_4\in \F_q,$ $B_3\in \F_q^\times,$ and $B_2\in \F_q-\{c^2:c\in\F_q^\times\}.$

\item[(b)] If $\chi(1)=\frac{q^7}{3},$ then there exists $t\in T_0$ such that ${}^t\chi=\chi_{\frac{q^7}{3}}^{b_2,b_4,B_2},$ for some $b_3,b_4\in \F_3$ and $B_2\in \{c^2:c\in \F_q^\times\}.$

\end{itemize}

\end{theorem}

%************************************************************************************
%************************************************************************************

\section{Sylow $5$-subgroups of the Chevalley groups $E_8(5^f)$} \label{Section-E8}

Let $\F_q$ be a field of order $q$ and characteristic 5. We study $E_8(q)$ by its Lie root system. Let $\Sig:=\la \alpha_1,\alpha_2,\alpha_3,\alpha_4, \alpha_5, \alpha_6, \alpha_7, \alpha_8 \ra$ be the root system
of $E_8,$ see Carter \cite{cart2}, Chapter 3. The Dynkin diagram of $\Sig$ is

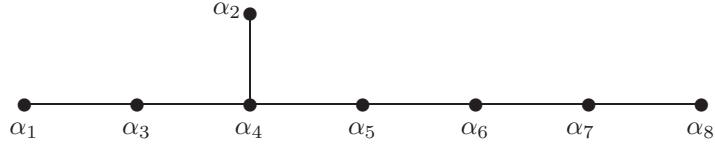
\begin{center}
\setlength{\unitlength}{1cm}
\begin{picture}(10,2)
\thinlines
\put(0.5,0.5){\circle*{0.17}}
\put(2,0.5){\circle*{0.17}}
\put(3.5,0.5){\circle*{0.17}}
\put(5,0.5){\circle*{0.17}}
\put(6.5,0.5){\circle*{0.17}}
\put(8,0.5){\circle*{0.17}}
\put(9.5,0.5){\circle*{0.17}}

\put(3.5,1.7){\circle*{0.17}}

\put(0.5,0.5){\line( 1, 0){1.5}}
\put(2,0.5){\line( 1, 0){1.5}}
\put(3.5,0.5){\line( 1, 0){1.5}}
\put(5,0.5){\line( 1, 0){1.5}}
\put(6.5,0.5){\line( 1, 0){1.5}}
\put(8,0.5){\line( 1, 0){1.5}}

\put(3.5,1.7){\line( 0, -1){1.2}}

\put(0.3,0.1){$\alpha_1$}
\put(1.8,0.1){$\alpha_3$}
\put(3.3,0.1){$\alpha_4$}
\put(4.8,0.1){$\alpha_5$}
\put(6.3,0.1){$\alpha_6$}
\put(7.7,0.1){$\alpha_7$}
\put(9.3,0.1){$\alpha_8$}

\put(3,1.7){$\alpha_2$}
\end{picture}
\end{center}

The positive roots are those roots which can be written as linear
combinations of the simple roots $\alpha_1, \alpha_2, ..., \alpha_8$ with nonnegative coefficients and we write $\Sig^+$ for the
set of positive roots. Here, $|\Sig^+|=120.$ We use the notation $\Eeight23465432$ for the
root $2\alpha_1 + 3\alpha_2 + 4 \alpha_3 + 6\alpha_4+ 5\alpha_5+ 4\alpha_6+ 3\alpha_7+ 2\alpha_8$ and we use a
similar notation for the remaining positive roots. Let $X_\alpha := \langle x_\alpha(t) \, | \, t \in \F_q \rangle$ be the
root subgroup corresponding to $\alpha \in \Sig.$ The group generated by all $X_\alpha$ for $\al\in\Sig^+$ is a Sylow $5$-subgroup of the Chevalley group $E_8(q),$ which we call $UE_8.$

In this section, we are going to construct irreducible characters $\chi$ of degree $\frac{q^{16}}{5}$ by considering the following special family of irreducible characters of $UE_8.$
 \[\FF_8:=\{\chi\in Irr(UE_8):\chi|_{X_\alpha}=\chi(1)\phi_a,\ ht(\alpha)=6,\ a\in \F_q^\times\}.\]
Let $\psi$ be an affording representation of some $\chi\in \FF_8.$ Using the same argument as in Section \ref{Section-D4} for all positive roots $\alpha$ with height greater than 6 to obtain $X_\alpha\subset ker(\chi).$ %Therefore,
Let $K_7$ be the normal subgroup of $UE_8$ generated by all root subgroups of root heights greater than $6$. Thus only the factor group $U := UE_8/K_7$ acts on a module affording $\chi.$ By the nature of the canonical map from $UE_8$ to $U,$ we can identify all root groups of root heights less than or equal 6 to their image groups. There are 43 positive roots of height less than or equal $6.$ These $43$
roots are given in Table~\ref{tab:posrootsE8}.

For positive roots, we use the abbreviation
$x_i(t) = x_{\alpha_i}(t),$ $i=1,2,\dots,43.$ Hence, $Z(U)=X_{37}X_{38}X_{39}X_{40}X_{41}X_{42}X_{43}=\la X_\beta:ht(\beta)=6\ra.$ Each element $u \in U$ can be written uniquely as
\begin{center}
$u = x_2(t_2) x_1(t_1) x_3(t_3) x_4(t_4) x_5(t_5)  \cdots x_{43}(t_{43})$
where $x_i(t_i)\in X_i.$
\end{center}
So we write $\prod_{i=1}^{43}x_i(t_i)$ as this order. It is noted that there is a permutation of $x_2.$

\begin{table}[!ht]
\caption{Positive roots of the root system $\Sigma$ of type $\widetilde{E_8}.$}
\label{tab:posrootsE8}
\vspace{-0.2in}
\begin{center}
\begin{tabular}{c|lll}
\hline
\rule{0cm}{0.4cm}
Height & Roots &&
\rule[-0.1cm]{0cm}{0.4cm}\\
\hline% \cline{1-3} \hline
\rule{0cm}{0.5cm}
6 &$\alpha_{43} := $ \Eeight00111111 &&\\
&&&\vspace{-0.1in}\\
&$\alpha_{40} := $ \Eeight01121100&$\alpha_{41} := $ \Eeight01111110&$\alpha_{42} := $ \Eeight01011111
\rule[-0.2cm]{0cm}{0.4cm}\\
&&&\vspace{-0.1in}\\
&$\alpha_{37} := $ \Eeight11121000&$\alpha_{38} := $ \Eeight11111100&$\alpha_{39} := $ \Eeight10111110\\
\hline
\rule{0cm}{0.4cm}
5 &$\alpha_{36} := $ \Eeight00011111&& \\
&&&\vspace{-0.1in}\\
&$\alpha_{33} := $ \Eeight01111100&$\alpha_{34}:=$  \Eeight01011110  & $\alpha_{35} := $ \Eeight00111110\\
&&&\vspace{-0.1in}\\
& $\alpha_{30}:=$ \Eeight11111000&$\alpha_{31} := $ \Eeight10111100 &$\alpha_{32} := $ \Eeight01121000
\rule[-0.2cm]{0cm}{0.4cm}\\
\hline
\rule{0cm}{0.4cm}
4 & $\alpha_{29} := $ \Eeight00001111&&\\
&&&\vspace{-0.1in}\\
&$\alpha_{26} := $ \Eeight01011100&$\alpha_{27} := $ \Eeight00111100&$\alpha_{28} := $ \Eeight00011110
\rule[-0.2cm]{0cm}{0.4cm}\\
&&&\vspace{-0.1in}\\
&$\alpha_{23} := $ \Eeight11110000&$\alpha_{24} := $ \Eeight10111000&$\alpha_{25} := $ \Eeight01111000\\
\hline
\rule{0cm}{0.4cm}
3 & $\alpha_{22} := $ \Eeight00000111 && \\
&&&\vspace{-0.1in}\\
&$\alpha_{19} := $ \Eeight00111000 &$\alpha_{20} := $ \Eeight00011100&$\alpha_{21} := $ \Eeight00001110
\rule[-0.2cm]{0cm}{0.4cm}\\
&&&\vspace{-0.1in}\\
&$\alpha_{16} := $ \Eeight10110000&$\alpha_{17} := $ \Eeight01110000&$\alpha_{18} := $ \Eeight01011000\\
\hline
\rule{0cm}{0.4cm}
2 & $\alpha_{15} := $ \Eeight00000011&&\\
&&&\vspace{-0.1in}\\
&$\alpha_{12} := $ \Eeight00011000& $\alpha_{13} := $ \Eeight00001100&$\alpha_{14} := $ \Eeight00000110 \\
&&&\vspace{-0.1in}\\
&$\alpha_{9} := $ \Eeight10100000&$\alpha_{10} := $ \Eeight01010000&$\alpha_{11} := $ \Eeight00110000
\rule[-0.2cm]{0cm}{0.4cm}\\
\hline
\rule{0cm}{0.4cm}
1 & $\alpha_2$\hspace{0.9cm}  $\alpha_1$\hspace{0.9cm} $\alpha_3$  & $\alpha_4$ \hspace{0.9cm} $\alpha_5$\hspace{0.9cm}$\alpha_6$ &$\alpha_7$ \hspace{0.9cm} $\alpha_8$
\rule[-0.2cm]{0cm}{0.4cm}\\
\hline
\end{tabular}
\end{center}
\end{table}

For $\alpha,\beta\in \Sig,$ the commutator formula $[x_{\alpha}(r), x_{\beta}(s)] = x_{\alpha+\beta}(-C_{\alpha,\beta}rs)$ if $\alpha+\beta\in\Sig,$ $=1$ otherwise, see Cater \cite{cart2}, Theorem~5.2.2. For each extraspecial pair $(\alpha,\beta),$ we choose the coefficient $C_{\alpha,\beta}:=-1.$
By computing directly or using MAGMA \cite{MAG}, all nontrivial commutators are given in Table~\ref{tab:commrelE8}.

\begin{table}[!ht]
\caption{Commutator relations for type $\widetilde{E_8}.$}
\label{tab:commrelE8}

%\begin{center}
\begin{tabular}{lll}

$\left[x_1(t), x_3(u)\right]=x_9(tu),$ &
$\left[x_2(t), x_4(u)\right]=x_{10}(tu),$ &
$\left[x_3(t), x_4(u)\right]=x_{11}(tu),$ \\
$\left[x_4(t), x_5(u)\right]=x_{12}(tu),$ &
$\left[x_5(t), x_6(u)\right]=x_{13}(tu),$ &
$\left[x_6(t), x_7(u)\right]=x_{14}(tu),$\\
$\left[x_7(t), x_8(u)\right]=x_{15}(tu),$&
$\left[x_1(t), x_{11}(u)\right]=x_{16}(tu),$&
$\left[x_4(t), x_9(u)\right]=x_{16}(-tu),$ \\
$\left[x_2(t), x_{11}(u)\right]=x_{17}(tu),$&
$\left[x_3(t), x_{10}(u)\right]=x_{17}(tu),$&
$\left[x_2(t), x_{12}(u)\right]=x_{18}(tu),$\\
$\left[x_5(t), x_{10}(u)\right]=x_{18}(-tu),$ &
$\left[x_3(t), x_{12}(u)\right]=x_{19}(tu),$&
$\left[x_5(t), x_{11}(u)\right]=x_{19}(-tu),$ \\
$\left[x_4(t), x_{13}(u)\right]=x_{20}(tu),$&
$\left[x_6(t), x_{12}(u)\right]=x_{20}(-tu),$ &
$\left[x_5(t), x_{14}(u)\right]=x_{21}(tu),$\\
$\left[x_7(t), x_{13}(u)\right]=x_{21}(-tu),$ &
$\left[x_6(t), x_{15}(u)\right]=x_{22}(tu),$ &
$\left[x_8(t), x_{14}(u)\right]=x_{22}(-tu),$ \\
$\left[x_1(t), x_{17}(u)\right]=x_{23}(tu),$ &
$\left[x_2(t), x_{16}(u)\right]=x_{23}(tu),$ &
$\left[x_9(t), x_{10}(u)\right]=x_{23}(tu),$\\
$\left[x_1(t), x_{19}(u)\right]=x_{24}(tu),$ &
$\left[x_5(t), x_{16}(u)\right]=x_{24}(-tu),$&
$\left[x_9(t), x_{12}(u)\right]=x_{24}(tu),$\\
$\left[x_2(t), x_{19}(u)\right]=x_{25}(tu),$ &
$\left[x_3(t), x_{18}(u)\right]=x_{25}(tu),$&
$\left[x_5(t), x_{17}(u)\right]=x_{25}(-tu),$\\
$\left[x_2(t), x_{20}(u)\right]=x_{26}(tu),$ &
$\left[x_6(t), x_{18}(u)\right]=x_{26}(-tu),$ &
$\left[x_{10}(t), x_{13}(u)\right]=x_{26}(tu),$ \\
$\left[x_3(t), x_{20}(u)\right]=x_{27}(tu),$ &
$\left[x_6(t), x_{19}(u)\right]=x_{27}(-tu),$&
$\left[x_{11}(t), x_{13}(u)\right]=x_{27}(tu),$\\
$\left[x_4(t), x_{21}(u)\right]=x_{28}(tu),$ &
$\left[x_7(t), x_{20}(u)\right]=x_{28}(-tu),$ &
$\left[x_{12}(t), x_{14}(u)\right]=x_{28}(tu),$ \\
$\left[x_5(t), x_{22}(u)\right]=x_{29}(tu),$ &
$\left[x_8(t), x_{21}(u)\right]=x_{29}(-tu),$&
$\left[x_{13}(t), x_{15}(u)\right]=x_{29}(tu),$\\
$\left[x_1(t), x_{25}(u)\right]=x_{30}(tu),$ &
$\left[x_2(t), x_{24}(u)\right]=x_{30}(tu),$ &
$\left[x_5(t), x_{23}(u)\right]=x_{30}(-tu),$  \\
$\left[x_9(t), x_{18}(u)\right]=x_{30}(tu),$&
$\left[x_1(t), x_{27}(u)\right]=x_{31}(tu),$&
$\left[x_6(t), x_{24}(u)\right]=x_{31}(-tu),$\\
$\left[x_9(t), x_{20}(u)\right]=x_{31}(tu),$ &
$\left[x_{13}(t), x_{16}(u)\right]=x_{31}(-tu),$ &
$\left[x_4(t), x_{25}(u)\right]=x_{32}(tu),$ \\
$\left[x_{10}(t), x_{19}(u)\right]=x_{32}(-tu),$ &
$\left[x_{11}(t), x_{18}(u)\right]=x_{32}(-tu),$&
$\left[x_{12}(t), x_{17}(u)\right]=x_{32}(-tu),$\\
$\left[x_2(t), x_{27}(u)\right]=x_{33}(tu),$ &
$\left[x_3(t), x_{26}(u)\right]=x_{33}(tu),$ &
$\left[x_6(t), x_{25}(u)\right]=x_{33}(-tu),$\\
$\left[x_{13}(t), x_{17}(u)\right]=x_{33}(-tu),$&
$\left[x_2(t), x_{28}(u)\right]=x_{34}(tu),$&
$\left[x_7(t), x_{26}(u)\right]=x_{34}(-tu),$\\
$\left[x_{10}(t), x_{21}(u)\right]=x_{34}(tu),$ &
$\left[x_{14}(t), x_{18}(u)\right]=x_{34}(-tu),$&
$\left[x_3(t), x_{28}(u)\right]=x_{35}(tu),$ \\
$\left[x_7(t), x_{27}(u)\right]=x_{35}(-tu),$&
$\left[x_{11}(t), x_{21}(u)\right]=x_{35}(tu),$ &
$\left[x_{14}(t), x_{19}(u)\right]=x_{35}(-tu),$\\
$\left[x_4(t), x_{29}(u)\right]=x_{36}(tu),$ &
$\left[x_8(t), x_{28}(u)\right]=x_{36}(-tu),$ &
$\left[x_{12}(t), x_{22}(u)\right]=x_{36}(tu),$ \\
$\left[x_{15}(t), x_{20}(u)\right]=x_{36}(-tu),$&
$\left[x_1(t), x_{32}(u)\right]=x_{37}(tu),$&
$\left[x_4(t), x_{30}(u)\right]=x_{37}(tu),$\\
$\left[x_{10}(t), x_{24}(u)\right]=x_{37}(-tu),$ &
$\left[x_{12}(t), x_{23}(u)\right]=x_{37}(-tu),$ &
$\left[x_{16}(t), x_{18}(u)\right]=x_{37}(-tu),$ \\
$\left[x_1(t), x_{33}(u)\right]=x_{38}(tu),$ &
$\left[x_2(t), x_{31}(u)\right]=x_{38}(tu),$&
$\left[x_6(t), x_{30}(u)\right]=x_{38}(-tu),$\\
$\left[x_9(t), x_{26}(u)\right]=x_{38}(tu),$ &
$\left[x_{13}(t), x_{23}(u)\right]=x_{38}(-tu),$ &
$\left[x_1(t), x_{35}(u)\right]=x_{39}(tu),$ \\
$\left[x_7(t), x_{31}(u)\right]=x_{39}(-tu),$ &
$\left[x_9(t), x_{28}(u)\right]=x_{39}(tu),$ &
$\left[x_{14}(t), x_{24}(u)\right]=x_{39}(-tu),$\\
$\left[x_{16}(t), x_{21}(u)\right]=x_{39}(tu),$ &
$\left[x_4(t), x_{33}(u)\right]=x_{40}(tu),$ &
$\left[x_6(t), x_{32}(u)\right]=x_{40}(-tu),$ \\
$\left[x_{10}(t), x_{27}(u)\right]=x_{40}(-tu),$ &
$\left[x_{11}(t), x_{26}(u)\right]=x_{40}(-tu),$&
$\left[x_{17}(t), x_{20}(u)\right]=x_{40}(tu),$\\
$\left[x_2(t), x_{35}(u)\right]=x_{41}(tu),$ &
$\left[x_3(t), x_{34}(u)\right]=x_{41}(tu),$ &
$\left[x_7(t), x_{33}(u)\right]=x_{41}(-tu),$ \\
$\left[x_{14}(t), x_{25}(u)\right]=x_{41}(-tu),$ &
$\left[x_{17}(t), x_{21}(u)\right]=x_{41}(tu),$&
$\left[x_2(t), x_{36}(u)\right]=x_{42}(tu),$\\
$\left[x_8(t), x_{34}(u)\right]=x_{42}(-tu),$&
$\left[x_{10}(t), x_{29}(u)\right]=x_{42}(tu),$  &
$\left[x_{15}(t), x_{26}(u)\right]=x_{42}(-tu),$ \\
$\left[x_{18}(t), x_{22}(u)\right]=x_{42}(tu),$&
$\left[x_3(t), x_{36}(u)\right]=x_{43}(tu),$&
$\left[x_8(t), x_{35}(u)\right]=x_{43}(-tu),$\\
$\left[x_{11}(t), x_{29}(u)\right]=x_{43}(tu),$&
$\left[x_{15}(t), x_{27}(u)\right]=x_{43}(-tu),$ &
$\left[x_{19}(t), x_{22}(u)\right]=x_{43}(tu),$
\end{tabular}
%\end{center}
\end{table}

Let $H:=\la X_\al:\, \al_4\neq \al\in\Sig^+, (\al,\al_5)> 0\ra=H_6H_5H_4H_3H_2$ where $H_6=Z(U),$ $H_5=\prod_{i=30}^{36}X_i,$ $H_4=\prod_{i=24}^{29}X_i,$ $H_3=\prod_{i=18}^{21}X_i$ and $H_2=X_{12}X_{13}.$ Let $T:=\la X_1,X_3,X_4,X_2,X_6,X_7,X_8 \ra=T_4T_3T_2T_1$ where $T_4=X_{23},$ $T_3=X_{16}X_{17}X_{22},$ $T_2=X_9X_{10}X_{11}X_{14}X_{15}$ and $T_1=X_1X_3X_4X_2X_6X_7X_8.$ It is clear that $|H|=q^{26},$ $|T|=q^{16},$ $H_k$ is generated by all root groups in $H$ of root height $k,$ as same as for $T_k$ generated by all root subgroups in $T$ of root height $k,$ and $T$ is a transversal of $HX_5$ in $U.$ Both $H$ and $HX_5$ are elementary abelian and normal in $U.$ $T$ is isomorphic to $UA_4(q)\times UA_3(q),$ where $UA_k(q)$ is the unipotent subgroup of the standard Borel subgroup of the general linear group $GL_{k+1}(q).$ It is noted that by letting $\{\beta_1,\beta_2,\beta_3,\beta_4\}$ be a simple root set of type $A_4,$ the isomorphism from $\la X_1,X_3,X_4,X_2\ra$ to $UA_4(q)$  sends $x_1(t)$ to $x_{\beta_1}(t),$ $x_3(t)$ to $x_{\beta_2}(t),$ $x_4(t)$ to $x_{\beta_3}(t),$ and $x_2(t)$ to $x_{\beta_4}(-t)$ for all $t\in \F_q.$

We consider linear characters $\lam\in Irr(H)$ such that $\lam|_{X_i}=\phi_{a_i}$ for $37\leq i\leq 43$ and $\lam|_{X_j}=\phi_{b_j}$ for all appropriate $j\leq 36$ where $a_i\in\F_q^\times$ and $b_j\in \F_q.$ Since the maximal split torus of the Chevalley group $E_8(q)$ acts transitively on $\otimes_{i=37}^{43}Irr(X_i)^\times,$ it suffices to suppose that $\lam|_{X_i}=\phi$ for all $37\leq i\leq 43.$

\begin{definition}
\label{def:stabilizers}
For $b_i\in \F_q$ where $i\in[12..13,18..21,24..36]$ we define
\begin{itemize}
\item[(a)] $B_5:=b_{30}+b_{31}-b_{32}-b_{33}-2b_{34}+2b_{35}+2b_{36}.$

\item[(b)] $B_4:=2b_{24}-2b_{25}+b_{26}-b_{27}-b_{28}+b_{29}.$

\item[(c)] $B_3:=b_{18}-b_{19}-b_{20}+b_{21}.$

\item[(d)] $B_2:=b_{12}-b_{13}.$

\item[(e)] $R_5:=\{r_5(v):=x_{30}(v)x_{31}(v)x_{32}(-v)x_{33}(-v)x_{34}(-2v)x_{35}(2v)x_{36}(2v):v\in \F_q\}.$

\item[(f)] $R_4:=\{r_4(v):=x_{24}(2v)x_{25}(-2v)x_{26}(v)x_{27}(-v)x_{28}(-v)x_{29}(v): v\in \F_q\}.$

\item[(g)] $R_3:=\{r_3(v):=x_{18}(v)x_{19}(-v)x_{20}(-v)x_{21}(v): v\in \F_q\}.$

\item[(h)] $R_2:=\{r_2(v):=x_{12}(v)x_{13}(-v): v\in \F_q\}.$

\item[(i)] $L_1:=\{l_1(u):=x_2(2u)x_1(u)x_3(-2u)x_4(u)x_6(u)x_7(2u)x_8(-2u):\, u\in \F_q\},$ $S_1:=L_1T_2T_3T_4.$

\item[(j)] $L_2:=\{l_2(u):=l_1(u)x_9(u^2)x_{10}(-u^2)x_{11}(u^2)x_{14}(-u^2)x_{15}(2u^2):t\in \F_q\},$ $S_2:=L_2T_3T_4.$

\item[(k)] $L_3:=\{l_3(u):=l_2(u)x_{16}(4u^3)x_{17}(2u^3)x_{22}(3u^3): u\in \F_q\},$ $S_3:=L_3T_4$

\item[(l)] $S_4:=\{l_4(u):=l_3(u)x_{23}(3u^4): u\in\F_q\}.$

\item[(m)] If $B_2=c^4\in \F_q^\times,$ $F_4:=\{s_4(u c):u\in \F_5 \}$ and $F_5:=\{x_5(v c_\phi): v\in\F_5\}.$
\end{itemize}

\end{definition}

It is easy to check that for $k\in [2..5],$ $R_k\leq H_k$ of order $q,$ $S_k\leq S_{k-1}\leq T$ with $S_5=\{1\},$ and $F_4\leq S_4,$ $F_5\leq X_5$ of order $5.$ It is noted that all $B_i$ are defined for each $\lambda$ as above, hence $B_i=B_i(\lambda).$ Since $R_k\cong \F_q,$ for each $a\in \F_q$ we define $\phi_a(r_k(t))=\phi_a(t)$ for all $r_k(t)\in R_k.$ Hence, $Irr(R_k)=\{\phi_a:a\in\F_q\}.$ Since $S_4 \cong \F_q,$ we can define $\phi_a(s_4(t)) = \phi_a(t)$ for all $s_4(t) \in S_4$.
When $B_2=c^4\in\F_q^\times,$ for each linear $\xi\in Irr(F_4)$ there is $b_4\in \{t a_4: t\in \F_5\}\cong (\F_5,+)$ such that $\xi=\phi_{b_4}|_{F_4}$ where $\phi_{b_4}\in Irr(S_4)$ and $\phi(a_4 c)\neq 1.$ Use the same argument for $F_5\leq X_5,$ for each $\xi\in Irr(F_5)$ there is $b_5\in\{t a_5: t\in \F_5\}\cong (\F_5,+)$ such that $\xi=\phi_{b_5}|_{F_5},$ where $\phi_{b_5}\in Irr(X_5)$ and $\phi(a_5 c_\phi)\neq 1.$

Let $\overline{H_5}$ be the normal closure of $H_5$ in $HX_5S_1.$ All properties of $\lambda$'s are known as follows. It is clear that $X_5\subset Stab_U(\lambda).$

\begin{lemma} \label{lem:5ext1}
The following are true
\begin{itemize}
\item[(a)] $R_5=\{x\in H_5: |\lambda^U(x)|=\lambda^U(1)\}$ and $S_1 =Stab_{T}(\lambda|_{H_6H_5}).$ Moreover, $\lambda^U|_{R_5}=\lambda^U(1)\phi_{B_5}.$

\item[(b)] If $B_5\neq 0,$ then $Stab_T(\lambda)=\{1\}.$ Hence, if $\eta$ is an extension of $\lambda$ to $HX_5,$ then $I_U(\eta)=HX_5.$ Furthermore, if $\eta,\eta'$ are two extensions of  $\lambda|_{H_6H_5H_4}$ to $HX_5,$ then $\eta^U=\eta'^U$ iff $B_i(\eta)=B_i(\eta')$ for $i=2,3$ and $\eta|_{X_5}=\eta'|_{X_5}.$

\item[(c)]  If $B_5=0,$ then there exists $x\in T$ such that ${}^{x}\lambda|_{X_i}=1_{X_i}$ for all $X_i\subset H_5.$  Furthermore, $\overline{H_5}\subset ker({}^x\lambda)^{HX_5S_1}$ and the induction map from $Irr(HX_5S_1,{}^x\lambda)$ to $Irr(U,\lambda)$ is bijective.
\end{itemize}
\end{lemma}
\noindent{\it Proof.} See Subsection \ref{proof-lem:5ext1}$.~\Box$

\medskip
{\bf Remark} When $(q,5)=1,$ both $R_5$ and $Stab_T(\lambda)$ are trivial. Hence, $\lambda$ extends to $HX_5$ and induces irreducibly to $U$ of degree $[U:HX_5]=q^{16}.$

\medskip
Lemma \ref{lem:5ext1} (a) can be observed by the following figure.

\begin{center}
\begin{picture}(315,170)
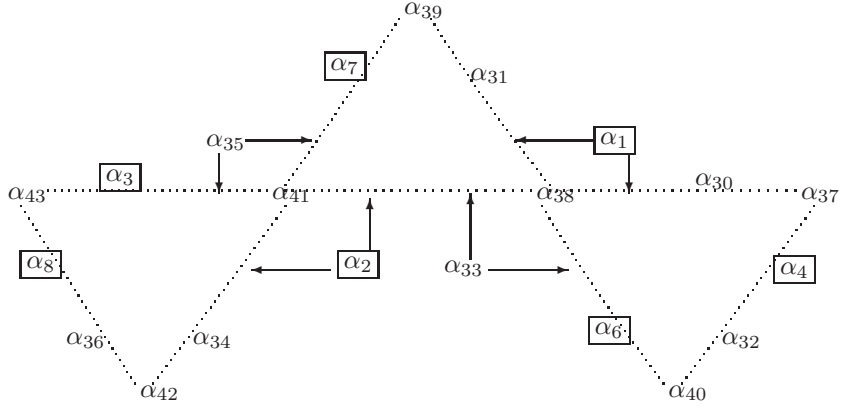
\put(-5,0){Figure $UE_8(q)$: Relations among roots of heights 5 in $H$ and 1 in $T$. }
\put(5,90){$\alpha_{43}$} %r
\put(55,15){$\alpha_{42}$}
\put(105,90){$\alpha_{41}$} %r
\put(205,90){$\alpha_{38}$} %r
\put(305,90){$\alpha_{37}$} %r
\put(255,15){$\alpha_{40}$}
\put(155,160){$\alpha_{39}$}

\qbezier[30](20,93)(60,93)(103,93)   %43 -> 41
\qbezier[30](110,93)(150,93)(203,93) %41 -> 38
\qbezier[30](210,93)(250,93)(303,93) %38 -> 37

\qbezier[30](53,20)(32,53)(10,87)   %42 -> 43
\qbezier[30](60,20)(85,53)(110,87) %42 -> 41
\qbezier[30](253,20)(228,53)(207,87)  %40 -> 38
\qbezier[30](260,20)(285,53)(310,87)  %40 -> 37
\qbezier[30](110,95)(132,126)(153,158)  %41 -> 39
\qbezier[30](210,95)(185,126)(165,158)  %38 -> 39

\put(80,110){$\alpha_{35}$} \put(95,112){\vector(1,0){25}} \put(85,107){\vector(0,-1){15}}
\put(227,107){\framebox(15,10){$\alpha_1$}} \put(227,112){\vector(-1,0){30}} \put(240,107){\vector(0,-1){15}}
\put(130,60){\framebox(15,10){$\alpha_2$}} \put(142,70){\vector(0,1){20}} \put(127,63){\vector(-1,0){30}}
\put(170,62){$\alpha_{33}$} \put(180,67){\vector(0,1){25}} \put(187,63){\vector(1,0){30}}

\put(75,35){{$\alpha_{34}$}}                %on 41-43
\put(40,93){\framebox(15,10){$\alpha_3$}}  %on 41-42
\put(125,135){\framebox(15,10){$\alpha_7$}}  %on 39-41
\put(180,135){{$\alpha_{31}$}}      %on 38-39
\put(265,95){{$\alpha_{30}$}}        %on 37-38
\put(225,35){\framebox(15,10){$\alpha_6$}}      %on 38-40
\put(275,35){{$\alpha_{32}$}}  %on 37-40
\put(295,58){\framebox(15,10){$\alpha_4$}}  %on 37-40
\put(27,35){{$\alpha_{36}$}}   %on 42-43
\put(10,60){\framebox(15,10){$\alpha_8$}}  %on 42-43
\end{picture}
\end{center}

\medskip
We have $q^{19}$ linear characters $\lambda$ of $H$ such that $\lambda|_{X_i}=\phi$ for all $X_i\subset Z(U).$ In these, there are $q^{18}$ linears with $B_5=0$ and $q^{18}(q-1)$ linears with $B_5\neq 0.$ By Lemma \ref{lem:5ext1} (a), it is clear that $B_5=B_5(\lambda)$ is invariant under the action of $T.$ Therefore, by Lemma \ref{lem:5ext1} (b), these $q^{18}(q-1)$ linears with $B_5\neq 0$ extend to $HX_5$ and induce irreducibly to $U.$ Thus, we obtain $\frac{q^{19}(q-1)}{q^{16}}=q^3(q-1)$ irreducible characters of $U$ of degree $q^{16}$ which are parametrized by $(b_5,B_2,B_3,B_5)$ where  $b_5,B_2,B_3\in \F_q$ and $B_5\in \F_q^\times.$ Hence, we denote them by $\chi_{q^{16}}^{b_5,B_2,B_3,B_5}.$

Since $[U:HX_5S_1]=q^6,$ by Lemma \ref{lem:5ext1}(c), the constituents of $q^{18}$ linears of $H$ with $B_5=0$ inducing to $U$ correspond to the ones of $q^{12}$ linears $\lambda$ of $H$ inducing to  $HX_5S_1,$ where $\lambda|_{X_i}=\phi$ for all $X_i\subset H_6,$ $\lambda|_{X_i}=1_{X_i}$ for all $X_i\subset H_5,$ and $\lambda|_{X_i}=\phi_{b_i}$ for all $X_i\subset H_4H_3H_2$ where $b_i\in\F_q.$ Let $\lambda\in Irr(H)$ be one of these above $q^{12}$ linears. %such that $\lambda|_{X_i}=\phi$ for all $X_i\subset H_6,$ $\lambda|_{X_i}=1_{X_i}$ for all $X_i\subset H_5,$ and $\lambda|_{X_i}=\phi_{b_i}$ for all $X_i\subset H_4H_3H_2$ where $b_i\in\F_q.$
Now we consider $\lambda$ in the scenario of the subgroup $HX_5S_1.$ Let $\overline{H_5H_4}$ be the normal closure of $H_5H_4$ in $HX_5S_2.$

%(Work with $HX_5S_1/\overline{H_5}$)

\begin{lemma} \label{lem:5ext2}
The following are true
\begin{itemize}
\item[(a)] $R_4=\{x\in H_4: |\lambda^{HX_5S_1}(x)|=\lambda^{HX_5S_1}(1)\}$ and $S_2 =Stab_{S_1}(\lambda|_{H_6H_5H_4}).$ Moreover, $\lambda^{HX_5S_1}|_{R_4}=\lambda^{HX_5S_1}(1)\phi_{B_4}.$

\item[(b)] If $B_4\neq 0,$ then $Stab_{S_1}(\lambda)=\{1\}.$ Hence, if $\eta$ is an extension of $\lambda$ to $HX_5,$ then $I_{HX_5S_1}(\eta)=HX_5.$  Furthermore, if  $\eta,\eta'$ are two extensions of $\lambda|_{H_6H_5H_4H_3}$ to $HX_5,$ then $\eta^{HX_5S_1}=\eta'^{HX_5S_1}$ iff $B_2(\eta)=B_2(\eta')$ and $\eta|_{X_5}=\eta'|_{X_5}.$

\item[(c)]  If $B_4=0,$ then there exists $x\in S_1$ such that ${}^{x}\lambda|_{X_i}=1_{X_i}$ for all $X_i\subset H_5H_4.$  Furthermore, $\overline{H_5H_4}\subset ker({}^x\lambda)^{HX_5S_2}$ and the induction map from $Irr(HX_5S_2,{}^x\lambda)$ to $Irr(HX_5S_1,\lambda)$ is bijective.
\end{itemize}
\end{lemma}
\noindent{\it Proof.} See Subsection \ref{proof-lem:5ext2}$.~\Box$

\medskip
The main idea of Lemma \ref{lem:5ext2} (a) can be visualized in the following figure.
\begin{center}
\begin{picture}(300,115)
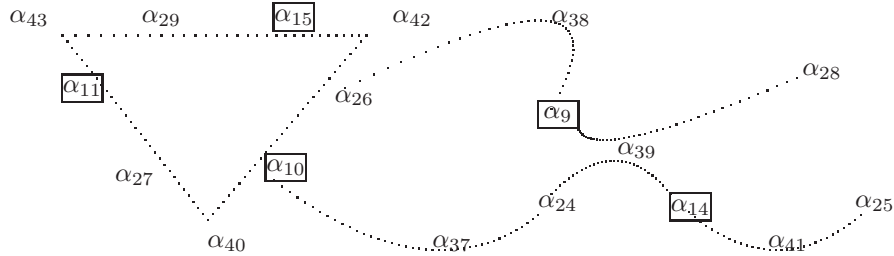
\put(10,-3){Figure $UE_8(q)$:  Relations of between root heights 4 in $H$ and 2 in $T$. }
\put(0,100){$\alpha_{43}$}
\put(50,100){$\alpha_{29}$}
\put(100,97){\framebox(15,10){$\alpha_{15}$}}
\put(145,100){$\alpha_{42}$}
\put(20,70){\framebox(15,10){$\alpha_{11}$}}
\put(40,40){$\alpha_{27}$}
\put(75,15){$\alpha_{40}$}
\put(97,40){\framebox(15,10){$\alpha_{10}$}}
\put(123,70){$\alpha_{26}$}
%%\put(80,72){\circle{15}}
%%\put(75,70){$\alpha_3$}
\qbezier[30](20,95)(95,95)(135,95)
\qbezier[30](135,95)(105,60)(75,25)
\qbezier[30](75,25)(48,60)(20,95)

\put(160,15){$\alpha_{37}$}
\put(287,15){$\alpha_{41}$}
\put(205,100){$\alpha_{38}$}
\put(230,50){$\alpha_{39}$}

\put(200,60){\framebox(15,10){$\alpha_9$}}
\put(300,80){$\alpha_{28}$}
\qbezier[40](215,60)(220,45)(298,79)
\qbezier[40](127,75)(240,130)(205,67)

\put(200,30){$\alpha_{24}$}
\qbezier[40](100,40)(160,-5)(200,27)
\put(250,25){\framebox(15,10){$\alpha_{14}$}}
\qbezier[30](205,35)(230,60)(250,35)
\put(320,30){$\alpha_{25}$}
\qbezier[30](257,28)(290,0)(322,27)

\end{picture}
\end{center}

\medskip
Recall that we have $q^{12}$ linear characters $\lambda$ of $H$ such that $\lambda|_{X_i}=\phi$ for all $X_i\subset Z(U)$ and $\lambda|_{X_i}=1_{X_i}$ for all $X_i\subset H_5.$ In these, there are $q^{11}$ linears with $B_4=0$ and $q^{11}(q-1)$ linears with $B_4\neq 0.$ By Lemma \ref{lem:5ext2} (a), it is clear that $B_4=B_4(\lambda)$ is invariant under the action of $S_1.$ Therefore, by Lemma \ref{lem:5ext2} (b), these $q^{11}(q-1)$ linears with $B_4\neq 0$ extend to $HX_5$ and induce irreducibly to $HX_5S_1.$ Thus, we obtain $\frac{q^{12}(q-1)}{|S_1|}=q^2(q-1)$ irreducible characters of $HX_4S_1$ of degree $|S_1|=q^{10}$ which are parametrized by $(b_5,B_2,B_4)$ where  $b_5,B_2\in \F_q$ and $B_4\in \F_q^\times.$ By Lemma \ref{lem:5ext1} (c), we obtain $q^2(q-1)$ characters of $U$ of degree $q^{16}$ which can be denoted by $\chi_{q^{16}}^{b_5,B_2,B_4}.$

Since $[HX_5S_1:HX_5S_2]=5,$ by Lemma \ref{lem:5ext2}(c), the constituents of $q^{11}$ linears of $H$ with $B_4=0$ inducing to $HX_5S_1$ correspond to the ones of $q^{6}$ linears $\lambda$ of $H$ inducing to $HX_5S_2,$ where $\lambda|_{X_i}=\phi$ for all $X_i\subset H_6,$ $\lambda|_{X_i}=1_{X_i}$ for all $X_i\subset H_5H_4,$ and $\lambda|_{X_i}=\phi_{b_i}$ for all $X_i\subset H_3H_2$ where $b_i\in\F_q.$ Let $\lambda\in Irr(H)$ be one of above $q^6$ linears of $H$. %such that $\lambda|_{X_i}=\phi$ for all $X_i\subset H_6,$ $\lambda|_{X_i}=1_{X_i}$ for all $X_i\subset H_5H_4,$ and $\lambda|_{X_i}=\phi_{b_i}$ for all $X_i\subset H_3H_2$ where $b_i\in\F_q.$
 Now we consider $\lambda$ in the scenario of the subgroup $HX_5S_2.$ Let $\overline{H_5H_4H_3}$ be the normal closure of $H_5H_4H_3$ in $HX_5S_3.$

%(Work with $HX_5S_2/\overline{H_5H_4}$)

\begin{lemma} \label{lem:5ext3}
The following are true
\begin{itemize}
\item[(a)] $R_3=\{x\in H_3: |\lambda^{HX_5S_2}(x)|=\lambda^{HX_5S_2}(1)\}$ and $S_3=Stab_{S_2}(\lambda|_{H_6H_5H_4H_3}).$ Moreover, $\lambda^{HX_5S_2}|_{R_3}=\lambda^{HX_5S_2}(1)\phi_{B_3}.$

\item[(b)] If $B_3\neq 0,$ then $Stab_{S_2}(\lambda)=\{1\}.$ Hence, if $\eta$ is an extension of $\lambda$ to $HX_5,$ then $I_{HX_5S_2}(\eta)=HX_5.$ Furthermore, if $\eta,\eta'$ are two extensions of $\lambda$ to $HX_5,$ then $\eta^{HX_5S_2}=\eta'^{HX_5S_2}$ iff $\eta|_{X_5}=\eta'|_{X_5}.$

\item[(c)]  If $B_3=0,$ then there exists $x\in S_2$ such that ${}^{x}\lambda|_{X_i}=1_{X_i}$ for all $X_i\subset H_5H_4H_3.$ Furthermore, $\overline{H_5H_4H_3}\subset ker({}^x\lambda)^{HX_5S_3}$ and the induction map from $Irr(HX_5S_3,{}^x\lambda)$ to $Irr(HX_5S_2,\lambda)$ is bijective.
\end{itemize}
\end{lemma}
\noindent{\it Proof.} See Subsection \ref{proof-lem:5ext3}$.~\Box$

\medskip
The main idea of Lemma \ref{lem:5ext3} (a) can be described as follows.
\begin{center}
\begin{picture}(300,80)
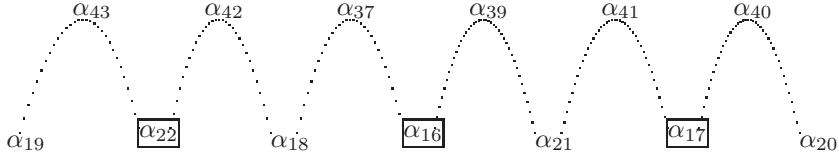
\put(5,0){Figure $UE_8(q)$:  Relations of between root heights 3 in $H$ and 3 in $T$. }

\put(25,70){$\alpha_{43}$}
\put(75,70){$\alpha_{42}$}
\put(125,70){$\alpha_{37}$}
\put(175,70){$\alpha_{39}$}
\put(225,70){$\alpha_{41}$}
\put(275,70){$\alpha_{40}$}

\put(0,20){$\alpha_{19}$}
\put(50,20){\framebox(15,10){$\alpha_{22}$}}
\put(100,20){$\alpha_{18}$}
\put(150,20){\framebox(15,10){$\alpha_{16}$}}
\put(200,20){$\alpha_{21}$}
\put(250,20){\framebox(15,10){$\alpha_{17}$}}
\put(300,20){$\alpha_{20}$}

\qbezier[30](5,25)(30,110)(50,27)
\qbezier[30](62,27)(80,110)(100,25)
\qbezier[30](110,25)(130,110)(150,27)
\qbezier[40](162,27)(180,110)(200,25)
\qbezier[40](210,25)(230,110)(250,27)
\qbezier[40](262,27)(280,110)(300,25)

\end{picture}
\end{center}

\medskip
Recall that we have $q^{6}$ linear characters $\lambda$ of $H$ such that $\lambda|_{X_i}=\phi$ for all $X_i\subset Z(U)$ and $\lambda|_{X_i}=1_{X_i}$ for all $X_i\subset H_5H_4.$ In these, there are $q^{5}$ linears with $B_3=0$ and $q^{5}(q-1)$ linears with $B_3\neq 0.$ By Lemma \ref{lem:5ext3} (a), it is clear that $B_3=B_3(\lambda)$ is invariant under the action of $S_2.$ Therefore, by Lemma \ref{lem:5ext3} (b), these $q^{5}(q-1)$ linears with $B_3\neq 0$ extend to $HX_5$ and induce irreducibly to $HX_5S_2.$ Thus, we obtain $\frac{q^{6}(q-1)}{|S_2|}=q(q-1)$ irreducible characters of $HX_4S_2$ of degree $|S_2|=q^{5}$ which are parametrized by $(b_5,B_3)$ where  $b_5\in \F_q$ and $B_3\in \F_q^\times.$ By Lemma \ref{lem:5ext2} (c) and Lemma \ref{lem:5ext1} (c),  we obtain $q(q-1)$ characters of $U$ of degree $q^{16}$ which can be denoted by $\chi_{q^{16}}^{b_5,B_3}.$

Since $[HX_5S_2:HX_5S_3]=3,$ by Lemma \ref{lem:5ext3}(c), the constituents of $q^{5}$ linears of $H$ with $B_3=0$ inducing to $HX_5S_2$ correspond to the ones of $q^{2}$ linears $\lambda$ of $H$ inducing to $HX_5S_3,$ where $\lambda|_{X_i}=\phi$ for all $X_i\subset H_6,$ $\lambda|_{X_i}=1_{X_i}$ for all $X_i\subset H_5H_4H_3,$ and $\lambda|_{X_i}=\phi_{b_i}$ for all $X_i\subset H_2$ where $b_i\in\F_q.$ Let $\lambda\in Irr(H)$ be one of above $q^2$ linears of $H$. % such that $\lambda|_{X_i}=\phi$ for all $X_i\subset H_6,$ $\lambda|_{X_i}=1_{X_i}$ for all $X_i\subset H_5H_4H_3,$ and $\lambda|_{X_i}=\phi_{b_i}$ for all $X_i\subset H_2$ where $b_i\in\F_q.$
Now we consider $\lambda$ in the scenario of the subgroup $HX_5S_3.$
%Let $\overline{H_5H_4H_3H_2}$ be the normal closure of $H_5H_4H_3H_2$ in $HX_5S_4.$

%(Work with $HX_5S_3/\overline{H_5H_4H_3}$)

\begin{lemma} \label{lem:5ext4}
The following are true.
\begin{itemize}
\item[(a)] $R_2=\{x\in H_2: |\lambda^{HX_5S_3}(x)|=\lambda^{HX_5S_3}(1)\}$ and $S_4=Stab_{S_3}(\lambda).$ Moreover, $\lambda^{HX_5S_3}|_{R_2}=\lambda^{HX_5S_3}(1)\phi_{B_2}.$

\item[(b)] If $B_2\not\in \{c^4:c\in \F_q^\times\}$ and let $\eta$ be an extension of $\lambda$ to $HX_5,$ then $I_{HX_5S_3}(\eta)=HX_5.$ Therefore, $S_4$ acts transitively and faithfully on all extensions of $\lambda$ to $HX_5.$

\item[(c)]  If $B_2=c^4\in \F_q^\times ,$  then $\lambda$ extends to $HX_5F_4$ and $HF_5S_4.$ Let $\lambda_1,\lambda_2$ be two extensions of $\lambda$ to $HX_5F_4.$ Then $I_{HX_5S_3}(\lambda_1)=HX_5F_4.$ Moreover, ${\lambda_1}^{HX_5S_3}={\lambda_2}^{HX_5S_3}$ iff $\lambda_1|_{F_4}=\lambda_2|_{F_4}$ and $\lambda_1|_{F_5}=\lambda_2|_{F_5}.$
\end{itemize}
\end{lemma}
\noindent{\it Proof.} See Subsection \ref{proof-lem:5ext4}$.~\Box$

\medskip
It is noted that $Stab_{S_3}(\lambda)=Stab_{T}(\lambda).$ The main idea of Lemma \ref{lem:5ext4} (a) can be visualized in the following figure.
\begin{center}
\begin{picture}(300,80)
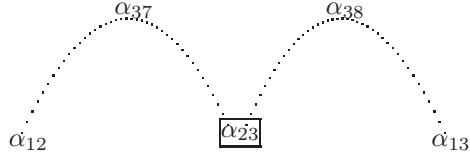
\put(-10,0){Figure $UE_8(q)$:  Relations of between root heights 2 in $H$ and 4 in $T$. }

\put(90,70){$\alpha_{37}$}
\put(170,70){$\alpha_{38}$}

\put(50,20){$\alpha_{12}$}
\put(130,20){\framebox(15,10){$\alpha_{23}$}}
\put(210,20){$\alpha_{13}$}

\qbezier[40](55,25)(95,110)(133,28)
\qbezier[40](140,28)(175,110)(215,25)

\end{picture}
\end{center}

\medskip
Recall that we have $q^{2}$ linear characters $\lambda$ of $H$ such that $\lambda|_{X_i}=\phi$ for all $X_i\subset Z(U)$ and $\lambda|_{X_i}=1_{X_i}$ for all $X_i\subset H_5H_4H_3.$ By Lemma \ref{lem:5ext4} (a), it is clear that $B_2=B_2(\lambda)$ is invariant under the action of $S_3.$ Since $F_q^\times$ is cyclic, we have $|\{c^4:c\in \F_q^\times\}|=\frac{q-1}{4}.$ Therefore, there are $\frac{q(q-1)}{4}$ linears with $B_2=c^4\in \F_q^\times,$ and there are $\frac{3q(q-1)}{4}$ linears with $B_2\not\in \{c^4: c\in \F_q^\times\}.$ Thus, these linears with $B_2\not\in \{c^4: c\in \F_q^\times\}$ extend to $HX_5$ and induce irreducibly to $HX_5S_3$ of degree $|S_3|=q^2.$ By Lemma \ref{lem:5ext3} (c), Lemma \ref{lem:5ext2} (c) and Lemma \ref{lem:5ext1} (c),  we obtain $\frac{3(q-1)}{4}$ characters of $U$ of degree $q^{16}$ which can be denoted by $\chi_{q^{16}}^{B_2}.$

Now we sum up all irreducibles of degree $q^{16}$ as counted above and denoted by $\chi_{q^{16}}^{b_5,B_2,B_3,B_5},$ $\chi_{q^{16}}^{b_5,B_2,B_4},$ $\chi_{q^{16}}^{b_5,B_3},$ and $\chi_{q^{16}}^{B_2}.$ Therefore, $\FF_8$ contains exactly $q^3(q-1)+q^2(q-1)+q(q-1)+\frac{3(q-1)}{5}$ characters $\chi$ of $U$ of degree $q^{16}$ such that $\chi|_{X_i}=\chi(1)\phi$ for all $X_i\subset Z(U).$

By Lemma \ref{lem:5ext4} (c),  let $\lambda_1$ be an extension of $\lambda$ to $HX_5F_4,$ then ${\lambda_1}^{HX_5S_3}$ is irreducible of degree $[HX_5S_3:HX_5F_4]=\frac{q^2}{5}.$ These ${\lambda_1}^{HX_5S_3}$ only depend on $B_2$ and their restrictions to $F_4,F_5.$ Hence, Lemma \ref{lem:5ext3} (c), Lemma \ref{lem:5ext2} (c) and Lemma \ref{lem:5ext1} (c), ${\lambda_1}^U\in Irr(U)$ of degree $\frac{q^{16}}{5}$ is denoted by $\chi_{\frac{q^{16}}{5}}^{b_4,b_5,B_2}$ where $b_4,b_5\in \F_5$ and $B_2\in\{c^4:c\in \F_q^\times\}.$ Therefore, $\FF_8$ has exactly $\frac{25(q-1)}{4}$ irreducibles of degree $\frac{q^{16}}{5}$ such that $\chi|_{X_i}=\chi(1)\phi$ for all $X_i\subset Z(U).$

By the transitivity of the conjugate action of the maximal split torus $T_0$ of the Chevalley group $E_8(q)$ on $\oplus_{i=37}^{43}Irr(X_i)^\times,$ there are $(q-1)^8(q^3+q^2+q+\frac{3}{4})$ characters $\chi\in \FF_8$ of degree $q^{16},$ and $\frac{25(q-1)^8}{4}$ characters $\chi\in \FF_8$ of degree $\frac{q^{16}}{5}$ such that $\chi|_{X_i}=\chi(1)\phi_{a_i},$ where $a_i\in\F_q^\times,37\leq i\leq 43.$  The following diagram summarizes all the above arguments with their assumptions.

\begin{center}
\begin{picture}(370,105)
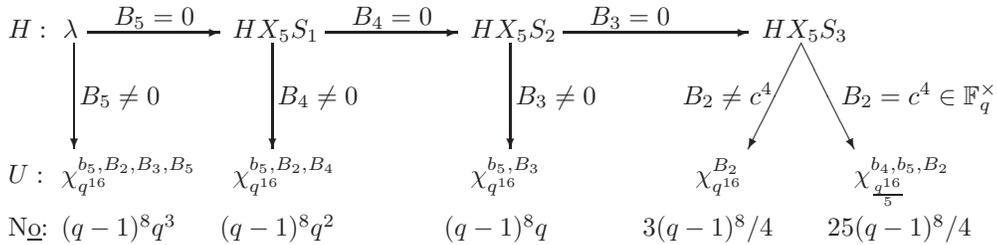
\put(70,0){Figure $UE_8(q)$:  Summary on the branching rules of $\lambda$. }
\put(0,35){$U:$}
\put(0,15){N\underline{o}:}

\put(0,90){$H:~\lambda$}
\put(30,92){\vector(1,0){50}} \put(40,95){$B_5=0$}
\put(25,88){\vector(0,-1){40}} \put(27,65){$B_5\neq 0$}
\put(20,35){$\chi_{q^{16}}^{b_5,B_2,B_3,B_5}$}
\put(20,15){$(q-1)^8q^3$}

\put(85,90){$HX_5S_1$}
\put(120,92){\vector(1,0){50}} \put(130,95){$B_4=0$}
\put(100,88){\vector(0,-1){40}} \put(102,65){$B_4\neq 0$}
\put(85,35){$\chi^{b_5,B_2,B_4}_{q^{16}}$}
\put(80,15){$(q-1)^8q^2$}

\put(175,90){$HX_5S_2$}
\put(210,92){\vector(1,0){70}} \put(220,95){$B_3=0$}
\put(190,88){\vector(0,-1){40}} \put(192,65){$B_3\neq 0$}
\put(175,35){$\chi^{b_5,B_3}_{q^{16}}$}
\put(165,15){$(q-1)^8q$}

\put(285,90){$HX_5S_3$}
\put(300,88){\vector(-1,-2){20}} \put(255,65){$B_2\neq c^4$}
\put(300,88){\vector(1,-2){20}} \put(315,65){$B_2=c^4\in\F_q^\times$}
\put(260,35){$\chi_{q^{16}}^{B_2}$}
\put(240,15){$3(q-1)^8/4$}
\put(320,35){$\chi^{b_4,b_5,B_2}_{\frac{q^{16}}{5}}$}
\put(310,15){$25(q-1)^8/4$}

\end{picture}
\end{center}

This gives the proof for the next theorem.

\begin{theorem}
Let $\chi\in \FF_8.$ The following are true.
\begin{itemize}
\item[(a)] If $\chi(1)=q^{16},$ then there exists $t\in T_0$ such that ${}^t\chi$ is an element of $ \{\chi^{b_5,B_2,B_3}_{q^{16}},\chi^{b_5,B_2}_{q^{16}},\chi^{b_5}_{q^{16}},\chi^{B_2}_{q^{16}}\}.$

\item[(b)] If $\chi(1)=q^{16}/5,$ then there exists $t\in T_0$ such that ${}^t\chi=\chi^{b_4,b_5,B_2}_{\frac{q^{16}}{5}}.$

\end{itemize}
\end{theorem}

%%%%%%%%%-------------------------------------%%%%%%%%%%%%%%%%%%%%
%%%%%%%%%-------------------------------------%%%%%%%%%%%%%%%%%%%%
%%%%%%%%%-------------------------------------%%%%%%%%%%%%%%%%%%%%
%%%%%%%%%-------------------------------------%%%%%%%%%%%%%%%%%%%%

\section{All proofs.} \label{Section-All-proofs}

In all proofs, we use the following technique:
\begin{itemize}
\item[(a)] For all the decomposition of the commutator formula into  product, we apply the formula $[a,bc]=[a,c][a,b]^c.$

\item[(b)] For $ H\leq G$ and  $L\unlhd G,$ for each $\lambda\in Irr(L),$ $Stab_G(\lambda):=\{x\in G:{}^x\lambda=\lambda\},$ and $Stab_{G}(\lambda)\subset Stab_G(\lambda|_H)=:K,$ hence, $Stab_G(\lambda)=Stab_{K}(\lambda).$

\item[(c)] For $K\leq G$ and $H\unlhd G,$ to extend a linear character $\lambda$ of $H$ to $HK$, we check if $[HK,HK]\subset ker(\lambda).$
\end{itemize}

%%%%%%%%%%%%%%%%%%%%%% PROOF OF FIELDS' PROPERTY %%%%%%%%%%%%%%%%%%%%%%

\subsection{Proof of Proposition \ref{prop:mainFp}}
\label{proof-prop:mainFp}
Let $a\in\F_q^\times.$ We are going to prove \ref{prop:mainFp} (c).

%(a) Consider $t$ as a variable, it is clear that $ca$ is a solution of $t^p-a^{p-1}t$ for all $c\in \F_p.$ Since $a\in\F_q^\times,$ all $\{ca:c\in\F_p\}$ are distinct. Therefore, the claim holds by the degree of two polynomials $t^p-a^{p-1}t$ and $\prod_{c\in \F_p}(t-ca).$

% \medskip
% (b) For each $u,v\in\F_q,$ we have
% \begin{center} $\begin{array}{ll} (u^p-a^{p-1}u)+(v^p-a^{p-1}v)&=(u^p+v^p)-a^{p-1}(u+v) \\ &=(u+v)^p-a^{p-1}(u+v)\in\T_a. \end{array}$ \end{center}
% Hence, $\T_a\leq (\F_q,+).$
%
% By (a), all $t\in\{u+ca:c\in\F_p\}$ give the same value in $\T_a.$ Therefore, $|\T_a|\leq |\F_q/\F_p|=\frac{q}{p}.$ For each $b\in\T_a,$ the equation $t^p-a^{p-1}t=b$ has at most $p$ solutions. Hence, $|\T_a|=\frac{q}{p}.$

Since $\gcd(q-1,p)=1,$ for each $b\in\F_q^\times,$ there is $s\in\F_q^\times$ such that $b=s^p.$ We have
%\begin{center}
$b\prod_{c\in\F_p}(t-ca)=s^p\prod_{c\in\F_p}(t-ca)=\prod_{c\in\F_p}(st-csa)\in\T_{sa}.$
%\end{center}
Hence, $\F_q^\times$ acts on $\{\T_a:a\in\F_q^\times\}$ such that $s^p\T_a=\T_{sa}.$ We claim that $b\T_a=\T_a$ iff $b\in\F_p^\times.$

By (a), it is clear if $b\in\F_p^\times.$ Suppose there exists $b\in\F_q-\F_p$ such that $b\T_a=\T_a.$ Since $b\in \F_q-\F_p,$ we have $b=s^p$ for some $s\in\F_q-\F_p.$ Hence, $s^{p-1}\neq 1.$ By induction, we have $\T_a=b\T_a=b^2\T_a=...=b^k\T_a$ for all $k\in\N.$ Therefore, for each $t\in\F_q,$ $t^p-(as^k)^{p-1}t\in\T_a$ for all $k\in\N.$ For each $l\in\N^+,$ we have
\begin{center}
 $\begin{array}{ll}
a^{p-1}t(s^{p-1}-1)^l&=a^{p-1}t\sum_{k=0}^l \frac{l!}{k!(l-k)!}s^{(p-1)k}(-1)^{l-k}-t^p(1-1)^l
\\
&=a^{p-1}t\sum_{k=0}^l \frac{l!}{k!(l-k)!} s^{(p-1)k}(-1)^{l-k}-t^p\sum_{k=0}^l \frac{l!}{k!(l-k)!} (-1)^{l-k}
\\
&=-\sum_{k=0}^l(-1)^{l-k} \frac{l!}{k!(l-k)!} (t^p-(as^k)^{p-1}t)\in\T_a.
\end{array}$
\end{center}

Since $s^{p-1}-1\in\F_q^\times,$ there exists $l\in\N^\times$ such that $(s^{p-1}-1)^l=1.$ Therefore, $a^{p-1}t\in\T_a$ for all $t\in\F_q.$ We have $(t^p-a^{p-1}t)+a^{p-1}t=t^p\in\T_a.$ This makes a contradiction since $\{t^p:t\in \F_q\}=\F_q\gvertneqq \T_a.$ So the claim holds.
Thus $\F_q^\times/\F_p^\times$ acts faithfully and transitively on $\{\T_a:a\in\F_q^\times\},$ and  $|\{\T_a:a\in\F_q\}|=\frac{q-1}{p-1}.$

It is easy to see that $ker\phi_a=ker\phi_{ca}$ for all $c\in\F_p^\times$ since $\F_p\cong \Z_p$ and  $\phi_a(ku)=\phi_a(u)^k$ for all $k\in\N.$ Therefore, $|\{ker\phi_a:a\in \F_q^\times\}|=\frac{q-1}{p-1}.$

Since $\{ker\phi_a : a\in \F_q^\times\}$ are all subgroups of index $p$ in $\F_q,$ $\{ker\phi_a : a\in \F_q^\times\}=\{\T_a:a\in \F_q^\times\}.$ Therefore, for each $a\in\F_q^\times,$ there exists $b\in\F_q^\times$ such that $b\T_a=ker\phi,$ and $cb\T_a=ker\phi$ iff $c\in\F_p^\times.~\Box$

%%%%%%%%%%%%%%%%%%%%%%%%%%%%%%%%%%%%%% PROOF OF REDUCTION LEMMA %%%%%%%%%%%%%%%%%%%%%%%%%%

\subsection{Proof of Proposition \ref{lem:Reduction}}
\label{proof-lem:Reduction}
(a) Suppose $\chi\in Irr(N/Y,\lambda),$ we are going to show that $\chi^G\in Irr(G)$ by showing that the inertia group $I_G(\chi)=N.$

Since $Y\subset ker(\chi)$ and $Z\subset Z(N),$ we have $\chi|_{ZY}=\chi(1) \lambda.$ Since $X\subset N_G(ZY),$ for each $x\in X,$ ${}^x\lambda\in Irr(ZY).$ Hence, for any $u\neq v\in X$ we have
\begin{center}
${}^u\chi|_{ZY}=\chi(1)\,{}^u\lambda\neq \chi(1)\,{}^v\lambda ={}^v\chi|_{ZY},$ i.e. ${}^u\chi\neq {}^v\chi.$
\end{center}
Therefore, $x\in X$ such that ${}^x\chi=\chi$ iff $x=1.$ Since $X$ is a transversal of $N$ in $G,$ this shows that the inertia group $I_G(\chi)=N.$

The above argument also proves that for $\chi_1,\chi_2\in Irr(N/Y,\lambda)$ and $u\neq v\in X$ we have
%${}^u\chi_1|_{ZY}=\chi_1(1){}^u\lambda\neq \chi_2(1){}^v\lambda ={}^v\chi_2|_{ZY},$ i.e.
${}^u\chi_1\neq {}^v\chi_2.$  So by the Mackey formula for the double coset $N\backslash G/N=G/N$ represented by $X,$ we have
\[
({\chi_1}^G,{\chi_2}^G)=({\chi_1}^G|_N,\chi_2)=\sum_{x\in X}({}^x\chi_1,\chi_2)=(\chi_1,\chi_2)
=\left\{\begin{array}{ll} 1&\mbox{ if } \chi_1=\chi_2\\ 0 &\mbox{otherwise}\end{array}\right.
\]

%\medskip
(b) It is enough to show that the induction map is surjective, i.e. for each $\xi\in Irr(G,\lambda)$ there exists $\chi\in Irr(N/Y,\lambda)$ such that $\xi=\chi^G.$

Suppose $\xi|_N=\sum_{\chi_i\in S}a_i\chi_i$ where $a_i\in \N^\times$ and $S\subset Irr(H).$ By Frobenious reciprocity, $0\neq (\xi,\lambda^G)=(\xi|_{ZY},\lambda),$ there exists at least a constituent $\chi_0$ of $\xi|_N$ such that $(\chi_0|_{ZY},\lambda)\neq 0,$ i.e. $\chi_0\in Irr(N,\lambda).$

Since $\lambda|_Y=\lambda(1) 1_Y$ and $(\chi_0|_Y,\lambda|_Y)\geq(\chi_0|_{ZY},\lambda)>0,$ we have $\chi_0$ is a constituent of ${1_Y}^N.$ Since $Y\unlhd N,$ all constituents of ${1_Y}^N$ are $Irr(N/Y).$ Therefore, $\chi_0\in Irr(N/Y,\lambda).$
By (a), ${\chi_0}^G\in Irr (G),$ hence it forces $\xi={\chi_0}^G.~\Box$

%************************************************************************************

\subsection{Proofs of Sylow 2-subgroups of $D_4(2^f)$}

%\bigskip
\subsubsection{Proof of Lemma \ref{lem:2extensions}}
\label{proof-lem:2extensions}
Set $\lambda=\lambda_{b_3,b_5,b_6,b_7}^{a_8,a_9,a_{10}}$ for the whole proof.

(a)First we show that $Stab_{T}(\lambda)=S_{124}.$ Since ${}^y\lambda(x)=\lambda(x)$ iff $\lambda(x^{-1}x^y)=\lambda([x,y])=1$ and $X_8X_9X_{10}\subset Z(U),$ it suffices to check for $[X_5X_6X_7,T].$ For all $t_i,s_j\in \F_q,$ we have
%
%$[x_5x_6x_7,x_1x_2x_4]=[x_6,x_1][x_5,x_2][x_7,x_1][x_5,x_4][x_7,x_2][x_6,x_4].$
\[[x_5(t_5)x_6(t_6)x_7(t_7),x_1(s_1)x_2(s_2)x_4(s_4)]=x_8(t_6t_1+t_5t_2)x_9(t_7t_1+t_5t_4)x_{10}(t_7t_2+t_6t_4).\]

Therefore, $x_1(s_1)x_2(s_2)x_4(s_4)\in Stab_{T}(\lambda)$ iff for all $t_5,t_6,t_7\in\F_q,$ %\[\phi(a_8(t_6s_1+t_5s_2)+a_9(t_7s_1+t_5s_4)+a_{10}(t_7s_2+t_6s_4))=1.\]

%We have
\begin{center}
$\begin{array}{ll}
1&=\phi(a_8(t_6s_1+t_5s_2)+a_9(t_7s_1+t_5s_4)+a_{10}(t_7s_2+t_6s_4))\\
&=\phi(t_5(a_8s_2+a_9s_4)+t_6(a_8s_1+a_{10}s_4)+t_7(a_9s_1+a_{10}s_2))
\end{array}$
\end{center}
iff $a_8s_2+a_9s_4=a_8s_1+a_{10}s_4=a_9s_1+a_{10}s_2=0,$ i.e. $\frac{s_1}{a_{10}}=\frac{s_2}{a_9}=\frac{s_4}{a_8}.$ So $Stab_{T}(\lambda)=S_{124}.$

To find all scalar points of $X_5X_6X_7$ on $\lam^U,$ since $X_3T$ is a transversal of $H$ in $U$ and $[X_3,X_5X_6X_7]=\{1\},$ it is enough to find ones of $X_5X_6X_7$ on $T,$ i.e. find $x_5x_6x_7\in X_5X_6X_7$ such that $\lam([x_5x_6x_7,x_1x_2x_4])=1$ for all $x_1x_2x_4\in T.$  Use above computation, for all $s_i\in \F_q,$ we need
\begin{center}
$\begin{array}{ll}
1&=\phi(a_8(t_6s_1+t_5s_2)+a_9(t_7s_1+t_5s_4)+a_{10}(t_7s_2+t_6s_4))\\
&=\phi(s_1(a_8t_6 +a_9t_7)+s_2(a_8t_5+a_{10}t_7)+s_4(a_9t_5+a_{10}t_6))\end{array}$
\end{center}
iff $a_8t_6 +a_9t_7=a_8t_5+a_{10}t_7=a_9t_5+a_{10}t_6=0,$ i.e. $\frac{t_5}{a_{10}}=\frac{t_6}{a_9}=\frac{t_7}{a_8}.$ Hence, $\prod_{i=5}^7x_i(t_i)=x_{567}(\frac{t_7}{a_8})\in S_{567}.$ So $S_{567}=\{x\in X_5X_6X_7:|\lambda^U(x)|=\lambda^U(1)\}.$

Now, to prove that $\lambda^U|_{S_{567}}=q^4\phi_{At_0},$ it suffices to check that $\lambda(x_{567}(t))=\phi_{At_0}(t).$ For each  $x_{567}(t)=x_5(a_{10}t)x_6(a_9t)x_7(a_8t)\in S_{567},$ we have
\[\lambda(x_{567}(t))=\phi(t(b_5a_{10}+b_6a_9+b_7a_8))=\phi(tAt_0)=\phi_{At_0}(t).\]

\medskip
(b)  We study $Irr(U,\lambda)$ by two following ways. Let $K_1:=H X_3 F_{124}$ and  $K_2:=H S_{124}F_3.$ Since $H=[U,U],$ it is clear that $H_1,K_1\unlhd U.$

\[
\begin{array}{rcl}
&H&\\
\swarrow &&\searrow \\
HX_3&&HS_{124}\\
\downarrow \hspace{0.4cm}&&\hspace{0.4cm}\downarrow\\
K_1=H X_3 F_{124}&&K_2=H S_{124} F_3\\
 \searrow  &&\swarrow \\
& U &
\end{array}
\]

Since $HX_3$ is abelian, $\lambda$ extends to $HX_3$ as $\eta_1.$
By (a), $S_{124}=Stab_{T}(\lambda),$ for all $x\in H,x_{124}\in S_{124},$ $\lambda([x,x_{124}])=1,$ hence $\lambda$ extends to $HS_{124}$ as $\eta_2.$ To show that $\lambda$ extends to $K_1$ and $K_2,$ we prove that $[K_1,K_1]\subset ker(\eta_1), [K_2,K_2]\subset ker(\eta_2).$
We have %By using $[a,bc]=[a,c][a,b]^c,$
\\
$[x_3(t_3),x_1(s_1)x_2(s_2)x_4(s_4)]
%=[x_3(t_3),x_4(s_4)][x_3(t_3),x_1(a_{10}t)x_2(a_9t)]^{x_4(t_4)}
%=x_7(t_3t_4)([x_3(t_3),x_2(t_2)][x_3(t_3),x_1(t_1)]^{x_2(t_2)})^{x_4(t_4)}=x_7(t_3t_4)(x_6(t_2t_3)x_5(t_1t_3)x_8(t_1t_2t_3))^{x_4(t_4)}
=x_5(s_1t_3)x_6(s_2t_3)x_7(s_4t_3)x_8(s_1s_2t_3)x_9(s_1s_4t_3)x_{10}(s_2s_4t_3),$ and $\lambda(x_5(s_1t_3)x_6(s_2t_3)x_7(s_4t_3)x_8(s_1s_2t_3)x_9(s_1s_4t_3)x_{10}(s_2s_4t_3))$
\\
$=\phi(t_3(b_5s_1+b_6s_2+b_7s_4+a_8s_1s_2+a_9s_1s_4+a_{10}s_2s_4))=(*).$

Plug $s_1=a_{10}t,s_2=a_9t,s_4=a_8t$ into $(*),$ we have

$(*)=\phi(t_3(t(b_5a_{10}+b_6a_9+b_7a_8)+t^2a_8a_9a_{10}))=\phi(t_3At(t_0+t)).$

Now we divide into two cases where $t_0=0$ and $t_0\neq 0.$ First, if $t_0=0,$ $\phi(t_3At^2)=1$ for all $t_3$ iff $t=0,$ hence, $Stab_{T}(\eta_1)=\{1\}=F_{124},$ i.e. $I_U(\eta_1)=HX_3.$ And $\phi(t_3At^2)=0$ for all $t$ iff $t_3=0,$ hence, $Stab_{X_3}(\eta_2)=\{1\}=F_3,$ i.e. $I_U(\eta_2)=HS_{124}.$

If $t_0\neq 0,$ then $\phi(t_3At(t_0+t))=1$ for all $t_3$ iff $t\in\{0,t_0\}.$ Therefore, $[K_1,K_1]\subset ker\lambda.$ For each $\eta\in Irr(HX_3,\lambda),$ $Stab_{T}(\eta)=\{1,x_{124}(t_0)\}=F_{124}.$

We have $\phi(t_3At(t_0+t))=1$ for all $t$ iff  $t_3\in\{0,\frac{(t_0)_\phi}{A}\},$ by Proposition \ref{prop:mainFp}. Hence, $[K_2,K_2]\subset ker(\lambda).$ For each $\gamma\in Irr(HS_{124},\lambda),$ $Stab_{X_3}(\gamma)=\{1,x_3(\frac{(t_0)_\phi}{A})\}=F_3.$

So $\lambda$ extends to $K_1$ and $K_2.$ For each $\lambda_i\in Irr(K_i,\lambda),$ $I_U(\lambda_i)=K_i,~i=1,2.$

\medskip
(c) Let $\lambda_1,\lambda_2$ be extensions of $\lambda$ to $K_1.$ Let $\eta$ be an extension of $\lambda$ to $K_2.$ By (b), we have ${\lambda_1}^U,{\lambda_2}^U,\eta^U\in Irr(U,\lambda).$

We choose $1\in S\subset T$ as a representative set of the double coset $K_1\backslash U/K_2,$ by Mackey formula, since $K_1\cap K_2=HF_3F_{124}$ and $K_1\unlhd U,$ we have
\begin{center}
$\begin{array}{ll}
({\lambda_1}^U, {\eta}^U)&=\sum_{s\in S}({}^s\lambda_1|_{{}^sK_1\cap K_2},\eta|_{{}^sK_1\cap K_2})
\\
&=\sum_{s\in S}({}^s\lambda_1|_{HF_3F_{124}},\eta|_{HF_3F_{124}}).
\end{array}$
\end{center}

For each $s\in S,$ if ${}^s\lambda_1|_{HF_3F_{124}}=\eta|_{HF_3F_{124}},$ then ${}^s\lambda_1|_H=\eta|_H.$ Since both are extensions of $\lambda$ from $H,$ we have ${}^s\lambda=\lambda,$ i.e. $s\in Stab_{T}(\lambda)=S_{124}.$ There is unique $s=1\in S\cap S_{124}$ since $S$ is a representative set of  $K_1\backslash U/K_2.$ Therefore, $({\lambda_1}^U, {\eta}^U)=(\lambda_1|_{HF_3F_{124}},\eta_2|_{HF_3F_{124}})=1$ iff $\lambda_1|_{F_i}=\eta|_{F_i},i\in \{124,3\}.$

So ${\lambda_1}^U=\eta^U={\lambda_2}^U\in Irr(U,\lambda)$ iff $\lambda_1|_{F_i}=\lambda_2|_{F_i},i\in\{124,3\}.~\Box$

\medskip
It is remarked that since $K_1,K_2\unlhd U,$ the double coset $K_1\backslash U/K_2$ equals $U/K_1K_2=U/HX_3S_{124}.$ Hence, we can pick above $S=X_1X_2$ as a transversal of $U/K_1K_2.$

%************************************************************************************

%%%%%%%%%%%%%%%%%%%%%%%%%%%%%%%%%%%%%%%%%%%%%%%%%%%%%%%%%%%%%%%%%%%%%%%%%%%%%%%%%%%%%%%%%%%%%%%

%\bigskip
\subsubsection{Proof of Theorem \ref{thm:Main2Chars}}\label{proof-thm:Main2Chars}
Fix $a_8,a_9,a_{10}\in\F_q^\times$ and set $\lambda=\lambda_{b_5,b_6,b_7}^{a_8,a_9,a_{10}}$ for some $b_5,b_6,b_7\in\F_q$  in the whole proof. By Lemma \ref{lem:2extensions} and using the same notations, we mainly find the generic character values: in (a) $\chi_{8,9,10,q^3}^{a_8,a_9,a_{10}}={\eta_1}^U$ where $t_0=0,$ and in (b) $\chi_{8,9,10,\frac{q^3}{2}}^{b_{124},b_3,t_0,a_8,a_9,a_{10}}={\eta_1}^U$ where $b_{124},b_3\in \F_2,$ $t_0\in\F_q^\times.$

(a)  Suppose $t_0=0$ and $F_{124}=\{1\}.$ Call $\eta$ an extension of $\lambda$ to $HX_3.$ By Lemma \ref{lem:2extensions} (b), $I_U(\eta)=HX_3.$ Therefore, ${\eta}^U\in Irr(U)$ and ${\eta}^U(1)=q^3.$ By Lemma \ref{lem:2extensions} (a), $S_{567}X_8X_9X_{10}\subset Z({\eta}^U),$ hence $|{\eta}^U(x)|=q^3$ for all $x\in S_{567}X_8X_9X_{10}.$ We have $|S_{567}X_8X_9X_{10}|q^3q^3=q^{10}=|U|.$ By the scalar product $({\eta}^U,{\eta}^U)=1,$ it forces ${\eta}^U(x)=0$ if $x\notin S_{567}X_8X_9X_{10}.$ So we have the formula as stated.

%It is remarked that ${\lambda_1}^U$ is the unique character in $Irr(U,\lambda).$ Since $X_1X_2X_4$ is a transversal of $U/HX_3,$ for fixed $a_8,a_9,a_{10}\in\F_q^\times,$ $X_1X_2X_4$ acts transitively on the set of all linear characters $\lambda_{-,b_5,b_6,b_7}^{a_8,a_9,a_{10}}$ of $HX_3$ such that $b_5a_{10}+b_6a_9+b_7a_8=0.$

\medskip
(b) Suppose $t_0\neq 0,$  and $|F_3|=|F_{124}|=2.$ By Lemma \ref{lem:2extensions} (b), let $\eta_1,\eta_2$ be extensions of $\lambda$ to $K_1:=HX_3F_{124}$ and $K_2:=HS_{124}F_3$ respectively such that $\eta_1|_{F_i}=\eta_2|_{F_i}=\phi_{b_i},$ where $b_i\in\F_2,~ i\in\{124,3\}.$ By the proof of Lemma \ref{lem:2extensions} (c), ${\eta_1}^U={\eta_2}^U.$

We choose $V\subset T$ as a transversal of $K_1$ in $U,$ and  $1\in S\subset X_3$ such that $SX_1X_2$ is a transversal of $K_2$ in $U,$ so $|S|=q/2.$  Since $K_1\unlhd U,$ we have ${\eta_1}^U(\prod_{i=1}^{10}x_i)=\sum_{x\in V}{}^x\eta_1(\prod_{i=1}^{10}x_i)=0$ if $x_1x_2x_4\notin K_1.$ Since $T$ is abelian, $[x,y]=1$ for all $x\in V$ and $y\in F_{124}.$ Therefore, $F_{124}\subset Z({\eta_1}^U)$ and we have
\begin{center}
${\eta_1}^U(\prod_{i=1}^{10}x_i(t_i))=\delta_{a_8t_1,a_{10}t_4}\delta_{a_8t_2,a_9t_4}\phi(b_{124}\frac{t_1}{a_{10}}){\eta_1}^U(x_3(t_3)\prod_{i=5}^{10}x_i(t_i)).$
\end{center}

Since $K_2\unlhd U,$ we have ${\eta_2}^U(x_3\prod_{i=5}^{10}x_i)=\sum_{x\in SX_1X_2}{}^x\eta_2(x_3\prod_{i=5}^{10}x_i)=0$ if $x_3\notin F_3.$ Since $X_8X_9X_{10}\subset Z(U),$ we need to compute the two following cases: ${\eta_2}^U(\prod_{i=5}^7x_i)$ and ${\eta_2}^U(x_3\prod_{i=5}^7x_i)$ with $x_3\in F_3^\times.$

Since $[X_3,X_5X_6X_7]=\{1\},$ we have ${\eta_2}^U(x_5x_6x_7)=\sum_{x\in SX_1X_2}{}^x\eta_2(x_5x_6x_7)=\frac{q}{2}\sum_{x_1x_2\in X_1X_2}{}^x\eta_2(x_5x_6x_7).$ Since $(x_5x_6x_7)^{x_1x_2}=x_5x_6x_7[x_5,x_2][x_6,x_1][x_7,x_1][x_7,x_2],$ plug in $x_5(t_5)x_6(t_6)x_7(t_7)$ and $x_1(s_1)x_2(s_2),$ we have

\begin{center}
$\begin{array}{l}
{\eta_2}^U(\prod_{i=5}^7x_i(t_i))\\
=\frac{q}{2}\sum_{s_1,s_2}\eta_2(x_5(t_5)x_6(t_6)x_7(t_7)x_8(t_5s_2+t_6s_1)x_9(t_7s_1)x_{10}(t_7s_2))
\\
=\frac{q}{2}\eta_2(x_5(t_5)x_6(t_6)x_7(t_7))\sum_{s_1,s_2}\phi(a_8(t_5s_2+t_6s_1)+a_9t_7s_1+a_{10}t_7s_2)
\\
=\frac{q}{2}\eta_2(x_5(t_5)x_6(t_6)x_7(t_7))\sum_{s_1,s_2}\phi(s_1(a_8t_6+a_9t_7)+s_2(a_8t_5+a_{10}t_7)).
\end{array}$
\end{center}

Since $\sum_{t\in\F_q}\phi(t)=0,$ to obtain non-zero values, it forces $a_8t_6+a_9t_7=0$ and $a_8t_5+a_{10}t_7=0.$ Hence, $\frac{t_5}{a_{10}}=\frac{t_6}{a_9}=\frac{t_7}{a_8},$ and $\prod_{i=5}^7x_i(t_i)=x_{567}(\frac{t_7}{a_8})\in S_{567}.$ By Lemma \ref{lem:2extensions} (a), we have
\begin{center}
$\begin{array}{ll}
{\eta_2}^U(\prod_{i=5}^7x_i(t_i))&=\delta_{a_8t_5,a_{10}t_7}\delta_{a_8t_6,a_9t_7}\frac{q^3}{2}\eta_2(x_{567}(\frac{t_7}{a_8}))
\\
&=\delta_{a_8t_5,a_{10}t_7}\delta_{a_8t_6,a_9t_7}\frac{q^3}{2}\phi(At_0\frac{t_7}{a_8}).
\end{array}$
\end{center}

Therefore, ${\eta_2}^U(\prod_{i=1}^{10}x_i(t_i))=\frac{q^3}{2}\phi(b_{124}\frac{t_1}{a_{10}}+At_0\frac{t_7}{a_8}+\sum_{i=8}^{10}a_it_i)$  if $\prod_{i=1}^{10}x_i(t_i)\in F_{124}S_{567}X_8X_9X_{10}=Z,$ as stated in the theorem.

Now we compute ${\eta_2}^U(x_3\prod_{i=5}^7x_i)$ with $x_3\in F_3^\times=\{x_3(t_0^\phi)\}$ where $t_0^\phi=\frac{(t_0)_\phi}{A}.$ Since $I_U(\eta_2)=K_2\unlhd U$ and $SX_1X_2$ is a representative set of $U/K_2,$ $({}^x\eta_2)^U={\eta_2}^U\in Irr(U)$ for all $x\in SX_1X_2.$
For each $x_2(s)\in X_2,$ we have
\begin{center}
${}^{x_2(s)}\eta_2(x_5(t))=\eta_2(x_5(t)x_8(ts))=\phi(b_5t+a_8ts)=\phi(t(b_5+a_8s)).$
\end{center}
So instead of choosing $s=\frac{b_5}{a_8},$ we suppose that $\eta_2$ has $b_5=0,$ i.e. $\eta_2(x_5)=1$ for all $x_5\in X_5.$ It is easy to check that $t_0,$ $\eta_2|_{F_{124}}=\phi_{b_{124}}$ and $\eta_2|_{F_3}=\phi_{b_3}$ are invariant under this conjugate action.

We have
$[x_3(t_3)x_5(t_5)x_6(t_6)x_7(t_7),x_1(s_1)x_2(s_2)]$
\\
$=x_3(t_3)x_5(t_5+t_3s_1)x_6(t_6+t_3s_2)x_7(t_7)x_8(t_3s_1s_2+t_5s_2+t_6s_1)x_9(t_7s_1)x_{10}(t_7s_2).$
Therefore,
\\
$\begin{array}{l}
{\eta_2}^U(x_3(t_3)x_5(t_5)x_6(t_6)x_7(t_7))=\sum_{x\in SX_1X_2}{}^x\eta_2(x_3(t_3)x_5(t_5)x_6(t_6)x_7(t_7))
\\
=\frac{q}{2}\sum_{x\in X_1X_2}{}^x\eta_2(x_3(t_3)x_5(t_5)x_6(t_6)x_7(t_7))
\\
=\frac{q}{2}\sum_{s_1,s_2}\eta_2(x_3(t_3)x_5(t_5+t_3s_1)x_6(t_6+t_3s_2)x_7(t_7)x_8(t_3s_1s_2+t_5s_2+t_6s_1)x_9(t_7s_1)x_{10}(t_7s_2))
\\
=\frac{q}{2}\eta_2(x_3(t_3)\prod_{i=5}^7x_i(t_i))\sum_{s_1,s_2}\phi(b_6t_3s_2+a_8(t_3s_1s_2+t_5s_2+t_6s_1)+a_9t_7s_1+a_{10}t_7s_2)
\\
=\frac{q}{2}\eta_2(x_3(t_3)x_6(t_6)x_7(t_7))\sum_{s_1,s_2}\phi(s_1(a_8t_3s_2+a_8t_6+a_9t_7)+s_2(b_6t_3+a_{10}t_7+a_8t_5)).
\end{array}$

Set $C(t_5,t_6,t_7)=\sum_{s_1,s_2}\phi(s_1(a_8t_3s_2+a_8t_6+a_9t_7)+s_2(b_6t_3+a_{10}t_7+a_8t_5)).$ We have
\begin{center}
$\begin{array}{ll}
C(t_5,t_6,0)&=\sum_{s_1,s_2}\phi(s_1(t_3s_2+t_6)a_8+s_2(b_6t_3+a_8t_5))
\\
&=q\sum_{s_2=\frac{t_6}{t_3}}\phi(\frac{t_6}{t_3}(b_6t_3+a_8t_5))
\\
&=q\phi(b_6t_6+\frac{a_8t_5t_6}{t_3}).
\end{array}$
\end{center}

Therefore, we get
\begin{center}
$\begin{array}{ll}
{\eta_2}^U(x_3(t_3)x_5(t_5)x_6(t_6))&=\frac{q^2}{2}\eta_2(x_3(t_3)x_6(t_6))\phi(b_6t_6+\frac{a_8t_5t_6}{t_3})
\\
&=\frac{q^2}{2}\eta_2(x_3(t_3))\phi(\frac{a_8t_5t_6}{t_3}).
\end{array}$
\end{center}

Since $x_5(t_5)x_6(t_6)x_7(t_7)=x_5(t_5+\frac{a_{10}t_7}{a_8})x_6(t_6+\frac{a_9t_7}{a_8})x_{567}(\frac{t_7}{a_{10}}),$ where $x_{567}(\frac{t_7}{a_8})=x_5(a_{10}\frac{t_7}{a_8})x_6(a_9\frac{t_7}{a_8})x_7(a_8\frac{t_7}{a_8})\in S_{124}\subset Z({\eta_2}^U)$ and $\eta_2(x_{567}(\frac{t_7}{a_8}))=\phi_{At_0}(\frac{t_7}{a_8}),$ we have
\begin{center}
$\begin{array}{ll}
{\eta_2}^U(x_3(t_3)x_5(t_5)x_6(t_6)x_7(t_7))&=\phi_{At_0}(\frac{t_7}{a_8}){\eta_2}^U(x_3(t_3)x_5(t_5+\frac{a_{10}t_7}{a_8})x_6(t_6+\frac{a_9t_7}{a_8}))
\\
&=\phi(At_0\frac{t_7}{a_8})\frac{q^2}{2}\eta_2(x_3(t_3))\phi(\frac{a_8}{t_3}(t_5+\frac{a_{10}t_7}{a_8})(t_6+\frac{a_9t_7}{a_8}))
\\
&=\frac{q^2}{2}\phi(b_3t_3+At_0\frac{t_7}{a_8}+\frac{a_8^2a_9a_{10}}{(t_0)_\phi}(t_5+\frac{a_{10}t_7}{a_8})(t_6+\frac{a_9t_7}{a_8}))
\\
&=\frac{q^2}{2}\phi(b_3t_0^\phi+At_0\frac{t_7}{a_8}+\frac{A^2}{(t_0)_\phi}(\frac{t_5}{a_{10}}+\frac{t_7}{a_8})(\frac{t_6}{a_9}+\frac{t_7}{a_8})).~\Box
\end{array}$
\end{center}

%******************************************************************************

%%%%%%%%%%%%%%%%%%%%%%%%%%%%%%%%%%%%%%%%%%%%%%%%%%%%%%%%%%%%%%%%%%%%%%%%%%%%%%%%%%%%%%%%%%%%%%%

\subsection{Proof of Sylow $3$-subgroups of $E_6(3^f)$}

\subsubsection{Proof of Lemma \ref{lem:3extensions}}
\label{proof-lem:3extensions}

Set $\lambda=\lambda_{b_8,b_9,b_{10}}^{b_{12},b_{13},b_{14},b_{15},b_{16}}$ for the whole proof.

(a) Recall $H_3=\prod_{i=12}^{16}X_i\unlhd H$ is elementary abelian and $H_4H_3\unlhd U.$ First, we show that $R_3=\{x\in H_3:|\lambda^U(x)|=\lambda^U(1)\}$ and $\lambda^U(r_3(t))=q^8\phi_{B_3}(t)$ for all $r_3(t)\in R_3.$

Since $\lambda$ is linear and $\lambda(x)\in \C$ for all $x\in H,$ by the induction formula, we have $|\lambda^U(x)|=\lambda^U(1)$ iff %$\lambda^U(x)=\lambda^U(1)\lambda(x),$ i.e.
${}^y\lambda(x)=\lambda(x)$ for all $y\in TX_4$ which is a transversal of $H$ in $U.$
Since ${}^y\lambda(x)=\lambda(x)$ iff $\lambda([x,y])=1,$ we are going to find all $x\in H_3$ such that $\lambda([x,y])=1$ for all $y\in TX_4.$ It is clear that $[X_i,X_4]=\{1\}=[X_i,X_7X_{11}]$ for all $12\leq i\leq 16.$ Here, we write $\prod_{j=1}^6x_i(u_j)\in T$ with $u_4=0,$ it suffices to check for all $y=\prod_{j=1}^6x_j(u_j)\in T.$  For $t_i,u_j\in \F_q,$ we have
\begin{center}
$\begin{array}{l}
[\prod_{i=12}^{16}x_i(t_i),\prod_{j=1}^6x_j(u_j) ]=\\
\left[x_{12}(t_{12}),x_2(u_2)\right] [x_{15}(t_{15}),x_2(u_2)][x_{16}(t_{16}),x_2(u_2)][x_{13}(t_{13}),x_1(u_1)]\\
\left[x_{15}(t_{15}),x_1(u_1)\right][x_{14}(t_{14}),x_3(u_3)][x_{16}(t_{16}),x_3(u_3)][x_{12}(t_{12}),x_5(u_5)]\\
\left[x_{13}(t_{13}),x_5(u_5)\right][x_{14}(t_{14}),x_6(u_6)][x_{15}(t_{15}),x_6(u_6)]\\
=x_{17}(-t_{12}u_2)x_{19}(-t_{15}u_2)x_{20}(-t_{16}u_2)x_{17}(-t_{13}u_1)x_{18}(-t_{15}u_1)x_{19}(-t_{14}u_3)\\
x_{21}(-t_{16}u_3)x_{18}(t_{12}u_5)x_{19}(t_{13}u_5)x_{20}(t_{14}u_6)x_{21}(t_{15}u_6)
\end{array}$
\end{center}

Since $\lambda(x_i(t))=\phi(t), 17\leq i\leq 21,$ for all $u_j\in\F_q$ it forces
\begin{center}
$(-t_{12}-t_{15}-t_{16})u_2+(-t_{13}-t_{15})u_1+(-t_{14}-t_{16})u_3+(t_{12} +t_{13})u_5+( t_{14}+t_{15})u_6=0.$
\end{center}
So we have a system with variables $t_i$:
\begin{center}
$\left\{\begin{array}{rl}
-t_{12}-t_{15}-t_{16} &=0\\
 -t_{13}-t_{15}        &=0\\
      -t_{14}-t_{16} &=0\\
 t_{12} +t_{13}       &=0\\
  t_{14}+t_{15}        &=0
\end{array}\right.$
\end{center}

Since $\gcd(q,3)=3,$ we have $t_{12}=t_{16}=t_{15},t_{13}=t_{14}=-t_{15}$ for all $t_{15}=t\in\F_q.$ So $x\in H_3$ satisfies $|\lambda^U(x)|=\lambda^U(1)$ iff $x=x_{12}(t)x_{13}(-t)x_{14}(-t)x_{15}(t)x_{16}(t)=r_3(t)\in R_3$ for $t\in \F_q,$ i.e. $R_3=\{x\in H_3:|\lambda^U(x)|=\lambda^U(1)\}.$

By the above computation, to show that $\lambda^U|_{R_3}=\lambda^U(1)\phi_{B_3}$, it is enough to check that $\lambda(r_3(t)=\phi_{B_3}(t).$ For each $r_3(t)\in R_3,$ we have
\begin{center}
$\lambda(r_3(t))=\phi(t(b_{12}-b_{13}-b_{14}+b_{15}+b_{16}))=\phi_{B_3}(t).$
\end{center}

Now we show that $S_1=Stab_{T}(\lambda|_{H_4H_3}).$ Since $[H_4,T]=[H_3,X_7X_{11}]=\{1\},$ it suffices to find $y\in X_2X_1X_3X_5X_6$ such that $\lambda([x,y])=1$ for all $x\subset H_3.$ Using the above computation of $[\prod_{i=12}^{16}x_i(t_i),\prod_{j=1}^6x_j(u_j) ],$ for all $t_j\in\F_q$ it forces

\begin{center}
$(-u_2+u_5)t_{12}+(-u_1+u_5)t_{13}+(-u_3+u_6)t_{14}+(u_1+u_2-u_6)t_{15}+(-u_2-u_3)t_{16}=0.$
\end{center}
So we have a system with variables $u_j:$
\begin{center}
$\left\{\begin{array}{rl}
 -u_2+u_5 &=0\\
 -u_1+u_5 &=0\\
 -u_3+u_6 &=0\\
 -u_1-u_2+u_6 &=0\\
 -u_2-u_3 &=0.
\end{array}\right.$
\end{center}

Since $\gcd(q,3)=3,$ we have $u_1=u_5=u_2,u_3=u_6=-u_2$ for all $u_2=t\in\F_q.$ So $\prod_{j=1}^6x_j(u_j)=x_2(t)x_1(t)x_3(-t)x_5(t)x_6(-t)=s_1(t)\in S_1.$

\medskip
(b) Since $H=Z(U)H_3X_8X_9X_{10},$ to find $Stab_T(\lambda)\subset Stab_{T}(\lambda|_{H_4H_3})=S_1,$ by (a), it is enough to find $s_1\in S_1$ such that ${}^{s_1}\lambda(x_i)=\lambda(x_i)$ for $i=8,9,10.$ Again, for each $x_i(t_i)\in X_i,$ $i=8,9,10$ and $s_1=s_1(t,r,s)%=x_2(t)x_1(t)x_3(-t)x_5(t)x_6(-t)x_7(r)x_{11}(s)
\in S_1,$ we compute $[x_i,s_1]$
%\begin{center}
%$\begin{array}{ll}
\\
$ [x_8(t_8),s_1 ]=x_{20}(t_8s)x_{17}(-t_8r)x_{14}(t_8t)x_{20}(-t_8t^2)x_{13}(t_8t)x_{19}(t_8t^2).$
\\
$ [x_9(t_9),s_1 ]=x_{21}(t_9s)x_{15}(t_9t)x_{21}(-t_9t^2)x_{12}(-t_9t)x_{18}(-t_9t^2)x_{13}(-t_9t)x_{19}(-t_9t^2)x_{17}(t_9t^2).$
\\
$ [x_{10}(t_{10}),s_1 ]=x_{18}(-t_{10}r)x_{16}(-t_{10}t)x_{15}(t_{10}t)x_{21}(-t_{10}t^2)x_{14}(-t_{10}t)x_{20}(t_{10}t^2)x_{19}(-t_{10}t^2)$
%\end{array}$

%\end{center}

Since $\lambda(x_i(t))=\phi(b_it), 12\leq i\leq 16$ and $\lambda(x_i(t))=\phi(t),17\leq i\leq 21,$ from $\lambda([x_9(t_9),s_1])=\lambda([x_{10}(t_{10}),s_1])=1$ for all $t_9,t_{10}\in \F_q,$ we have
\begin{center}
$s=2t^2+b_{12}t+b_{13}t-b_{15}t$ and  $r=2t^2-b_{14}t+b_{15}t-b_{16}t.$
\end{center}
From $\lambda([x_8(t_8),s_1])=1$ for all $t_8\in\F_q,$ we have
$s-r+b_{14}t+b_{13}t=0.$ Therefore, $s_1\in Stab_T(\lambda)$ iff $r,s$ as above and
\begin{center}
$\begin{array}{ll}
s-r+b_{14}t+b_{13}t&=2t^2+b_{12}t+b_{13}t-b_{15}t-(2t^2-b_{14}t+b_{15}t-b_{16}t)+b_{14}t+b_{13}t
\\
&=t(b_{12}-b_{13}-b_{14}+b_{15}+b_{16})
\\
&=tB_3=0.
\end{array}$
\end{center}
Therefore, if $B_3\neq 0,$ then $Stab_T(\lambda)=\{1\}.$

\medskip
(c) By (a), $T/S_1$ acts faithfully on the set of all extensions of $\lambda|_{H_4}$ to $H_4H_3$ with the same $B_3.$ Since $|T/S_1|=q^4=|H_3/R_3|,$ this action is transitive. Therefore, with $B_3=0,$ there exists $x\in T$ such that ${}^x\lambda=\lambda_{b_8',b_9',b_{10}'}^{0,0,0,0,0}$ for some $b_8',b_9',b_{10}'\in \F_q.$

%Continue the computation in (b), when $B_3= 0,$ we have $Stab_{HT}(\lambda)=HS_2$ if $\lambda=\lambda_{b_8,b_9,b_{10}}^{0,0,0,0,0}.$ Therefore, we show that there exists $x\in T$ such that ${}^x\lambda=\lambda_{b_8',b_9',b_{10}'}^{0,0,0,0,0}$ for some $b_8',b_9',b_{10}'\in \F_q.$

% Let $\lambda=\lambda_{b_8,b_9,b_{10}}^{b_{12},b_{13},b_{14},b_{15},b_{16}}.$ For each $x\in T$ and $y\in H_3,$ we have ${}^x\lambda(y)=\lambda(y^x)=\lambda(y[y,x])=\lambda([y,x])\lambda(y).$ By the computation of  $[\prod_{i=12}^{16}x_i(t_i),\prod_{j=1}^6x_j(u_j) ]$ in (a), by setting ${}^x\lambda=\lambda_{b_8',b_9',b_{10}'}^{b'_{12},b'_{13},b'_{14},b'_{15},b'_{16}}$ we have
% \begin{center}
% $\left(\begin{array}{l} b'_{12}\\b'_{13}\\b'_{14}\\b'_{15}\\b'_{16} \end{array}\right)
% =\left(\begin{array}{l} b_{12}\\b_{13}\\b_{14}\\b_{15}\\b_{16} \end{array}\right)+\left(\begin{array}{rrrrr}
% -1&0&0&1&0\\
% 0&-1&0&1&0\\
% 0&0&-1&0&1\\
% -1&-1&0&0&1\\
% -1&0&-1&0&0
% \end{array}\right) \left(\begin{array}{l} u_2\\u_1\\u_3\\u_5\\u_6 \end{array}\right)$
% \end{center}
%
% Thus, with $ \left(\begin{array}{r} u_2\\u_1\\u_3\\u_5\\u_6 \end{array}\right)= \left(\begin{array}{r} t\\-b_{14}+b_{15}+b_{16}+t\\b_{16}-t\\-b_{13}-b_{14}+b_{15}+b_{16}+t\\-b_{14}+b_{16}-t \end{array}\right),$ we have $\left(\begin{array}{l} b'_{12}\\b'_{13}\\b'_{14}\\b'_{15}\\b'_{16} \end{array}\right)=\left(\begin{array}{l} 0\\0\\0\\0\\0 \end{array}\right)$ for all $t\in \F_q.$

Now set $\lambda=\lambda_{b_8,b_9,b_{10}}^{0,0,0,0,0},$ and $\overline{H_3}$ is the normal closure of $H_3$ in $HX_4S_1.$ To show that $\overline{H_3}\subset ker(\lambda^{HX_4S_1})\unlhd HX_4S_1,$ it suffices to show that $H_3\subset ker(\lambda^{HX_4S_1}).$ By (a) $Stab_{TX_4}(\lambda|_{H_4H_3})=S_1X_4$ which is a transversal of $H$ in $HX_4S_1,$ the claim holds by the induction formula and $H_3\subset ker(\lambda).$

By Lemma \ref{lem:Reduction} for $G=U$ with $N=M=HX_4S_1,$ $X=X_1X_3X_5X_6,$ $Y=\overline{H_3}$ and $Z=H_4,$ the induction map from $Irr(HX_4S_1/\overline{H_3},\lambda)$ to $Irr(U,\lambda)$ is bijective. Since $\overline{H_3}\subset \lambda^{HX_4S_1},$ $Irr(HX_4S_1/\overline{H_3},\lambda)=Irr(HX_4S_1,\lambda).~\Box$

%%%%%%%%%%%%%%%%%%%%%%%%%%%%%%%%%%%%%%%%%%%%%%%%%%%%%%%%%%%%%%%%%%%%%%%%%%%%%%%%%%%%%%%%%%%%%%%

\subsubsection{Proof of Lemma \ref{lem:3ext_B1=0}}
\label{proof-lem:3ext_B1=0}
Recall  $R_2=\{r_2(t) := x_8(-t)x_9(t)x_{10}(t) : t\in \F_q\}\leq H_2=X_8X_9X_{10}$ and $\lambda=\lambda_{b_8,b_9,b_{10}}^{0,0,0,0,0}.$ By Lemma \ref{lem:3extensions} (c), it suffices to work with the quotient group $HX_4S_1/\overline{H_3}.$

(a) The fact $S_2 = Stab_{S_1}(\lambda)$ comes directly from Lemma \ref{lem:3extensions} (b) with $B_3=0.$ % that $Stab_T(\lambda)=S_2.$
 Since $X_4S_1$ is a transversal of $H$ in $HX_4S_1$ and $[H_2,X_4]=\{1\},$ to show $R_2=\{x\in H_2:|\lambda^{HX_4S_1}(x)|=\lambda^{HX_4S_1}(1)\}$ we are going to find all $x\in H_2$ such that $\lambda([x,y])=1$ for all $y\in S_1.$ Since $H\unlhd HX_4S_1$ is abelian, using the computation in Lemma \ref{lem:3extensions} (b), for $s_1(t,r,s)\in S_1$ and $x_8(t_8)x_9(t_9)x_{10}(t_10)\in H_2$ we have
\begin{center}
$\begin{array}{l}
[x_8(t_8)x_9(t_9)x_{10}(t_{10}),s_1(t,r,s)]\\
=[x_8(t_8),s_1(t,r,s)][x_9(t_9),s_1(t,r,s)][x_{10}(t_10),s_1(t,r,s)]\\
=x_{20}(t_8s)x_{17}(-t_8r)x_{20}(-t_8t^2)x_{19}(t_8t^2)x_{21}(t_9s)x_{21}(-t_9t^2)x_{18}(-t_9t^2)x_{19}(-t_9t^2)\\
x_{17}(t_9t^2)x_{18}(-t_{10}r)x_{21}(-t_{10}t^2)x_{20}(t_{10}t^2)x_{19}(-t_{10}t^2)
\end{array}$
\end{center}
Therefore, with $\lambda|_{X_i}=\phi$ for all $17\leq i\leq 21,$ for all $t,r,s\in \F_q$ we need
\begin{center}
$(t_8+t_9)s-(t_8+t_{10})r+(t_9-t_{10})t^2=0$
\end{center}
So $t_9=t_{10}=u$ and $t_8=-u$ for all $u\in \F_q,$ i.e. $x=r_2(u)\in R_2.$
%% $x\in H_2$ such that $\lambda([x,s_1(t,r,s)])=1$ for all $r,s,t\in\F_q$ iff $x=r_2(t)\in R_2.$

To show that $\lambda^{HX_4S_1}|_{R_2}=\lambda^{HX_4S_1}(1)\phi_{B_2},$ it is enough to check that $\lambda(r_2(t))=\phi_{B_2}(t).$ For each $r_2(t)\in R_2$ we have
\begin{center}
$\lambda(r_2(t))=\phi(t(-b_8+b_9+b_{10}))=\phi_{B_2}(t).$
\end{center}

\medskip
(b) Suppose that $B_2\not\in \{c^2:c\in \F_q^\times\}.$ Let $\eta$ be an extension of $\lambda$ to $HX_4.$
%Since $HX_4S_1=HX_4\rtimes S_1,$ we have $Stab_{S_1}(\eta)\subset Stab_{S_1}(\lambda)=\{1\}.$ Therefore, $I_U(\eta)=HStab_{TX_4}(\eta)=HX_4.$
By (a) that $S_2=Stab_{S_1}(\lambda),$ hence $Stab_{S_1}(\eta|_{H})=S_2.$
%Stab_{S_2}(\eta).$
Since $S_1$ is a tranversal of $HX_4$ in $HX_4S_1,$ to find $Stab_{S_1}(\eta),$ it is enough to find all $s_2(t)\in S_2$ such that $\eta([x_4,s_2(t)])=1$ for all $x_4\in X_4.$ For each $s_2(t)
%=x_2(t)x_1(t)x_3(-t)x_5(t)x_6(-t)x_7(2t^2)x_{11}(2t^2)
\in S_2,$ we have
\begin{center}
$\begin{array}{ll}[x_4(t_4),s_2(t)]=&%x_{16}(t_42t^2)x_{12}(-t_42t^2)
x_{10}(t_4t)%x_{18}(-t_4t^2)x_{16}(-t_4t^2)
x_9(t_4t)x_{21}(2t_4t^3)%x_{15}(t_4t^2)$\\
x_{21}(-t_4t^3)x_8(-t_4t)\\&x_{20}(-2t_4t^3)x_{17}(2t_4t^3)%x_{14}(-t_4t^2)
x_{20}(t_4t^3)%x_{13}(-t_4t^2)
x_{19}(-t_4t^3).
\end{array}$
\end{center}
Since $\eta(x_i(t))=\phi(b_it),8\leq i\leq 10$ and $\eta(x_i(t))=\phi(t), 17\leq i\leq 21,$ for all $t_4\in\F_q,$ $\eta([x_4(t_4),s_3(t)])=1$ forces
$$t_4(t^3-B_2t)=t_4(t^3-B_2t)\in ker\phi.$$

Since $t_4(t^3-B_2t)\in ker\phi$ for all $t_4\in\F_q,$ we have $0=t^3-B_2t=t(t^2-B_2).$ Since $B_2\notin \{c^2:c\in \F_q^\times\},$ the equation $t(t^2-B_2)=0$ only has trivial solution $t=0$ in $\F_q.$ Therefore, $s_2(t)=1,$ i.e. $Stab_{S_1}(\eta)=\{1\}.$ Hence, $I_{HX_4S_1}(\eta)=HX_4.$

\medskip
(c) Suppose $B_2=c^2\in\F_q^\times$ and let $\eta$ be an extension of $\lambda$ to $HX_4.$ Using the computation in (b), we continue with the analysis for the solutions of $t$ to obtain $t_4t(t^2-B_2)\in ker\phi$ for all $t_4\in \F_q.$ So it forces $t(t^2-B_2)=0.$ This equation has 3 solutions $\{0,\pm c\}.$ Hence, $Stab_{S_1}(\eta)=\{1,s_2(\pm c)\}=F_2.$ So $I_{HX_4S_1}(\eta)=HX_4F_2.$
By the above argument, $[HX_4F_2,HX_4F_2]\subset ker(\eta),$ hence $\eta$ extends to $I_{HX_4S_1}(\eta).$ % by checking $[I_U(\eta),I_U(\eta)]\subset ker(\eta).$

To show that $\lambda$ extends to $HF_4S_2,$ we check $[HF_4S_2,HF_4S_2]\subset ker(\lambda).$  With the same argument, it is enough to check that $[s_2(t),x_4(t_4)]\in ker\lambda.$ By the computation in (b), we need $t_4(t^3-B_2t)=t_4(t^3-c^2t)\in ker\phi$ for all $t\in \F_q.$ By Proposition \ref{prop:mainFp}, since $t_4\in\{0,\pm c_\phi\},$ the claim holds.

Let $\lambda_1,\lambda_2$ be two extensions of $\lambda$ to $HX_4F_2,$ and $\gamma$ an extension of $\lambda$ to $HF_4S_2.$ Since the degree of all irreducible constituents of $\lambda^{HX_4S_1}$ is $\frac{q^3}{3},$ we have ${\lambda_1}^{HX_4S_1},$ ${\lambda_2}^{HX_4S_1},$ $\gamma^{HX_4S_1}\in Irr(HX_4S_1,\lambda).$

Choose $1\in S\subset S_1$ as a representative set of the double coset $HF_4S_2\backslash HX_4S_1/HX_4F_2.$ Since $HF_4S_2\cap HX_4F_2=HF_4F_2$ and $HX_4F_2\unlhd HX_4S_1,$ by Mackey formula,
\begin{center}
$\begin{array}{ll}
({\lambda_1}^{HX_4S_1},\gamma^{HX_4S_1})&=\sum_{s\in S}({}^s\lambda_1|_{{}^s(HX_4F_2)\cap HF_4S_2},\gamma|_{{}^s(HX_4F_2)\cap HF_4S_2})
\\
&=\sum_{s\in S}({}^s\lambda_1|_{HF_4F_2},\gamma|_{HF_4F_2})
\end{array}$
\end{center}

For each $s\in S,$ if ${}^s\lambda_1|_{HF_4F_2}=\gamma|_{HF_4F_2},$ then ${}^s\lambda_1|_H=\gamma|_H.$ Since both are extensions of $\lambda,$ we have ${}^s\lambda=\lambda,$ i.e. $s\in Stab_{S_1}(\lambda)=S_2.$ There is unique $1\in S\cap S_2$ since $S$ is a representative set of $HF_4S_2\backslash HX_4S_1/HX_4F_2.$ So $({\lambda_1}^{HX_4S_1},\gamma^{HX_4S_1})=(\lambda_1|_{HF_4F_2},\gamma|_{HF_4F_2})=1$ iff $\lambda_1|_{F_i}=\gamma|_{F_i},i\in \{2,4\}.$

Therefore, ${\lambda_1}^{HX_4S_1}=\gamma^{HX_4S_1}={\lambda_2}^{HX_4S_1}$ iff $\lambda_1|_{F_i}=\lambda_2|_{F_i},i\in \{2,4\}.~\Box$

%%%%%%%%%%%%%%%%%%%%%%%%%%%%%%%%%%%%%%%%%%%%%%%%%%%%%%%%%%%%%%%%%%%%%%%%%%%%%%
%%%%%%%%%%%%%%%%%%%%%%%%%%%%%%%%%%%%%%%%%%%%%%%%%%%%%%%%%%%%%%%%%%%%%%%%%%%%%%
%%%%%%%%%%%%%%%%%%%%%%%%%%%%%%%%%%%%%%%%%%%%%%%%%%%%%%%%%%%%%%%%%%%%%%%%%%%%%%

\subsection{Proofs of Sylow 5-subgroups of $E_8(5^f)$}

\subsubsection{Proof of Lemma \ref{lem:5ext1}}
\label{proof-lem:5ext1}
(a) First we find all $x\in H_5$ such that $|\lambda^U(x)|=\lambda^U(1).$ Since $TX_5$ is a transversal of $H$ in $U,$ $[H_5,X_5]=\{1\}=[H_5,T_k]$ for all $k\geq 2,$ and ${}^y\lambda(x)=\lambda(x)$ iff $\lambda([x,y])=1,$ it suffices to find all $x\in H_5$ such that $\lambda([y,x])=1$ where $y\in T_1.$ For each $y=\prod_{1}^8x_i(u_i)\in T_1$ with $u_5=0,$ and $x=\prod_{j=30}^{36}x_j(v_j)\in H_5,$ to write shortly for the decomposition we write $x_i=x_i(-)$ and plug in the parameters in $(-)$ latter, we have
%\begin{center}
%$\begin{array}{l}
\\
$[\prod_{j=30}^{36}x_j(v_j),\prod_{1}^8x_i(u_i)]=\left[x_{30},x_4\right][x_{30},x_6][x_{31},x_2][x_{31},x_7][x_{32},x_1][x_{32},x_6][x_{33},x_1]\\
\left[x_{33},x_4\right][x_{33},x_7][x_{34},x_3][x_{34},x_8][x_{35},x_2][x_{35},x_1][x_{35},x_8][x_{36},x_2][x_{36},x_3]\\
=x_{37}(-v_{30}u_4)x_{38}(v_{30}u_6)x_{38}(-v_{31}u_2)x_{39}(v_{31}u_7)x_{37}(-v_{32}u_1)x_{40}(v_{32}u_6)x_{38}(-v_{33}u_1)\\
x_{40}(-v_{33}u_4)x_{41}(v_{33}u_7)x_{41}(-v_{34}u_3)x_{42}(v_{34}u_8)x_{41}(-v_{35}u_2)x_{39}(-v_{35}u_1)x_{43}(v_{35}u_8)\\
x_{42}(-v_{36}u_2)x_{43}(-v_{36}u_3).$
%\end{array}$
%\end{center}

Since $\lambda|_{X_i}=\phi$ for all $i\in[37..43],$ for all $s_j$ we need
\\
$(-v_{31}-v_{35}-v_{36})u_2+(-v_{32}-v_{33}-v_{35})u_1+(-v_{34}-v_{36})u_3+(-v_{30}-v_{33})u_4+(v_{30}+v_{32})u_6+(v_{31}+v_{33})u_7+(v_{34}+v_{35})u_8=0.$
\\
Therefore, we obtain a system with variables $v_i$ as follows.
\begin{center}
$\left\{\begin{array}{rl}
-v_{31}-v_{35}-v_{36}&=0\\
-v_{32}-v_{33}-v_{35}&=0\\
-v_{34}-v_{36}&=0\\
-v_{30}-v_{33}&=0\\
v_{30}+v_{32}&=0\\
v_{31}+v_{33}&=0\\
v_{34}+v_{35}&=0.
\end{array}\right.$
\end{center}

Since $\gcd(q,5)=5,$ $(v_{30},v_{31},v_{32},v_{33},v_{34},v_{35},v_{36})=(v,v,-v,-v,-2v,2v,2v)$ for all $v\in \F_q.$ Hence, $x=r_5(v)\in R_5,$ i.e. $R_5=\{x\in H_5:|\lambda^U(x)|=\lambda^U(1)\}.$

To show that $\lambda^U|_{R_5}=\lambda^U(1)\phi_{B_5},$ it suffices to check that $\lambda(r_5(v))=\phi_{B_5}(v).$ For each $r_5(v)\in R_5,$ we have
\[\lambda(r_5(v))=\phi(v(b_{30}+b_{31}-b_{32}-b_{33}-2b_{34}+2b_{35}+2b_{36}))=\phi_{B_5}(v).\]

To show that $S_1=Stab_T(\lambda|_{H_6H_5}),$ we find all $y\in T$ such that $\lambda([x,y])=1$ for all $x\in H_6H_5.$ Since $H_6=Z(U)$ and $[H_5,T_k]=\{1\}$ for all $k\geq 2,$ it is enough to find $y\in T_1$ such that $\lambda([x,y])=1$ for all $x\in H_5.$ Using the above computation of $[\prod_{j=30}^{36}x_j(v_j),\prod_{1}^8x_i(u_i)]$, we find $u_i$ such that for all $v_j:$
\\
$(-u_4+u_6)v_{30}+(-u_2+u_7)v_{31}+(-u_1+u_6)v_{32}+(-u_1-u_4+u_7)v_{33}+(-u_3+u_8)v_{34}+(-u_2-u_1+u_8)v_{35}+(-u_2-u_3)v_{36}=0.$
\\
Therefore, we obtain a system with variables $u_i$ as follows.
\begin{center}
$\left\{\begin{array}{rl}
-u_4+u_6&=0\\
-u_2+u_7&=0\\
-u_1+u_6&=0\\
-u_1-u_4+u_7&=0\\
-u_3+u_8&=0\\
-u_2-u_1+u_8&=0\\
-u_2-u_3&=0.
\end{array}\right.$
\end{center}

Since $\gcd(q,5)=5,$ we have $(u_2,u_1,u_3,u_4,u_6,u_7,u_8)=(2u,u,-2u,u,u,2u,-2u)$ for all $u\in \F_q.$ So $y=l_1(u)\in L_1,$ i.e. $S_1=Stab_T(\lambda|_{H_6H_5}).$

\medskip
(b) Suppose $B_5\neq 0.$ To show that $Stab_T(\lambda)=\{1\},$ we are going to show that $Stab_{S_1}(\lambda|_{H_6H_5H_4})=T_3T_4,$ $Stab_{T_3T_4}(\lambda|_{H_6H_5H_4H_3})=T_4$ and $Stab_{T_4}(\lambda)=\{1\}.$

First, we show that $Stab_{S_1}(\lambda|_{H_6H_5H_4})=T_3T_4.$ By the root heights, it is clear that $[H_6H_5H_4,T_3T_4]=\{1\},$ hence, $T_3T_4\subset Stab_{S_1}(\lambda|_{H_6H_5H_4}).$ It suffices to show that $Stab_{L_1T_2}(\lambda|_{H_6H_5H_4})=\{1\},$ i.e. there is no nontrivial  $y\in L_1T_2$ such that $\lambda([h,y])=1$ for all $h\in H_5H_4.$  For each $y=\prod_{i=1}^{15}x_i(u_i)\in L_1T_2$ (with $u_5=u_{12}=u_{13}=0$ and $\prod_{i=1}^8x_i(u_i)=l_1(u)),$ and $h=\prod_{j=24}^{36}x_j(v_j)\in H_5H_4,$ we have
\\
$\left[\prod_{j=24}^{36}x_j(v_j),\prod_{i=1}^{15}x_i(u_i)\right]=[x_{24},x_{10}][x_{24},x_{14}][x_{25},x_{14}][x_{26},x_{15}][x_{26},x_{9}][x_{26},x_{11}]\\
\left[x_{27},x_{15}\right][x_{27},x_{10}][x_{29},x_{11}][x_{24},x_{2}][[x_{24},x_{2}],x_6][[x_{24},x_{2}],x_4][x_{24},x_{6}][[x_{24},x_{6}],x_7]\\
\left[x_{25},x_{1}\right][[x_{25},x_{1}],x_4][[x_{25},x_{1}],x_6][x_{25},x_{4}][[x_{25},x_{4}],x_6][x_{25},x_{6}][[x_{25},x_{6}],x_7][x_{26},x_{3}]\\
\left[[x_{26},x_{3}],x_4\right][[x_{26},x_{3}],x_7][x_{26},x_{7}][[x_{26},x_{7}],x_8][x_{27},x_{2}][[x_{27},x_{2}],x_1][[x_{27},x_{2}],x_4]\\
\left[[x_{27},x_{2}],x_7\right][x_{27},x_{1}][[x_{27},x_{1}],x_7][x_{27},x_{7}][[x_{27},x_{7}],x_8][x_{28},x_{2}][[x_{28},x_{2}],x_3]\\
\left[[x_{28},x_{2}],x_8\right][x_{28},x_{3}][[x_{28},x_{3}],x_8][x_{28},x_{8}][x_{29},x_{4}]\\
=x_{37}(v_{24}u_{10})x_{39}(v_{24}u_{14})x_{41}(v_{25}u_{14})x_{42}(v_{26}u_{15})x_{38}(-v_{26}u_{9})x_{40}(v_{26}u_{11})x_{43}(v_{27}u_{15})\\
x_{40}(v_{27}u_{10})x_{39}(-v_{28}u_{9})x_{42}(-v_{29}u_{10})x_{43}(-v_{29}u_{11})x_{30}(-2v_{24}u)x_{38}(-2v_{24}u^2)x_{37}(2v_{24}u^2)\\
x_{31}(v_{24}u)x_{39}(v_{24}2u^2)x_{30}(-v_{25}u)x_{37}(v_{25}u^2)x_{38}(-v_{25}u^2)x_{32}(-v_{25}u)x_{40}(-v_{25}u^2)x_{33}(v_{25}u)\\
x_{41}(2v_{25}u^2)x_{33}(2v_{26}u)x_{40}(-2v_{26}u^2)x_{41}(4v_{26}u^2)x_{34}(2v_{26}u)x_{42}(-4v_{26}u^2)x_{33}(-2v_{27}u)\\
x_{38}(2v_{27}u^2)x_{40}(2v_{27}u^2)x_{41}(-4v_{27}u^2)x_{31}(-v_{27}u)x_{39}(-2v_{27}u^2)x_{35}(v_{27}2u)x_{43}(-4v_{27}u^2)\\
x_{34}(-2v_{28}u)x_{41}(-4v_{28}u^2)x_{42}(4v_{28}u^2)x_{35}(2v_{28}u)x_{43}(-4v_{28}u^2)x_{36}(-v_{28}2u)x_{36}(-v_{29}u).$

Since $\lambda|_{X_i}=\phi$ for all $i\in[37..43]$ and $\lambda|_{X_i}=\phi_{b_i}$ for the others, after evaluating the above with $\lambda$ to get 1,  for all $v_j,$ we need
\\
$v_{24}(u_{10} + u_{14} - 2b_{30}u  +b_{31}u + 2u^2)
+ v_{25}(u_{14} -b_{30}u  - b_{32}u  + b_{33}u + u^2)
+ v_{26}(u_{15} - u_{9} + u_{11} + 2b_{33}u - 2u^2  + 2b_{34}u)
+ v_{27}(u_{15} + u_{10} - 2b_{33}u  - b_{31}u  + 2b_{35}u - u^2)
+ v_{28}(-u_{9} - 2b_{34}u  + 2b_{35}u + u^2 - 2b_{36}u)
+ v_{29}(-u_{10} - u_{11} - b_{36}u)=0.$

Hence, we have a system with variables $u_i$ and $u:$

\begin{center}
$\left\{\begin{array}{rl}
u_{10} + u_{14} - 2b_{30}u  +b_{31}u + 2u^2&=0\\
u_{14} -b_{30}u  - b_{32}u  + b_{33}u + u^2&=0\\
u_{15} - u_{9} + u_{11} + 2b_{33}u - 2u^2  + 2b_{34}u&=0\\
u_{15} + u_{10} - 2b_{33}u  - b_{31}u  + 2b_{35}u - u^2&=0\\
-u_{9} - 2b_{34}u  + 2b_{35}u + u^2 - 2b_{36}u&=0\\
-u_{10} - u_{11} - b_{36}u&=0.
\end{array}\right.$
\end{center}
It is equivalent to:
\begin{center}
$(*)\left\{
\begin{array}{l}
u_9=u^2+(3b_{34}+2b_{35}+3b_{36})u,\\
u_{10}=-u^2+(b_{30}-b_{31}-b_{32}+b_{33})u, \\
u_{11}= u^2+(-b_{30}+b_{31}+b_{32}-b_{33}-b_{36})u,\\
u_{14}=-u^2+(b_{30}+b_{32}-b_{33})u,\\
u_{15}=2u^2+(-b_{30}+2b_{31}+b_{32}+b_{33}+3b_{35})u,\\
(b_{30}+b_{31}-b_{32} -b_{33}- 2b_{34}+2b_{35}+ 2b_{36})u=0.
\end{array}\right.$
\end{center}

The last equation is actually $B_5u=0.$ Since $B_5\neq 0,$ we have $u=0$ and $u_9=u_{10}=u_{11}=u_{14}=u_{15}=0,$ i.e. $Stab_{L_1T_2}(\lambda|_{H_6H_5H_4})=Stab_{T_1T_2}(\lambda|_{H_6H_5H_4})=\{1\}.$

Thus $T_1T_2$ acts faithfully on the set of all extensions of $\lambda|_{H_6}$ to $H_6H_5H_4$ with the same $B_5\neq 0,$ which is invariant under the action of $T,$ i.e. $B_5(\lambda)=B_5({}^x\lambda)$ for all $x\in T.$ Since $|H_5H_4/R_5|=q^{12}=|T_1T_2|,$ this action is transitive. Therefore, we choose $\lambda|_{X_i}=\phi$ for all $i\in [37..43],$ $\lambda|_{X_{36}}=\phi_{B_5/2},$ and $\lambda|_{X_i}=1_{X_i}$ for the others $X_i\subset H_5H_4.$ By the root heights, we have $R_5\prod_{i=37}^{43}X_i\subset Z(HX_5T_4T_3),$ $HX_5T_4T_3\unlhd U$ and $H_4\prod_{i=30}^{35}X_i\unlhd HX_5T_4T_3.$ By Lemma \ref{lem:Reduction} for $G=U$ with $N=M=HX_5T_3T_4,$ $X=T_1T_2,$ $Z=H_6R_5$ and $Y=H_4\prod_{i=30}^{35}X_i,$ the induction map from $Irr(HX_5T_4T_3/Y,\lambda)$ to $Irr(U,\lambda)$ is bijective. Since $X_5T_4T_3=Stab_{X_5T}(\lambda|_{H_6H_5H_4})$ is a transversal of $H$ in $HX_5T_4T_3,$ we have $\lambda^{HX_5T_4T_3}|_{Y}=[HX_5T_4T_3:H]\lambda|_Y=|X_5T_4T_3|1_Y.$ Hence, $Irr(HX_5T_4T_3/Y,\lambda)=Irr(HX_5T_4T_3,\lambda).$

Now we find $Stab_{T_4T_3}(\lambda|_{H_6H_5H_4H_3}).$ Since $[H_6H_5H_4H_3,T_3T_4]=[H_3,T_3],$  we find $y\in T_3$ such that $\lambda([x,y])=1$ for all $x\in H_3.$ For each $y=\prod_{j=16,17,22}x_{j}(u_{j})\in T_3$ and $x=\prod_{i=18}^{21}x_i(v_i)\in H_3$ we have
\\
$[x,y]=[x_{18},x_{16}][x_{18} ,x_{22} ][x_{19} ,x_{22}][x_{20} ,x_{17}][x_{20},x_{16}][x_{21} ,x_{17} ]\\
=x_{37}(v_{18}u_{16})x_{42}(v_{18}u_{22})x_{43}(v_{19}u_{22})x_{40}(-v_{20}u_{17})x_{39}(-v_{21}u_{16})x_{41}(-v_{21}u_{17}).$

Since $\lambda|_{X_i}=\phi$ for all $i\in [37..43],$ for all $v_i$ we need
\[v_{18}(u_{16}+u_{22})+v_{19}u_{22}-v_{20}u_{17}+v_{21}(-u_{17}-u_{16})=0.\]
The only solution is $(u_{16},u_{17},u_{22})=(0,0,0),$ i.e. $Stab_{T_4T_3}(\lambda|_{H_6H_5H_4H_3})=T_4.$

Next, we find  $Stab_{T_4}(\lambda).$ Since $[H,T_4]=[H_2,T_4],$ we find $y\in T_4$ such that $\lambda([x,y])=1$ for all $x\in H_2.$ For each $y=x_{23}(u_{23})\in T_4$ and $x=x_{12}(v_{12})x_{13}(v_{13})\in H_2$ we have
\begin{center}
$\left[x_{12}(v_{12})x_{13}(v_{13}),x_{23}(u_{23})\right] =[x_{12},x_{23}][x_{13},x_{23}]=x_{37}(-v_{12}u_{23})x_{38}(-v_{13}u_{23}).$
\end{center}
Evaluate with $\lambda,$ for all $v_i$ we need
$(-v_{12}-v_{13})u_{23}=0.$
Therefore, the only solution is $u_{23}=0,$ i.e. $Stab_{T_4}(\lambda)=\{1\}.$ So we finish the proof of $Stab_T(\lambda)=\{1\}.$

Let $\eta,\eta'$ be two extensions of $\lambda|_{H_6H_5H_4}$ to $HX_5.$ By the bijection of the induction map from $Irr(HX_5T_4T_3,\lambda)$ to $Irr(U,\lambda),$ it suffices to show that $\eta^{HX_5T_4T_3}=\eta'^{HX_5T_4T_3}$ iff $\eta|_{R_j}=\eta'|_{R_j}$ for $j=2,3$ and $\eta|_{X_5}=\eta'|_{X_5}.$ By the Mackey formula for the double coset $HX_5\backslash HX_5T_4T_3/HX_5=HX_5T_4T_3/HX_5$ represented by $T_4T_3$ we have
\[(\eta^{HX_5T_4T_3},\eta'^{HX_5T_4T_3})=\sum_{y\in T_4T_3}({}^y\eta,\eta')\]
Since $[X_5,T_3T_4]\subset H_4\prod_{i=30}^{35}X_i\subset ker(\lambda),$ we have ${}^y\eta|_{X_5}=\eta|_{X_5}.$ Therefore, the restrictions to $X_5$ of both $\eta,$ $\eta'$ are clear for the proof. To show for the restrictions to $R_k$ with $k=2,3,$ we are going to prove that $R_2R_3=\{x\in H_2H_3:|\lambda^{HX_5T_4T_3}(x)|=\lambda^{HX_5T_4T_3}(1)\}$ and $T_4T_3=Stab_{T_4T_3}(\lambda|_{R_2R_3}).$ Then by $Stab_{T_4T_3}(\lambda)=\{1\}$ and $|T_4T_3|=q^4=|H_3H_2/R_2R_3|,$ the claim holds.

By the above computations of $[H_3,T_3]$ and $[H_2,T_4]$ we find all $x\in H_2H_3$ such that $\lambda([x,y])=1$ for all $y\in T_4T_3.$ For $y=x_{16}(u_{16})x_{17}(u_{17})x_{22}(u_{22})x_{23}(u_{23})\in T_3T_4$ and $x=x_{12}(v_{12})x_{13}(v_{13})\prod_{i=18}^{21}x_i(v_i)\in H_2H_3,$ we solve for $v_i$ in the following.
\[ u_{16}(v_{18}-v_{21})  +u_{17}(-v_{20}-v_{21}) + u_{22}(v_{18}+v_{19})+ u_{23}(-v_{12}-v_{13})=0\]

We have a system with variables $v_i:$

\begin{center}
$\left\{\begin{array}{rl}
v_{18}-v_{21}&=0\\
-v_{20}-v_{21}&=0\\
v_{18}+v_{19}&=0\\
-v_{12}-v_{13}&=0.
\end{array}\right.$
\end{center}

We obtain solutions $(v_{18},v_{19},v_{20},v_{21})=(v,-v,-v,v)$ and $(v_{12},v_{13})=(s,-s)$ for all $v,s\in\F_q.$ Therefore, $x\in R_2R_3.$ Hence, ${}^y\lambda|_{R_2R_3}=\lambda|_{R_2R_3}$ for all $y\in T_4T_3.$
%It is easy to check that $\lambda(r_k(v))=\phi_{B_k}(v)$ for all $r_k(v)\in R_k$ with $k=2,3.$

\medskip
(c) Suppose that $B_5=0.$ By (a), $T/S_1$ acts faithfully on the set of all extensions of $\lambda|_{H_6}$ to $H_6H_5$ with the same $B_5.$ Since $|H_5|/|R_5|=q^6=|T/S_1|,$ this action is transitive. Hence, there exists $x\in T$ such that ${}^x\lambda|_{X_i}=1_{X_i}$ for all $X_i\subset H_5.$ Let $\lambda$ be this linear. So $S_1X_5=Stab_{TX_5}(\lambda|_{H_6H_5}),$ a transversal of $H$ in $HX_5S_1,$ and $\lambda^{HX_5S_1}|_{H_5}=\lambda^{HX_5S_1}(1)\lambda|_{H_5},$ i.e. $H_5\subset ker(\lambda^{HX_5S_1})$ and so is its normal closure $\overline{H_5}$ in $HX_5S_1.$

By Lemma \ref{lem:Reduction} with  $G=U,$ $N=M=HX_5S_1,$ $X=\prod_{i=1}^4X_iX_6X_7,$ $Z=H_6$ and $Y=\overline{H_5},$ the induction map from $Irr(HX_5S_1/\overline{H_5},\lambda)$ to $Irr(U,\lambda)$ is bijective. Since $\overline{H_5}\subset  ker(\lambda^{HX_5S_1}),$ we have $Irr(HX_5S_1/\overline{H_5},\lambda)=Irr(HX_5S_1,\lambda).~\Box$

\subsubsection{Proof of Lemma \ref{lem:5ext2}}
\label{proof-lem:5ext2}
Recall that $\lambda$ is a linear character of $H$ such that $\lambda|_{X_i}=\phi$ for all $X_i\subset H_6,$ $\lambda|_{X_i}=1_{X_i}$ for all $X_i\subset H_5,$ and $\lambda|_{X_i}=\phi_{b_i}$ for the others $X_i\subset H_4H_3H_2$ where $b_i\in \F_q.$
By Lemma \ref{lem:5ext1} (c), we work with the quotient group $HX_5S_1/\overline{H_5}.$ Abusing the notation of root groups, we call them root groups in the quotient group. %Since $B_5=0,$  $Stab_T(\lambda|_{H_6H_5H_4})=S_2=L_2T_3T_4$ by Lemma \ref{lem:5ext1} (b).

\medskip
(a) By computation (*) in Lemma \ref{lem:5ext1} (b) with $B_5=0$, $S_2=Stab_{S_1}(\lambda|_{H_6H_5H_4}).$ Now we show that $R_4=\{x\in H_4:|\lambda^{HX_5S_1}(x)|=\lambda^{HX_5S_1}(1)\}.$ For each $l_1y_2y_3y_4\in L_1T_2T_3T_4=S_1$ and $h_4\in H_4,$ we have  $[h_4,l_1y_2y_3y_4]=[h_4,l_1y_2].$ Hence, we are going to find all $h_4\in H_4$ such that $\lambda([h_4,l_1y_2])=1$ for all $l_1y_2\in L_1T_2.$ Using the computation of $[\prod_{j=24}^{36}x_j(v_j),\prod_{i=1}^{15}x_i(u_i)]$ in Lemma \ref{lem:5ext1} (b) with $b_j=0$  for $j\in[30..36]$, we solve for $v_i$ in the following equation:
\\
$u_{9}(-v_{26}-v_{28}) + u_{10}(v_{24}+v_{27}-v_{29}) + u_{11}(v_{26}-v_{29})   + u_{14}(v_{24}+v_{25}) + u_{15}(v_{26}+v_{27}) + u^2(2v_{24}+v_{25}-2v_{26}-v_{27}+v_{28}) =0.$

So we have a system with variables $v_i:$
\begin{center}
$\left\{\begin{array}{rl}
-v_{26}- v_{28}&=0\\
v_{24}+v_{27}-v_{29}&=0\\
v_{26}-v_{29}&=0\\
v_{24}+v_{25}&=0\\
v_{26}+v_{27}&=0\\
2v_{24}+v_{25}-2v_{26}-v_{27}+v_{28}&=0\\
\end{array}\right.$
\end{center}
We obtain the solution $(v_{24},v_{25},v_{26},v_{27},v_{28},v_{29})=(2v,-2v,v,-v,-v,v)$ for all $v\in \F_q,$ i.e. $\lambda([h_4,l_1y_2])=1$ for all $l_1y_2\in L_1T_2$ iff $h_4=r_4(v)\in R_4.$

It is clear that $\lambda^{HX_5S_1}(r_4(v))=\lambda^{HX_5S_1}(1)\phi_{B_4}(v)$ for all $r_4(v)\in R_4$ by checking directly that $\lambda(r_4(v))=\phi_{B_4}(v).$

\medskip
(b) Suppose that $B_4\neq 0.$
%By the action of $L_2$ on $\lambda|_{H_4}$ we can pick $x\in L_2$ such that ${}^x\lambda|_{X_{26}}=\phi_{B_4}$ and ${}^x\lambda|_{X_i}=1_{X_i}$ for the other $X_i\subset H_4.$ Therefore, we choose $\lambda$ such that $\lambda|_{X_i}=\phi$ for all $X_i\subset H_6=Z(U),$ $\lambda|_{X_{26}}=\phi_{B_4},$  $\lambda|_{X_i}=1_{X_i}$ for all others $X_i\subset H_5H_4,$ and $\lambda|_{X_i}=\phi_{b_i}$ for the others $X_i\subset H_3H_2$ where $b_i\in \F_q.$
%
Since $Stab_T(\lambda|_{H_6H_5H_4})=S_2=L_2T_3T_4,$ we are going to show that $Stab_{S_2}(\lambda|_{H_6H_5H_4H_3})=T_4,$ and then, $Stab_{T_4}(\lambda)=1$ is done by using the same argument in Lemma \ref{lem:5ext1} (b). It means that we find all $y\in S_2$ such that $\lambda([x,y])=1$ for all $x\in H_3$ since $\lambda([H_6H_5H_4,S_2])=\{1\}.$

 It is clear that $T_4\subset Stab_{S_2}(\lambda|_{H_6H_5H_4H_3}).$ So by (a) and $|H_3|=q^4=|L_2T_3|,$ it suffices to show that $L_2T_3$ acts faithfully on the set of all  extensions of $\lambda|_{H_6H_5H_4}$ to $H_6H_5H_4H_3,$ i.e. $Stab_{L_2T_3}(\lambda|_{H_6H_5H_4H_3})=\{1\}.$
 By the root heights and $H$ is abelian, $[H_3,L_2T_3]=[H_3,T_3][H_3,L_2],$ where $[H_3,T_3]$ is computed in Lemma \ref{lem:5ext1} (b). Since we work with $HX_5S_1/\overline{H_5},$ for each $x=\prod_{i=18}^{21}x_i(v_i)\in H_3$ and $y=\prod_{j=1}^{11}x_j(u_j)\prod_{j=14}^{17}x_j(u_j)x_{22}(u_{22})\in L_2T_3,$ we have
\\
$[x,y]=[x,x_{16}x_{17}x_{22}][x_{18},x_3][[[x_{18},x_3],x_4],x_6][[[x_{18},x_3],x_6],x_7][[x_{18},x_3],x_{14}][x_{18},x_6]\\
\left[[[x_{18},x_6],x_7],x_8\right][[x_{18},x_6],x_{15}][[x_{18},x_6],x_{11}][[x_{18},x_6],x_9][x_{19},x_2][[[x_{19},x_2],x_1],x_4]\\
\left[[[x_{19},x_2],x_1],x_6\right][[[x_{19},x_2],x_4],x_6][[[x_{19},x_2],x_6],x_7][[x_{19},x_2],x_{14}][x_{19},x_1][[[x_{19},x_1],x_6],x_7]\\
\left[[x_{19},x_1],x_{10}\right][[x_{19},x_1],x_{14}][x_{19},x_6][[[x_{19},x_6],x_7],x_8][[x_{19},x_6],x_{10}][[x_{19},x_6],x_{15}][x_{20},x_2]\\
\left[[[x_{20},x_2],x_3],x_4\right][[[x_{20},x_2],x_3],x_7][[[x_{20},x_2],x_7],x_8][[x_{20},x_2],x_{15}][[x_{20},x_2],x_{11}]\\
\left[[x_{20},x_2],x_9\right][x_{20},x_3][[[x_{20},x_3],x_7],x_8][[x_{20},x_3],x_{10}][[x_{20},x_3],x_{15}][x_{20},x_7][[x_{20},x_7],x_9]\\
\left[x_{21},x_4\right][[x_{21},x_4],x_9][x_{21},x_8][[x_{21},x_8],x_{10}][[x_{21},x_8],x_{11}]\\
=x_{37}(v_{18}u_{16})x_{42}(v_{18}u_{22})x_{43}(v_{19}u_{22})x_{40}(-v_{20}u_{17})x_{39}(-v_{21}u_{16})x_{41}(-v_{21}u_{17})\\
x_{25}(2v_{18}u)x_{40}(-2v_{18}u^3)x_{41}(4v_{18}u^3)x_{41}(-2v_{18}u^3)x_{26}(v_{18}u)x_{42}(-4v_{18}u^3)x_{42}(2v_{18}u^3)\\
x_{40}(v_{18}u^3)x_{38}(-v_{18}u^3)x_{25}(-2v_{19}u)x_{37}(-2v_{19}u^3)x_{38}(2v_{19}u^3)x_{40}(2v_{19}u^3)x_{41}(-4v_{19}u^3)\\
x_{41}(2v_{19}u^3)x_{24}(-v_{19}u)x_{39}(-2v_{19}u^3)x_{37}(v_{19}u^3)x_{39}(v_{19}u^3)x_{27}(v_{19}u)x_{43}(-4v_{19}u^3)\\
x_{40}(-v_{19}u^3)x_{43}(2v_{19}u^3)x_{26}(-2v_{20}u)x_{40}(4v_{20}u^3)x_{41}(-8v_{20}u^3)x_{42}(8v_{20}u^3)x_{42}(-4v_{20}u^3)\\
x_{40}(-2v_{20}u^3)x_{38}(2v_{20}u^3)x_{27}(2v_{20}u)x_{43}(-8v_{20}u^3)x_{40}(-2v_{20}u^3)x_{43}(4v_{20}u^3)x_{28}(2v_{20}u)\\
x_{39}(-2v_{20}u^3)x_{28}(-v_{21}u)x_{39}(v_{21}u^3)x_{29}(-2v_{21}u)x_{42}(-2v_{21}u^3)x_{43}(2v_{21}u^3).$

Evaluating with $\lambda$ to get 1, we need the following equation true for all $v_i:$
\\
$v_{18}(u_{16} +u_{22} +2b_{25}u +b_{26}u -2u^3)+ v_{19}(u_{22} -b_{24}u -2b_{25}u +b_{27}u -3u^3) + v_{20}(-u_{17} -2b_{26}u +2b_{27}u +2b_{28}u -3u^3 ) + v_{21}(-u_{16} -u_{17} -b_{28}u -2b_{29}u +u^3 )=0.$

So we have a system with variables $u_j$ and $u:$
\begin{center}
$\left\{\begin{array}{rl}
u_{16} +u_{22} +2b_{25}u +b_{26}u -2u^3 &=0,\\
u_{22} -b_{24}u -2b_{25}u +b_{27}u -3u^3 &=0,\\
-u_{17} -2b_{26}u +2b_{27}u +2b_{28}u -3u^3 &=0,\\
-u_{16} -u_{17} -b_{28}u -2b_{29}u +u^3  &=0.
\end{array}\right.$
\end{center}
It is equivalent to:
\begin{center}
$\left\{\begin{array}{l}
u_{22}=3u^3+(b_{24}+2b_{25}-b_{27})u,\\
u_{17}=2u^3+(3b_{26}+2b_{27}+2b_{28})u,\\
u_{16}=4u^3+(2b_{26}-2b_{27}-3b_{28})u,\\
(2b_{24}-2b_{25}+b_{26}-b_{27}-b_{28}+b_{29})u=0.
\end{array}\right.$
\end{center}

The last equation in the system is actually $B_4u=0.$ Since $B_4\neq 0,$ the only solution of this system is $(u_{16},u_{17},u_{22})=(0,0,0),$ i.e. $Stab_{L_2T_3}(\lambda|_{H_6H_5H_4H_3})=\{1\}.$ Hence,  $Stab_{S_2}(\lambda|_{H_6H_5H_4H_3})=T_4$ and $Stab_{S_1}(\lambda)=\{1\}.$

The above argument also proves that $L_1T_2T_3$ acts transitively on the set of all extensions of $\lambda|_{H_6H_5H_4}$ to $H_6H_5H_4H_3$ with the same $B_4\neq 0.$ The number of these extensions is $|H_4H_3|/|R_4|.$ Therefore, there exists $x\in L_1T_2T_3$ such that ${}^x\lambda|_{X_i}=\phi$ for all $X_i\subset H_6,$ ${}^x\lambda|_{X_{29}}=\phi_{B_4},$  ${}^x\lambda|_{X_i}=1_{X_i}$ for the others $X_i\subset H_5H_4H_3.$ Let $\lambda$ be this  linear character.
By Lemma \ref{lem:Reduction} with $G=HX_5S_1,$ $N=M=HX_5T_4,$ $X=L_1T_2T_3,$ $Z=H_6R_4,$ $Y=H_3\prod_{i=24}^{28}X_i,$ the induction map from $Irr(HX_5T_4/Y,\lambda)$ to $Irr(HX_5S_1,\lambda)$ is bijective. Let $\eta,\eta'$ be two extensions of $\lambda|_{H_6H_5H_4H_3}$ to $HX_5.$ We have $\eta^{HX_5S_2},\eta^{HX_5S_2}\in Irr(HX_5S_2/Y,\lambda).$ Using the same argument in Lemma \ref{lem:5ext1} (b), we obtain $(\eta^{HX_5S_2},\eta'^{HX_5S_2})=1$ iff $\eta|_{R_2}=\eta'|_{R_2}$ and $\eta|_{X_5}=\eta'|_{X_5}.$

\medskip
(c) Suppose that $B_4=0.$ By (a), $S_1/S_2$ acts faithfully on the set of all extensions of $\lambda|_{H_6H_5}$ to $H_6H_5H_4$ with the same $B_4.$  Since $|S_1/S_2|=q^5=|H_4/R_4|,$ this action is transitive. Hence, with $B_4=0,$ there exists $x\in S_1$ such that ${}^x\lambda|_{X_i}=1_{X_i}$ for all $X_i\subset H_5H_4.$ Let $\lambda$ be this linear character. Since $S_2X_5=Stab_{S_1X_5}(\lambda|_{H_6H_5H_4})$ is a transversal of $H$ in $HX_5S_2,$ we have $\lambda^{HX_5S_2}|_{H_4}=[HX_5S_2:H]\lambda|_{H_4}=|X_5S_2|1_{H_4}.$ So $H_4\subset ker(\lambda^{HX_5S_2}).$ By Lemma \ref{lem:Reduction} for $G=HX_5S_1$ with $N=M=HX_5S_2,$ $X=T_2,$ $Y=H_4$ and $Z=H_6,$ the induction map from  $Irr(HX_5S_2/\overline{H_5H_4},\lambda)$ to $Irr(HX_5S_1,\lambda)$ is bijective where $\overline{H_5H_4}$ is the normal closure of $H_5H_4$ in $HX_5S_2.$ Since $H_5H_4\subset ker(\lambda^{HX_5S_2}),$ we have $Irr(HX_5S_2/\overline{H_5H_4},\lambda)=Irr(HX_5S_2,\lambda).~\Box$

\subsubsection{Proof of Lemma \ref{lem:5ext3}}
\label{proof-lem:5ext3}
Recall that $\lambda$ is a linear character of $H$ such that $\lambda|_{X_i}=\phi$ for all $X_i\subset H_6,$ $\lambda|_{X_i}=1_{X_i}$ for all $X_i\subset H_5H_4,$ and $\lambda|_{X_i}=\phi_{b_i}$ for the others $X_i\subset H_3H_2$ where $b_i\in \F_q.$
By Lemma \ref{lem:5ext2} (c), we work with the quotient group $HX_5S_2/\overline{H_5H_4}.$ Abusing the notation of root groups, we continue to call them root groups in the quotient group.

\medskip
(a) By the computation in Lemma \ref{lem:5ext2} (b) with $B_4=0,$ it is clear that $S_3=Stab_{S_2}(\lambda|_{H_6H_5H_4H_3}).$ Now we show that $R_3=\{x\in H_3: |\lambda^{HX_5S_2}(x)|=\lambda^{HX_5S_2}(1)\}.$ Since $X_5S_2$ is a transversal of $H$ in $HX_5S_2,$ we are going to find  $x\in H_3$ such that $\lambda([x,y])=1$ for all $y\in S_2.$ Since $[H_3,X_5]=\{1\}=[H_3,T_4],$ it is enough to work with $x\in H_3$ and $y\in S_2T_3.$ For each $x=\prod_{i=18}^{21}x_i(v_i)\in H_3$ and $y=\prod_{j=1}^{11}x_j(u_j)\prod_{j=14}^{17}x_j(u_j)x_{22}(u_{22})\in S_2T_3,$ by the computation in Lemma \ref{lem:5ext2}, we find $(v_i)_{i\in[18..21]}$ satisfying for all $u_j$ and $u$ in the following equation:
\\
$u_{16}(v_{18}-v_{21} ) +u_{17}(-v_{20}-v_{21} )  +u_{22}(v_{18}+v_{19})   +u^3(-2v_{18} -3v_{19} -3v_{20}+v_{21} )=0.$

We have a system with variables $v_i:$
\begin{center}
$\left\{\begin{array}{rl}
v_{18}-v_{21} &=0,\\
-v_{20}-v_{21}&=0,\\
v_{18}+v_{19}&=0,\\
-2v_{18} -3v_{19} -3v_{20}+v_{21} &=0.
\end{array}\right.$
\end{center}

 Its solutions are $(v_{18},v_{19},v_{20},v_{21})=(u,-u ,-u ,u)$ for all $u\in \F_q,$ i.e. $x=r_3(u)\in R_3.$ Now to show that $\lambda^{HX_5S_2}|_{R_3}=[HX_5S_2:H]\phi_{B_3},$ it is enough to check $\lambda(r_3(t))=\phi_{B_3}(t)$ which is clear.

\medskip
(b) Suppose that $B_3\neq 0.$ By (a) we have $Stab_{S_2}(\lambda|_{H_6H_5H_4H_3})=S_3=L_3T_4.$
%By the action of $S_2$ on $\lambda|_{H_3},$ we can choose $\lambda$ such that $\lambda|_{X_{21}}=\phi_{-B_3}$ and $\lambda|_{X_i}=1_{X_i}$ for the others $X_i\subset H_5H_4H_3.$
To show that $Stab_{S_2}(\lambda)=\{1\},$ since $|L_3T_4|=q^2=|H_2|,$ we show that $L_3T_4$ acts faithfully on the set of all extensions of $\lambda|_{H_6H_5H_4H_3}$ to $H,$ i.e. proving that there is no nontrivial $y\in L_3T_4$ such that $\lambda([x,y])=1$ for all $x\in H_2.$

By the root heights, $[H_2,L_3T_4]=[H_2,T_4][H_2,L_3],$ where $[H_2,T_4]$ is computed in Lemma \ref{lem:5ext1}(b). For $x=x_{12}(v_{12})x_{13}(v_{13})\in H_2$ and $y=l_3(u)x_{23}(u_{23})\in L_3T_4,$ we have
\\
$[x,y]=[x,x_{23}][x,l_3]= [x_{12},x_{23}][x_{12},x_2]
 [[[[x_{12},x_2],x_3],x_4],x_6][[[[x_{12},x_2],x_3],x_6],x_7]\\
\left[[[[x_{12},x_2],x_6],x_7],x_8\right][x_{12},x_3]
 [[[[x_{12},x_3],x_6],x_7],x_8][x_{12},x_6] [[[x_{12},x_2],x_3],x_{14}]\\
\left[[[x_{12},x_2],x_6],x_9\right] [[[x_{12},x_2],x_6],x_{11}][[[x_{12},x_2],x_6],x_{15}]
 [[[x_{12},x_3],x_6],x_{10}][[[x_{12},x_3],x_6],x_{15}]
 \\
\left[[[x_{12},x_6],x_7],x_9\right][[x_{12},x_2],x_{22}] [[x_{12},x_2],x_{16}]
 [[x_{12},x_3],x_{22}][[x_{12},x_6],x_{17}] [[x_{12},x_9],x_{10}]\\
\left[[x_{12},x_9],x_{14}\right][x_{13},x_{23}][x_{13},x_4][x_{13},x_7]
 [[[x_{13},x_4],x_7],x_9][[[x_{12},x_7],x_8],x_{10}] \\
\left[[[x_{13},x_7],x_8],x_{11}\right][[x_{13},x_4],x_{17}]
 [[x_{13},x_7],x_{16}][[x_{13},x_7],x_{17}][[x_{13},x_{10}],x_{15}]\\
\left[[x_{13},x_{10}],x_{11}\right][[x_{13},x_{11}],x_{15}]\\
 =x_{37}(-v_{12}u_{23})x_{18}(-2v_{12}u)   x_{40}(4v_{12}u^4)x_{41}(-8v_{12}u^4)
 x_{42}(8v_{12}u^4)x_{19}(2v_{12}u)   x_{43}(-8v_{12}u^4)\\
 x_{20}(v_{12}u)x_{41}(4v_{12}u^4)x_{38}(2v_{12}u^4) x_{40}(-2v_{12}u^4)x_{42}(-4v_{12}u^4)
 x_{40}(-2v_{12}u^4)x_{43}(4v_{12}u^4)\\
 x_{39}(-2v_{12}u^4) x_{42}(-6v_{12}u^4)x_{37}(-8v_{12}u^4) x_{43}(6v_{12}u^4)x_{40}(-2v_{12}u^4)
  x_{37}(v_{12}u^4)x_{39}(v_{12}u^4)  \\
 x_{38}(-v_{13}u_{23}) x_{20}(-v_{13}u)x_{21}(2v_{13}u) x_{39}(2v_{13}u^4)x_{42}(-4v_{13}u^4)
 x_{43}(4v_{13}u^4) x_{40}(2v_{13}u^4) \\
x_{39}(-8v_{13}u^4) x_{41}(-4v_{13}u^4) x_{42}(2v_{13}u^4) x_{40}(v_{13}u^4)x_{43}(-2v_{13}u^4).$

Evaluating with $\lambda$ to get $1,$ we obtain in the following equation:
\\
$v_{12}(-s_{23} -2b_{18}u +b_{20}u +2b_{19}u  -2u^4) +v_{13}(-u_{23} -b_{20}u +2b_{21}u -2u^4)=0.$

We have a system with variables $u_j$ and $u:$
\begin{center}
$\left\{\begin{array}{rl}
-u_{23} -2b_{18}u +b_{20}u +2b_{19}u -2u^4&=0,\\
-u_{23} -b_{20}u +2b_{21}u -2u^4&=0.
\end{array}\right.$
\end{center}
It is equivalent to:
\begin{center}
$\left\{\begin{array}{l}
u_{23}=3u^4+ (-2b_{18} +b_{20} +2b_{19})u,\\
 (b_{18} -b_{20} -b_{19} +b_{21})u=0.
\end{array}\right.$
\end{center}

The last equation is actually $B_3u=0.$ Since $B_3\neq 0,$ the only solution is $(u_{23},u)=(0,0),$ i.e. $Stab_{L_3T_4}(\lambda)=\{1\}$ or
$L_3T_4$ acts faithfully on the set of all extensions of $\lambda|_{H_6H_5H_4H_3}$ to $H.$ Hence, we also get $Stab_{S_2}(\lambda)=\{1\}.$

%%Re-check condition of Reduction Lemma in this case

Therefore, there exists $x\in L_3T_4$ such that ${}^x\lambda|_{X_i}=\phi$ for all $Xi\subset H_6,$ ${}^x\lambda|_{X_{21}}=\phi_{B_3},$ ${}^x\lambda|_{X_i}=1_{X_i}$ for the others $X_i\subset H_5H_4H_3.$ Let $\lambda$ be this linear. By Lemma \ref{lem:Reduction} with $G=HX_5S_2,$ $N=M=HX_5,$ $X=S_2,$ $Y=\prod_{i=18}^{20}X_i$ and $Z=H_6X_{21},$ the induction map from $Irr(HX_5/Y,\lambda)$ to $Irr(HX_5S_2,\lambda)$ is bijective. Using the same technique in Lemma \ref{lem:5ext1} (c), the rest statement holds.

\medskip
(c) Suppose $B_3=0.$ By (a), $S_2/S_3$ acts faithfully on the set of all extensions of $\lambda|_{H_6H_5H_4}$ to $H_6H_5H_4H_3$ with the same $B_3.$ Since $|S_2/S_3|=q^3=|H_3/R_3|,$ this action is transitive. Hence, there exists $x\in S_2$ such that ${}^x\lambda|_{X_i}=1_{X_i}$ for all $X_i\subset H_5H_4H_3.$ Let $\lambda$ be this linear.
Since $X_5S_3$ is a transversal of $H$ in $HX_5S_3$ and $S_3=Stab_{S_2}(\lambda|_{H_6H_5H_4H_3}),$ we have $\lambda^{HX_5S_3}|_{H_5H_4H_3}=\lambda^{HX_5S_3}(1)\lambda|_{H_5H_4H_3}.$ Therefore, $H_5H_4H_3\subset ker(\lambda^{HX_5S_3}),$ so is its normal closure $\overline{H_5H_4H_3}$ in $HX_5S_3.$

By Lemma \ref{lem:Reduction} with $G=HX_5S_2,$ $N=M=HX_5S_3,$ $X=T_3,$ $Y=H_3$ $Z=H_6,$ the induction map from $Irr(HX_5S_3/Y,\lambda)$ to $Irr(HX_5S_2,\lambda).$ Since $Y\subset ker(\lambda^{HX_5S_3}),$ we have $Irr(HX_5S_3/Y,\lambda)=Irr(HX_5S_3,\lambda).$

\subsubsection{Proof of Lemma \ref{lem:5ext4}}
\label{proof-lem:5ext4}
Recall that $\lambda$ is a linear character of $H$ such that $\lambda|_{X_i}=\phi$ for all $X_i\subset H_6=Z(U),$ $\lambda|_{X_i}=1_{X_i}$ for all $X_i\subset H_5H_4H_3,$ and $\lambda|_{X_i}=\phi_{b_i}$ for the others $X_i\subset H_2$ where $b_i\in \F_q.$
By Lemma \ref{lem:5ext2} (c), we work with the quotient group $HX_5S_3/\overline{H_5H_4H_3}.$ Abusing the notation of root groups, we continue to call them root group in the quotient group.

\medskip
(a) By the computation in Lemma \ref{lem:5ext3} (b) with $B_3=0,$ $S_4=Stab_{S_3}(\lambda).$ Now we show that $R_2=\{x\in H_2: |\lambda^{HX_5S_3}(x)|=\lambda^{HX_5S_3}(1)\}.$ Since $X_5S_3$ is a transversal of $H$ in $HX_5S_3,$ we are going to find  $x\in H_2$ such that $\lambda([x,y])=1$ for all $y\in S_3.$ Since $[H_2,X_5]=\{1\},$ it is enough to work with $x\in H_2$ and $y\in S_4.$ For each $x=\prod_{i=12}^{13}x_i(v_i)\in H_2$ and $y=l_3(u)x_{23}(u_{23})\in S_3T_4,$ by the computation in Lemma \ref{lem:5ext3} (b), we find $(v_{12},v_{13})$ satisfying for all $u_{23}$ and $u$ in the following equation:
\[u_{23}(-v_{12}-v_{13})+2u^4(-v_{12}-v_{13})=0.\]
So $(v_{12},v_{13})=(v,-v)$ for all $v\in\F_q,$ i.e. $x=r_2(v).$ %So  $R_2=\{x\in H_2: |\lambda^{HX_5S_3}(x)|=\lambda^{HX_5S_3}(1)\}.$
Since $\lambda(r_2(v))=\phi_{B_2}(v)$ for all $r_2(v)\in R_2,$ we have $\lambda^{HX_5S_3}|_{H_2}=[HX_5S_3:H]\phi_{B_2}.$

\medskip
(b) Suppose that $B_2\in \F_q-\{c^4:c\in \F_q^\times\}.$ Let $\eta$ be an extension of $\lambda$ to $HX_5.$ Since $S_4=Stab_{S_3}(\lambda),$ to get $I_{HX_5S_3}(\eta)=HX_5,$ we show that $S_4$ acts transitively on the set of all extensions of $\lambda$ to $HX_5.$ Hence, we find all $l_4\in S_4$ such that $\lambda([hx_5,l_4])=1$ for all $h\in H$ and $x_5\in X_5.$ Since $S_4=Stab_{S_3}(\lambda),$ we have $\lambda([h,l_4])=1$ for all $h\in H, l_4\in S_4.$ Thus we compute $[x_5,l_4].$ Since we work with $HX_5S_3/\overline{H_5H_4H_3},$ for each $x_5(v_5)\in X_5$ and $l_4(u)\in S_4,$ we have
\\
$\left[x_5(v_5),l_4(u)\right]=[x_5,x_4][[[x_5,x_4],x_9],x_{10}][[[x_5,x_4],x_9],x_{14}][[[[x_5,x_4],x_6],x_7],x_9]\\
\left[[[x_5,x_4],x_6],x_{17}\right][[x_5,x_4],x_{23}][x_5,x_6][[[x_5,x_6],x_{10}],x_{11}][[[x_5,x_6],x_{10}],x_{15}]\\
\left[[[x_5,x_6],x_{11}],x_{15}\right][[[[x_5,x_6],x_7],x_8],x_{10}][[[[x_5,x_6],x_7],x_8],x_{11}][[[x_5,x_6],x_7],x_{16}]\\
\left[[[x_5,x_6],x_7],x_{17}\right][[x_5,x_6],x_{23}][[x_5,x_{14}],x_{16}][[x_5,x_{14}],x_{17}][[x_5,x_{11}],x_{22}]\\
\left[[x_5,x_{10}],x_{16}\right][[x_5,x_{10}],x_{22}]\\
=x_{12}(-v_5u)x_{37}(-v_5u^5)x_{39}(-v_5u^5)x_{39}(2v_5u^5)x_{40}(2v_5u^5)x_{37}(3v_5u^5)x_{13}(v_5u)x_{40}(v_5u^5)\\
x_{42}(2v_5u^5)x_{43}(-2v_5u^5)x_{42}(-4v_5u^5)x_{43}(4v_5u^5)x_{39}(-8v_5u^5)x_{41}(-4v_5u^5)x_{38}(-3v_5u^5)\\
x_{39}(4v_5u^5)x_{41}(2v_5u^5)x_{43}(-3v_5u^5)x_{42}(3v_5u^5)x_{37}(4v_5u^5)
$

Evaluating with $\lambda$  to get 1, for all $v_5$ we need
\[v_5(-(b_{12}-  b_{13})u +u^5 )\in ker(\phi),\]
which is $v_5(u^5-B_2u)\in ker(\phi)$ for all $v_5.$ Hence, we solve for $u:$ $u(u^4-B_2)=0.$ Since $B_2\in \F_q-\{c^4:c\in\F_q^\times\},$ this equation only has one trivial solution $u=0,$ i.e. $Stab_{S_4}(\eta)=\{1\},$ or $I_{HX_5S_3}(\eta)=HX_5.$

\medskip
(c) Suppose that $B_2=c^4\in \F_q^\times.$ Let $\eta$ be an extension of $\lambda$ to $HX_5.$ Continue the computation in (b),  the equation $u(u^4-B_2)=0$ has $5$ solutions $u\in\{ac:a\in\F_5\},$ i.e. $l_4(u)\in F_4.$ Hence, $I_{HX_5S_3}(\eta)=HX_5F_4.$ Since $[HX_5,F_4]\subset ker(\eta),$ $\eta$ extends to $HX_5F_4,$ i.e. $\lambda$ extends to $HX_5F_4.$

Since $S_4=Stab_{S_3}(\lambda)\cong \F_q,$ we have $[H,S_4]\subset ker(\lambda).$ So $\lambda$ extends to $HS_4\unlhd HX_5S_3.$ Let $\lambda'$ be an extension of $\lambda$ to $HS_4.$ We find $I_{HX_5S_3}(\lambda').$ Since $Stab_{X_5S_3}(\lambda')\subset Stab_{X_5S_3}(\lambda'|_H)=X_5S_4,$ it is enough to find all $x_5\in X_5$ such that $\lambda'([x_5,hl_4])=1$ for all $hl_4\in HS_4.$ Since $HX_5$ is abelian, we have $[x_5,hl_4]=[x_5,l_4].$ For each $x_5(v_5)\in X_5$ and $l_4(u)\in S_4$, by the computation in (b), we need
\begin{center}
$v_5(u^5-B_2u)\in ker(\phi),$ for all $u\in \F_q.$
\end{center}
By Proposition \ref{prop:mainFp}, there are 5 solutions $v_5\in \{ac_\phi:a\in\F_5\},$ i.e. $x_5(v_5)\in F_5.$ Hence, $I_{HX_5S_3}(\lambda')=HF_5S_4.$
Since $[F_5,S_4]\subset ker(\lambda'),$ $\lambda'$ extends to $HF_5S_4,$ i.e. $\lambda$ extends to $HF_5S_4.$

Let $\lambda_1,\lambda_2$ be two extensions of $\lambda$ to $HX_5F_4,$ and  $\gamma$ an extension of $\lambda$ to $HF_5S_3.$ Since the degree of all irreducible constituents of $\lambda^{HX_5S_3}$ is $\frac{q^2}{5},$ we have ${\lambda_1}^{HX_5S_3},$ ${\lambda_2}^{HX_5S_3},$ $\gamma^{HX_5S_3}\in Irr(HX_5S_3,\lambda).$

Choose $1\in S\subset S_3$ as a representative set of the double coset $HF_5S_4\backslash HX_5S_3/HX_5F_4.$ Since $HF_5S_4\cap HX_5F_4=HF_5F_4$ and $HX_5F_4\unlhd HX_5S_3,$ by Mackey formula,
\begin{center}
$\begin{array}{ll}
({\lambda_1}^{HX_5S_3},\gamma^{HX_5S_3})&=\sum_{s\in S}({}^s\lambda_1|_{{}^s(HX_5F_4)\cap HF_5S_4},\gamma|_{{}^s(HX_5F_4)\cap HF_5S_4})
\\
&=\sum_{s\in S}({}^s\lambda_1|_{HF_5F_4},\gamma|_{HF_5F_4})
\end{array}$
\end{center}

For each $s\in S,$ if ${}^s\lambda_1|_{HF_5F_4}=\gamma|_{HF_5F_4},$ then ${}^s\lambda_1|_H=\gamma|_H.$ Since both are extensions of $\lambda,$ we have ${}^s\lambda=\lambda,$ i.e. $s\in Stab_{S_3}(\lambda)=S_4.$ There is unique $1\in S\cap S_4$ since $S$ is a representative set of $HF_5S_4\backslash HX_5S_3/HX_5F_4.$ So $({\lambda_1}^{HX_5S_3},\gamma^{HX_5S_3})=(\lambda_1|_{HF_5F_4},\gamma|_{HF_5F_4})=1$ iff $\lambda_1|_{F_i}=\gamma|_{F_i},i\in \{4,5\}.$

Therefore, ${\lambda_1}^{HX_5S_3}=\gamma^{HX_5S_3}={\lambda_2}^{HX_5S_3}$ iff $\lambda_1|_{F_i}=\lambda_2|_{F_i},i\in \{4,5\}.~\Box$

\section*{Acknowledgement} The first author would like to express his heartfelt thanks to Professor Geoffrey Robinson for all of his
support during the preparation of this work.

\end{document}